\documentclass[leqno,11pt]{article}

\usepackage{amssymb}
\usepackage{euscript}
\usepackage[dvips]{graphicx}
\usepackage{flafter}
\usepackage{pstricks}
\usepackage[latin1]{inputenc}
\usepackage{amsmath}
\usepackage{amsfonts}
\usepackage{pstricks-add}
\usepackage{dsfont}
\usepackage{bm}
\usepackage{variations}
\usepackage{hyperref}

\usepackage[refpage]{nomencl}
\makenomenclature
\RequirePackage{ifthen}
 \renewcommand{\nomgroup}[1]{%
 \ifthenelse{\equal{#1}{E}}{\item[\textbf{The log-correlated Gaussian field}]}{%
 \ifthenelse{\equal{#1}{A}}{\item[\textbf{Processes}]}{}}{%
 \ifthenelse{\equal{#1}{B}}{\item[\textbf{Set of functions}]}{}}{%
 \ifthenelse{\equal{#1}{C}}{\item[\textbf{Elements or subsets of $[0,R]^\d$}]}{}}{%
 \ifthenelse{\equal{#1}{D}}{\item[\textbf{Constants}]}{}}{%
 \ifthenelse{\equal{#1}{F}}{\item[\textbf{The Lebesgue measure, intervals}]}{}}{%
 \ifthenelse{\equal{#1}{G}}{\item[\textbf{Modulus of continuity}]}{}}{%
 \ifthenelse{\equal{#1}{H}}{\item[\textbf{Events}]}{}}{%
 \ifthenelse{\equal{#1}{I}}{\item[\textbf{Various}]}{}}
 }
\nomgroup

\usepackage{makeidx}
\makeindex

\usepackage{mathrsfs} 

\evensidemargin\oddsidemargin

\begin{document}
 
\newcommand{\eqnsection}{
\renewcommand{\theequation}{\thesection.\arabic{equation}}
   \makeatletter
   \csname  @addtoreset\endcsname{equation}{section}
   \makeatother}
\eqnsection

\def\r{{\mathbb R}}
\def\e{{\mathbb E}}
\def\p{{\mathbb P}}
\def\P{{\bf P}}
\def\E{{\bf E}}
\def\Q{{\bf Q}}
\def\z{{\mathbb Z}}
\def\N{{\mathbb N}}
\def\T{{\mathbb T}}
\def\G{{\mathbb G}}
\def\L{{\mathbb L}}
\def\1{{\mathds{1}}}
\def\deg{\chi}

\def\ee{\mathrm{e}}
\def\d{\, \mathrm{d}}
\def\S{\mathscr{S}}
\def\bs{{\tt bs}}
\def\bbeta{{\bm \beta}}

\def\d{\mathtt{d}}
\def\ttheta{{\bm \theta}}
\def\t{{\bf{t}}}
\def\a{{\bf{a}}}
\def\deg{\chi}
\def\B{\mathfrak{B}}

\def\M{{\mathbb{M}}}

\def\mS{{\underline{S}}}
\def\MS{{\overline{S}}}

\def\Y{{\bf Y}}

\def\R{{\bf R}}
\def\o{\mathfrak{O}}
\def\Bt{ \blacktriangleright}
\def\Bti{\blacktriangleleft}
\def\Ut{\Delta} 
\def\BB{ \blacktriangle}
\def\CC{\ntriangleleft}
\def\BBi{\underline{\blacktriangle}}
\def\CCi{\ntrianglelefteq}

\def\DDelta{\underline{\Delta} }

\def\DD{ \vartriangleright}
\def\DDi{\trianglerighteq}

\def\Loz{ \blacklozenge}
\def\Squ{\blacksquare}

\def\Loza{ \underline{\blacklozenge}}
\def\Squa{\underline{\blacksquare}}

\def\Tr{\blacktriangledown}
\def\Tra{\underline{\blacktriangledown}}

\newtheorem{theorem}{Theorem}[section]

\newtheorem{Definition}[theorem]{Definition}
\newtheorem{Lemma}[theorem]{Lemma}
\newtheorem{Proposition}[theorem]{Proposition}
\newtheorem{Remark}[theorem]{Remark}
\newtheorem{corollary}[theorem]{Corollary}


\title{ Maximum of a log-correlated Gaussian field}
\author{Thomas Madaule, Universit\'e Paris XIII }
\maketitle

\vglue50pt

{\leftskip=2truecm \rightskip=2truecm \baselineskip=15pt \small

\noindent{\slshape\bfseries Abstract.}  We study the maximum of a Gaussian field on $[0,1]^\d$ ($\d \geq 1$) whose correlations decay logarithmically with the distance.  Kahane \cite{Kah85} introduced this model to construct mathematically the Gaussian multiplicative chaos in the subcritical case. Duplantier, Rhodes, Sheffield and Vargas \cite{DRSV12a} \cite{DRSV12b} extended Kahane's construction to the critical case and established the KPZ formula at criticality. Moreover, they made in \cite{DRSV12a} several conjectures on the supercritical case and on the maximum of this Gaussian field. In this paper we resolve Conjecture 12 in \cite{DRSV12a}: we establish the convergence in law of the maximum and show that the limit law is the Gumbel distribution convoluted  by the limit of the derivative martingale.} 

\bigskip
\bigskip

\section{Introduction}
We study the maximum of a Gaussian field on $[0,1]^\d$ ($\d \geq 1$) whose correlations decay logarithmically with the distance. This model was introduced by Kahane \cite{Kah85} to construct mathematically the Gaussian multiplicative chaos (GMC). This family of random fields has found many applications in various fields of science, especially in turbulence and in mathematical finance.

A series of work of Duplantier, Rhodes, Sheffield and Vargas has generated a renewed interest on this model. In \cite{DRSV12a} and \cite{DRSV12b} they extend Kahane's \cite{Kah85} construction of the Gaussian multiplicative chaos to the critical case and establish the KPZ formula at criticality. Their proofs are inspired by the latest advances in the study of the branching random walk (BRW) especially concerning the Seneta-Heyde norming for the additive martingale. Moreover they make several conjectures on the supercritical case and on the maximum of the log-correlated Gaussian field (see \cite{DRSV12a}).

In this paper we resolve the Conjecture 12 in \cite{DRSV12a}: we establish the convergence in law of the maximum and show that the limit law is the Gumbel distribution convoluted  by the limit of the derivative martingale. Moreover we believe that this result could lead to the resolution of the conjecture 11 \cite{DRSV12a} on the existence of the GMC in the supercritical case. Our proof is deeply inspired by a powerful method of Elie Aïdékon, developed in \cite{Aid11}, to show the convergence in law of the minimum of a real-valued branching random walk.

We treat the case of star scale invariant log-correlated fields. This is a general class of field with no restriction on the dimension. It generalizes the notion of branching structure in a continuous setting and may to prove the existence and the uniqueness of the lognormal $*-$scale invariant random measures, see \cite{ARV11}.

Let us mention that in the discrete setting, $\mathbb{Z}^2\cap [0,N]^2$, if we add the zero boundary condition, the model becomes the so-called Gaussian free field (GFF), which has attracted many recent attentions, see \cite{BDG01}, \cite{BDZ11} and \cite{BZe10}. In particular we mention \cite{BDZ13} where they proved the convergence in law of the maximum of GFF after a suitable normalization.

In the first sub-section we shall introduce the model of log-correlated Gaussian random field and state the main result of the paper. In the second sub-section we set out the strategy of the proof.

\subsection{Star scale invariants kernels}
We follow \cite{DRSV12a} to introduce the log-correlated Gaussian field that we will study throughout the paper. We consider a family of centered stationary  Gaussian processes $(X_s(x))_{s\geq 0,\, x\in \r^{\d}}$ $\d\geq 1$, with covariances \nomenclature[e1]{$(X_s(x))_{s\geq 0,\, x\in \r^{\d}}$}{: the log-correlated Gaussian field} \nomenclature[e2]{$\mathtt{k}$}{: the kernel function of $X_\cdot(\cdot)$}\nomenclature[e3]{$\mathtt{g}(\cdot)$}{$:=1-\mathtt{k}(\cdot)$}
\begin{equation}
\label{defX}
\E[X_t(0)X_t(x)]=\int_0^t\mathtt{k}(\ee^u x)du,\quad \forall t>0,\, x\in \r^\d.
\end{equation}
The kernel function $\mathtt{k}:\r^\d\to \r$ is $\mathcal{C}^1$, satisfying $\mathtt{k}(0)=1$ and $\mathtt{k}(x)=0$ if {\bf $x\notin B(0,1):=\{x: |x|\leq 1\}$} ($|x|:= \underset{i\in [1,\d]}{\max}|x_i|$). Such fields have been studied in \cite{ARV11} via a white noise decomposition. We also denote $\mathtt{g}(\cdot):=1- \mathtt{k}(\cdot)$ and introduce for any $t>0$,\nomenclature[e4]{$Y_t(x)$}{$:=X_t(x) -\sqrt{2\d} t $}\nomenclature[c1]{$B(0,1)$}{$ :=\{x:\, \mid x \mid \leq 1\} $}
\begin{equation}
\label{defY}
Y_t(x):=X_t(x)-\sqrt{2\d}t.
\end{equation}
For any $A\subset \r^\d$ bounded, we are interested in 
\begin{equation}
\mathbb{M}_t(A):= \underset{x\in A}{\sup}\,Y_t(x),\qquad t>0,
\end{equation}
the maximum of the Gaussian field on the domain $A$, at time $t$. Let $\mathcal{B}(\r^\d)$ the Borel on $\r^\d$, and $\mathcal{B}_b(\r^\d) $ its restriction to the bounded sets. \nomenclature[e5]{$\M_t([0,1]^\d)$}{$:=\underset{x\in [0,1]^\d}{\sup}\,Y_t(x)$} We introduce for $t>0$ and $\gamma>0$, the random measures $M_t'(dx)$ and $M_t^\gamma(dx)$ defined by:\nomenclature[e6]{$M'_t(A)$}{$:= \int_{A}(-Y_t(x))\ee^{\sqrt{2\d}Y_t(x)+\d t}dx$}\nomenclature[e7]{$M^\gamma_t(A) $}{$:= \int_{A} \ee^{\sqrt{2\d}Y_t(x)+\d t}dx$}
\begin{equation}
\label{1.4} M_t'(A):=\int_{A}(-Y_t(x))\ee^{\sqrt{2\d}Y_t(x)+\d t}dx, \, \, \, M^\gamma_t(A):=\int_{A} \ee^{\gamma Y_t(x)+\gamma \sqrt{2\d} t-\frac{\gamma^2}{2}t}dx, \,\, \forall A\in \mathcal{B}_b(\r^d).
\end{equation}
Kahane in \cite{Kah85} proved that for any $\gamma\in [0,\sqrt{2\d})$ (called subcritical case), there exists a random measure $M_\infty^\gamma$ such that
\begin{equation}
\label{Kah1}
M_t^\gamma(A)\overset{a.s}{\to} M_\infty^\gamma(A),\qquad \forall A\in \mathcal{B}_b(\r^\d),
\end{equation}
whereas for $\gamma\geq \sqrt{2\d}$ (called critical and supercritical case),
\begin{equation}
\label{Kah2}
M_t^\gamma(A)\overset{a.s}{\to} 0,\qquad \forall A\in \mathcal{B}_b(\r^\d).
\end{equation}
One motivation of (\ref{Kah1}) is to give a rigorous construction of  a standard {\it Gaussian multiplicative chaos} (GMC) in the subcritical case which is formally defined as (see \cite{DRSV12a}) a random measure such that for any set $A\in \mathcal{B}_b(\r^d)$,
\begin{equation}
M_\infty^\gamma(A)=\int_A\ee^{\gamma X-\frac{\gamma^2}{2}\E[X^2(x)]}dx,
\end{equation}
where $X$ is a centered log-correlated Gaussian field:
\begin{equation}
\E[X(x)X(y)]=\log_+\frac{1}{|x-y|}+ g(x,y),
\end{equation}
with $\log_+(x)=\max(\log x,0)$ and $g$ a continuous bounded function on $\r^\d\times \r^\d $. It is an important problem to extend the construction for $\gamma\geq \sqrt{2\d}$. In \cite{DRSV12a} the authors are able to construct the GMC in the critical case $\gamma=\sqrt{2\d}$, via the following theorem:
\\

 
{\it 
\noindent{\bf Theorem A (\cite{DRSV12a})}
For each bounded open set $A\subset \r^d$, the martingale $(M_t'(A))_{t\geq 0}$ converges almost surely towards a positive random variable denoted by $M'(A)$.}
 \\

Concerning the construction of the GMC in the case $\gamma>\sqrt{2\d}$, they (\cite{DRSV12a}) conjectured
\\

{\it 
\noindent{\bf Conjectures (\cite{DRSV12a}) }
\begin{eqnarray*}
\nonumber
&&(A)\qquad\qquad\qquad\qquad\qquad \quad t^\frac{3\gamma}{2\sqrt{2\d}} \ee^{t(\frac{\gamma}{\sqrt{2}}-\sqrt{\d})^2} M_t^\gamma(dx) \overset{law}{\to} c_\gamma N_{\frac{\sqrt{2\d}}{\gamma}},\quad \text{as  } t\to \infty, \qquad\qquad\qquad\quad
\end{eqnarray*}
with $c_\gamma$ a positive constant, and $N_\frac{\sqrt{2\d}}{\gamma}$ a known positive random measure.
\begin{equation}
\nonumber
(B) \qquad\qquad\quad\quad \underset{x\in[0,1]^\d}{\sup} Y_{t}(x) +\frac{3}{2\sqrt{2\d}} \log t\overset{law}{\to} G_\d,\qquad \text{as }t\to \infty,\qquad\qquad \qquad\qquad
\end{equation}
where the distribution of $G_\d $ is a Gumbel distribution convoluted with $M_\infty'([0,1]^\d)$.

}

The authors also explained how to obtain the Conjecture (B) from Conjecture (A). Here we do not study Conjecture (A), but we resolve directly Conjecture (B):
\begin{theorem}
\label{maintheorem}
There exists a constant $C^*\in (0,\infty)$ such that, for any  real $z$,
\begin{equation}
\label{eqmaintheorem}
\underset{t\to\infty}{\lim}\P\left( \M_t([0,1]^\d)\leq -\frac{3}{2\sqrt{2d}}\log t -z\right)=\E\left(\ee^{-C^*\ee^{\sqrt{2\d}z} M'([0,1]^\d)}\right).
\end{equation}
\end{theorem}
\paragraph{Remark:}  In this paper we have assumed that the kernel $\mathtt{k}$ has compact support. This hypothesis is essential for the section 3. However it would be possible to relax this hypothesis when the long-range correlations decrease rapidly. 
\\

We believe that this result and the methods developed here, could lead to establish Conjecture (A): Basically, when $\gamma> \sqrt{2\d}$, $M_t^\gamma$ concentrates its mass only on the particles close to the maximum $ \mathbb{M}_t([0,1]^\d) $, where here and in the sequel, by particle in the log-correlated Gaussian field we mean a point $x\in [0,1]^\d $. We expect to establish the convergence of the random measure formed by the particles near to the maximum, just like in the BRW case (see \cite{Aid11}, \cite{Mad11}, \cite{BRV12} and \cite{DRSV12a} for an explicit connection between branching random walk and this model). This direction will be explored elsewhere.  

As in the case of the branching Brownian motion, see \cite{ABK10}, \cite{ABK10a}, \cite{ABK11} and \cite{ABBS11}, our work could also lead to the ``genealogy of the extremal particles" which in our context corresponds to their spatial position. Indeed in Lemma \ref{eqhtoui} we use our understanding of the paths of the extremal particles to prove that they are concentrated in clusters.

As we mentioned before, recently in \cite{BDZ13} the authors  showed the convergence in law of the maximum of the GFF. Furthermore it is believed that there exists some universality between all the log-correlated Gaussian fields, see \cite{CLeD01}. For instance, it is interesting to extend our result to some kernels $\mathtt{k}$ which are not invariant by translation. 

%
%
%
%
%
%

\subsection{Strategy of proof}
Here we try to give a guiding thread for the proof of Theorem \ref{maintheorem}. We mention that this strategy of proof is similar to that used in \cite{Aid11} for the BRW and also in \cite{BDZ13} for the GFF.

%

We start by introducing some notations. It will be convenient to consider a log-correlated Gaussian field starting from an arbitrary $a\in \r$, whose the law is denoted by $\P_a$. The law of $(Y_s(x))_{s\geq 0,x\in \r^\d}$ under $\P_a$ is the same as the law of $(a+Y_s(x))_{s\geq 0,x\in \r^\d}$ under $\P$. For any $l>0$ we define
\begin{equation}
\label{increment}( Y_s^{(l)}(x))_{s\geq 0,\, x\in \r^\d}:= ( Y_{s+l}(x)-Y_l(x))_{s\geq 0,\, x\in \r^\d}.
\end{equation}
This process is independent of $( Y_s(x))_{s\leq l,\, x\in \r^\d}  $ and has the same law as $( Y_s(x\ee^l))_{s\geq 0,\, x\in \r^\d}$ under $\P$, as we will see in (\ref{scalingl}). Observe that for any $x\in \r^\d$, $(Y_s(x))_{s\geq 0}\overset{\text{law}}{=}(B_s-\sqrt{2\d}s)_{s\geq 0}$ with $(B_s)$ a standard Brownian motion. For any $a,\, b,\, l\in \r^+$, define  
\begin{eqnarray*} 
\mathcal{C}_R(l,a,b):=\Big\{ f:\, \underset{x,y\in [0,R]^\d,\, |x-y|\leq \frac{1}{l}}{\sup} \frac{|f(x)-f(y)|}{|x-y|^\frac{1}{3}}\leq 1 ,\,\,\, \underset{y\in[0,R]^\d}{\min}f(y)> a ,\text{  and  }   \underset{y\in[0,R]^\d}{\max}f(y)< b\Big\}.
 \end{eqnarray*}
For any process $(f_s)_{s\geq 0}$, $t_2>t_1\geq 0 $, let
\begin{eqnarray*}
\underline{f}_{t_1}:=\underset{s\leq t_1}{\inf} f_s,&&\qquad \underline{f}_{[t_1,t_2]}:=\underset{t_1\leq s\leq t_2}{\inf} f_s,\qquad \text{and }
\\
\overline{f}_{t_1}:=\underset{s\leq t_1}{\sup} f_s,&&\qquad \overline{f}_{[t_1,t_2]}:=\underset{t_1\leq s\leq t_2}{\sup} f_s.
\end{eqnarray*}
Similarly we also define
\begin{eqnarray*}
|f|_{t_1}:=\underset{s\leq t_1}{\sup} |f_s|&&\text{   and    } \quad |f|_{[t_1,t_2]}:=\underset{t_1\leq s\leq t_2}{\sup} |f_s|.
\end{eqnarray*}
 
 As shown in \cite{DRSV12a}, the typical order of $\M_t$ is $-\frac{3}{2\sqrt{2d}}\log t$, so it will be convenient to introduce
\begin{eqnarray}
\label{uUu} \kappa_\d=\frac{1}{4\sqrt{2\d}},\, &a_t:= -\frac{3}{2\sqrt{2d}}\log t& \text{and }  I_t(z):= [a_t+z-1,a_t+z],\quad z\geq 0.
\end{eqnarray}
\nomenclature[d1]{$ \kappa_\d$}{$:=\frac{1}{8\sqrt{2\d}} $}
\nomenclature[f1]{$a_t$}{$:=-\frac{3}{2\sqrt{2d}}\log t $}
\nomenclature[f2]{$I_t(z)$}{$:= [a_t+z-1,a_t+z],\quad z\in \r $}
\nomenclature[b]{$\mathcal{C}_R(a,b,c)$}{: a subset of $\mathcal{C}([0,R]^\d,\r)$}

For any $x$, $r>0$ let $B(x,r):= \{y,\in \r^d,\, |x-y|\leq r\} $. Let $\lambda$ be the Lebesgue measure on $\r^\d$ and for any $A\in \mathcal{B}(\r^\d)$,  $\lambda_{A}:= \lambda(A\cap \cdot)$. Let $O_1$, $O_1$ be two metric space, $\mathcal{C}(O_1,O_2)$ is the set of continuous functions from $O_1$ into $O_2$. Finally for any $R>0$, $\rho(\cdot) \in \mathcal{C}([0,R]^\d,\r^+)$ , let\nomenclature[e8]{$\mathtt{I}_\d(\rho)$}{$:= \int_{[0,R]^\d}\rho(x) \ee^{-\sqrt{2\d}\rho(x)}dx$}\nomenclature[f3]{$\lambda $}{: the Lebesque measure on $ \r^\d$}
\nomenclature[f4]{$\lambda_A(\cdot) $}{$:=\lambda(A\cap \cdot),\,  A\in \mathcal{B}(\r^\d)$}
\begin{equation}
\label{defId}
 \mathtt{I}_\d (\rho):=\int_{[0,R]^\d} \rho(x)\ee^{-\sqrt{2\d}\rho(x)}dx .
\end{equation}

The key step of the proof of Theorem \ref{maintheorem} is the following proposition  
\begin{Proposition}
\label{queuedistrib}
There exists a constant $C^*>0$ such that:

-for any $R\geq 1$, $\epsilon>0$, there exist $l>0$ and $T_0>0$  such that:

-for any $t>T_0$ and  $\rho(\cdot) \in \mathcal{C}_R(l, \kappa_\d \log l,\log t)$, we have: \nomenclature[d2]{$C^* $}{: defined in (\ref{eqqueuedistrib})}
\begin{equation}
\label{eqqueuedistrib} \left|\P\left(\exists x\in [0,R]^\d  ,\, Y_t(x)\geq a_t+\rho(x) \right) -C^*  \mathtt{I}_\d (\rho)\right|\leq \epsilon   \mathtt{I}_\d (\rho) .
\end{equation}
\end{Proposition}

\begin{figure}[t]
\centering
\caption{}
\label{Tux12}
\includegraphics[interpolate=true,width=10cm,height=10cm]{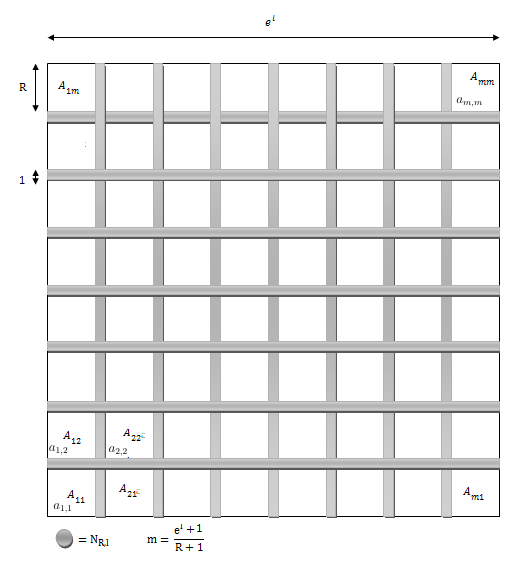} 
\end{figure}

\paragraph{Part I: Deduce Theorem \ref{maintheorem} from Proposition \ref{queuedistrib}.}
Fix $z\in \r$. Here we always assume $t\gg l\gg R>0$. By the Markov property at time $l$ and the scaling property (\ref{scalingl}) for $Y$,
\begin{eqnarray*}
\P\left( \mathbb{M}_t([0,1]^\d)\leq a_t-z\right)&=&  \P\Big( \forall x\in[0,1]^\d,\, Y^{(l)}_{t}(x) \leq a_t-z-Y_l(x)\Big)  
\\
&=&\E\Big( \P\left( \forall x\in[0,\ee^l]^\d,\, Y_{t-l}(x)\leq a_t-z-\rho(x)\right)_{\Big | \rho(\cdot)=Y_l(\cdot \ee^{-l})}   \Big).
\end{eqnarray*}
We write $[0,\ee^l]^\d= (\underset{{\bf i}}{\cup} A_{{\bf i}})\cup N_{R,l}$ with $A_{\bf i}$ and $N_{R,L}$ defined as in the figure \ref{Tux12} (pp 7), ${\bf i}\in \{1,...m\}^\d $. Clearly $\underset{R\to \infty}{\lim}\, \underset{ l\geq 1}{\sup}\, \lambda(\ee^{-l}N_{R,l})=0$, then we choose $R$ large enough such that
\begin{eqnarray*}
\P\left( \mathbb{M}_t([0,1]^\d)\leq a_t-z\right) &\simeq&\E\Big( \P\left( \forall x\in \underset{{\bf i}}{\cup} A_{{\bf i}}  ,\, Y_{t-l}(x)\leq a_t-z-\rho(x)\right)_{\Big | \rho(\cdot )=Y_l(\cdot \ee^{-l})}   \Big),
\end{eqnarray*}
where $a\simeq b$ means ``the amount $|a-b|$ can be neglected". Moreover, $x\in A_{\bf i}$ and $y \in A_{\bf j}$ with $i\neq j$ implies $|x-y|\geq 1$ and thus the processes $(Y_{t-l}(x))_{s\geq l}$ and $(Y_{t-s}(y))_{s\geq l}$ are independent.
Using the invariance by translation of $Y$ we get finally
\begin{eqnarray*}
\P\left( \mathbb{M}_t([0,1]^\d)\leq a_t-z\right) &\simeq&\E\Big( \underset{{\bf i}}{\prod} \P\left( \forall x\in    [0,R]^\d  ,\, Y_{t-l}(x)\leq a_t-z-\rho_{\bf i} (x)\right)_{\Big | \rho_{\bf i}(\cdot )=Y_l(a_{\bf i}+ \cdot \ee^{-l})}   \Big),
\end{eqnarray*}
with $a_{\bf i}:= (R+1)\left((i_1-1),...,(i_d-1)\right)$ (see figure \ref{Tux12}, (pp 7). For any ${\bf i}$ let us denote $P_{{\bf i},t}(Y_l):= \P\left( \exists x\in    [0,R]^\d  ,\, Y_{t-l}(x)\geq a_t-z-\rho_{\bf i} (x)\right)_{\Big | \rho_{\bf i}(\cdot )=Y_l(a_{\bf i}+ \cdot \ee^{-l})} $. As $\underset{x\in [0,1]^\d}{\sup} Y_l(x) \to -\infty$ when $l$ goes to infinity and $(\sup_{x\in [0,R]^\d} Y_{t-l}(x) -a_t)_{t\geq l}$ is tight (see pp 14 in \cite{DRSV12a}), we have  $\underset{l\to\infty}{\lim}\,\underset{t\to\infty}{\limsup}\, P_{{\bf i},t}(Y_l)=0$ for any ${\bf i}$, then
\begin{eqnarray*}
\P\left( \mathbb{M}_t([0,1]^\d)\leq a_t-z\right)=  \E\Big( \ee^{   \underset{{\bf i}}{\sum} \log [1- P_{{\bf i},t}(Y_l) ] } \Big)&\simeq&  \E\Big( \exp \big\{   -\underset{{\bf i}}{\sum}  P_{{\bf i},t}(Y_l)  \big\} \Big).
\end{eqnarray*}
Now we apply Proposition \ref{queuedistrib}, and get that
\begin{eqnarray*}
\P\left( \mathbb{M}_t([0,1]^\d)\leq a_t-z\right) &\simeq& \E\left( \exp \big\{ - C^*\underset{{\bf i}}{\sum}    \int_{[0,R]^\d} (-Y( a_{\bf i} +x\ee^{-l})-z)\ee^{\sqrt{2\d}[Y_l(a_{\bf i} +x\ee^{-l})+z]}dx \big\} \right)
\\
&=& \E\left( \exp \left\{ -C^* \ee^{\sqrt{2\d} z}(- z M_{l}({\cup} A_{{\bf i}} )+  M_{l}'( {\cup} A_{{\bf i}})) \right\} \right),
\end{eqnarray*}
where the last equality comes from a change of variables. Choosing  $R$ and $l$ large enough, and applying Theorem A and (\ref{Kah2}), we can affirm that
\begin{eqnarray*}
M_{l}'( {\cup} A_{{\bf i}}))\simeq M_{l}'(  [0,1]^\d)\simeq M_{\infty}'( [0,1]^\d)\quad \text{and  }\,\,   M_{l}({\cup} A_{{\bf i}} )\leq M_{l}([0,1]^\d )\simeq 0.
\end{eqnarray*}
Finally we have obtained that for $t \gg l \gg R$,
\begin{eqnarray*}
\P\left( \mathbb{M}_t([0,1]^\d)\leq a_t-z\right) &\simeq&  \E\left( \exp \left\{ -C^*\ee^{\sqrt{2\d} z}[ -z M_{l}({\cup} A_{{\bf i}} )+  M_{l}'( {\cup} A_{{\bf i}})] \right\} \right)
\\
&\simeq & \E\left( \exp \left\{ -C^*\ee^{\sqrt{2\d} z}  M_{l}'([0,1]^\d) \right\} \right).
\end{eqnarray*}
Thus we get Theorem \ref{maintheorem}. \hfill$\Box $

Before giving the main ideas to prove Proposition \ref{queuedistrib}, let us observe that Proposition \ref{queuedistrib} yields the tail distribution of $\mathbb{M}_t([0,R]^\d)$. Indeed by choosing $\rho(\cdot)=\rho$ (a constant function) we immediately obtain that
\begin{equation}
\label{tail} \underset{\rho \to \infty}{\lim}\, \underset{t\to\infty}{\lim}\, \frac{\ee^{\sqrt{2\d}\rho}}{\rho} \P\left(  \mathbb{M}_t([0,R]^\d) \geq a_t+\rho \right)=C^*R^\d.
\end{equation}
In the case of BRW, by using the ``optional lines", a result similar to (\ref{tail}) is enough to obtain the asymptotic distribution of the maximum. For our model, it is not clear whether there exists an analogue tool of ``optional lines", thus we need here a general statement as in Proposition \ref{queuedistrib}.

The proof of Proposition \ref{queuedistrib} relies on a fine understanding of the path of the particles near to the maximum (called, in the following, {\it the extremal particles}). Furthermore to establish the trajectory of an extremal particle $x$ at time $t$, we will also need to control the fluctuations of $(Y_s(y)-Y_s(x))_{s\leq t}$ for $y\in B(x,\ee^{-t})$ see Lemmas \ref{maxfx}, \ref{LemL} and Proposition \ref{tension*}.

\paragraph{Part II: Sketch of proof of Proposition \ref{queuedistrib}.}

Below are the three main steps:
\\

\noindent{\bf Step 1:} In Proposition \ref{tension*}, we establish a localization of the paths of the extremal particles. We prove that with probability close to $1$, any $x\in [0,R]^\d$ satisfying $Y_t(x)\geq  a_t+\rho(x) -1$, $x $ must also verify that $Y_{\cdot}(x)\in \DD_t^{L,\rho(x)}$ (when $L$ large) with \nomenclature[b1]{ $\DD_t^{L,\alpha}$}{$:=\left\{ f:\, \overline{f}_{t} \leq \alpha,\, \overline{f}_{[\frac{t}{2},t]}  \leq a_t+ \alpha+L,\, f_t \geq a_t+\alpha-1 \right\}$}
\begin{equation}
\label{dDdD}
 \DD_t^{\alpha,L}:= \left\{ (f_s)_{s\geq 0},\, \overline{f}_{t} \leq \alpha,\, \overline{f}_{[\frac{t}{2},t]}  \leq a_t+ \alpha+L,\, f_t \geq a_t+\alpha-1 \right\},\quad \forall L,\,\alpha,\, t>0.
\end{equation} 
See also figure \ref{Tux32} (pp 69).
\\

On the set $\{\exists x\in [0,R]^\d,\,  Y_t(x)\geq a_t+\rho(x) \}$, $\lambda_{[0,R]^\d}(\{x,\, Y_t(x)\geq a_t+ \rho(x)-1\} )>0$ and then we can write
\begin{eqnarray*}
1=  \frac{\int_{[0,R]^\d}\1_{\{ Y_t(y) \geq a_t+\rho(y)-1 \}}  dy}{\lambda_{[0,R]^\d}(\{x,\, Y_t(x)\geq a_t+ \rho(x)-1\} )} &\simeq& \int_{[0,R]^\d}\frac{\1_{\{Y_\cdot(y)\in \DD_t^{L,\rho(y)} \}}}{\lambda_{[0,R]^\d}(\{x,\, Y_t(x)\geq a_t+ \rho(x)-1\} )} dy.
\end{eqnarray*}
By taking the expectations we get that
\begin{eqnarray}
\label{dfA} A_{(\ref{dfA})}&:=&\P\left(\exists x\in [0,R]^\d,\, Y_t(x)\geq a_t+\rho(x)   \right)
\\
\nonumber &\simeq& \int_{[0,R]^\d}\E\left(  \frac{\1_{\{Y_\cdot(y)\in \DD_t^{\rho(y),L},\,  \exists x\in [0,R]^\d,\, Y_t(x)\geq a_t+\rho(x) \}}}{\lambda_{[0,R]^\d}(\{x,\, Y_t(x)\geq a_t+ \rho(x)-1\} )}\right) dy .
\end{eqnarray}

 \noindent{\bf Step 2:} The Lemma \ref{htoui} shows that on the set $\{ Y_\cdot(y)\in  \DD_t^{\rho(y),L}\}$, with an overwhelming probability, for $b$ large enough, $\{x \in [0,R]^\d,\, Y_t(x)\geq a_t+\rho(x)-1\}= \{x \in B(y,\ee^{b-t}),\,  Y_t(x)\geq a_t+\rho(x)-1\}$. In other words, only the particles close enough to $y$ are extremal. Moreover, as $\rho(\cdot )\in \mathcal{C}_R(l,\kappa_\d \log l,\log t)$ is a regular function, we may replace $\rho(x)$ by $ \rho(y)$ for any $x\in B(y,\ee^{b-t})$ and thus (\ref{dfA}) becomes
\begin{eqnarray}
\nonumber A_{(\ref{dfA})} &\simeq & \int_{[0,R]^\d}\E\left(\frac{\1_{\{Y_\cdot(y)\in \DD_t^{\rho(y),L},\,  \exists x\in B(y,\ee^{b-t}),\,   Y_t(x)\geq a_t+\rho(y) \}}}{\lambda_{ B(y,\ee^{b-t})}(\{x,\,  Y_t(x)\geq a_t+ \rho(y)-1\} )}\right) dy 
\\
\label{holes} &=& \int_{[0,R]^\d}\E_{-\rho(y)}\left(\frac{\1_{\{Y_\cdot(y)\in \DD_t^{0,L},\,  \exists x\in B(y,\ee^{b-t}),\,   Y_t(x)\geq a_t) \}}}{\lambda_{ B(y,\ee^{b-t})}(\{x,\,  Y_t(x)\geq a_t-1\} )}\right) dy=: A_{(\ref{holes})}.
\end{eqnarray}
 
 \noindent{\bf Step 3:}  We are now able to take profit from the two previous steps, using some elementary properties of $Y$. First, by the Markov property at time $t_b=t-b$, we get that
\begin{eqnarray*}
&&A_{(\ref{holes})} =
\\
&&  \int_{[0,R]^\d}\E_{-\rho(y)}\left( \1_{\{\overline{Y}_{t_b}(y)\leq 0,\, \overline{Y}_{[\frac{t}{2},t_b]}(y)\leq a_t+L \}}  \Phi_{t,y}\left[Y_{t_b}(y)-a_t-L, (Y_{t_b}(y+h)-Y_{t_b}(y))_{|h|\leq \ee^{b-t}}\right] \right) dy,
\end{eqnarray*}
with 
\begin{eqnarray*}
\Phi_{t,y}(z,(g(h))_{|h|\leq \ee^{b-t}})= \E_z\left( \frac{\1_{\{\overline{Y}_b^{(t_b)}(y)\leq 0,\, Y_{b}^{(t_b)}(y)\geq -L-1,\, \exists |h|\leq \ee^{b-t},\, Y_b^{(t_b)}(y+h) \geq -L+g(h) \}} }{\lambda_{B(0,\ee^{b-t})}( \{h,\, Y_b^{(t_b)}(y+h) \geq -L+g(h) \} ) } \right).
\end{eqnarray*}
Then using successively the following three properties of $Y$: a) the invariance by translation, b) the scaling property (\ref{scalingl}); c) Lemma \ref{rRhodesdec}; and  finally the Girsanov's transformation with density $\ee^{Y_{t_b}(0)+\d t_b}$, we obtain
\begin{eqnarray*}
A_{(\ref{holes})} &= &  \int_{[0,R]^\d}\ee^{-\sqrt{2\d} \rho(y)} t^{\frac{3}{2}}\E_{-\rho(y)}\left( \1_{\{\overline{B}_{t_b}\leq 0,\, \overline{B}_{[\frac{t}{2},t_b]}\leq a_t+L \}}  F\left[B_{t_b}(y)-a_t-L, (B_s-B_0)_{s\in [0,t]} )\right] \right) dy.
\end{eqnarray*}
with $F$ a function (defined in (\ref{defFLb})) which does not depend on $t$ and $y$ any more. Finally we conclude via a renewal theorem, see Proposition \ref{renouv}, to ensure that uniformly on $y\in [0,R]^\d$,
\begin{equation}
\label{renewouv}
\E_{-\rho(y)}\left( \1_{\{\overline{B}_{t_b}\leq 0,\, \overline{B}_{[\frac{t}{2},t_b]}\leq a_t+L \}}  F\left[B_{t_b}(y)-a_t-L, (B_s-B_0)_{s\in [0,t]} )\right] \right)\sim C^*\frac{ \rho(y)}{t^\frac{3}{2}}.
\end{equation}
From (\ref{holes}) and (\ref{renewouv}) we deduce Theorem \ref{queuedistrib}. \hfill$\Box $

\paragraph{Remark: } Let $G$ be a random variable independent of $M'([0,1]^\d)$ and satisfying $\P\left( G\leq -z\right)=  \ee^{-C^* \ee^{\sqrt{2\d}z}},\, \forall z\in \r$. By combining (\ref{tail}) and (\ref{eqmaintheorem}) we get
\begin{equation}
\label{deriveettail}
\underset{\rho \to \infty}{\lim} \frac{\ee^{\sqrt{2\d}\rho}}{\rho}\P\Big( G +\frac{1}{\sqrt{2d}}\log M'([0,1]^\d)\geq \rho\Big)= C^*.
\end{equation}
We could hope that (\ref{deriveettail}) may to identify the tail of distribution of $M'([0,1]^\d)$. As proved in \cite{BKNSW13}, in dimension one and for a particular model, we expect that $\P\left( M'([0,1]^\d)\geq \rho\right)\underset{\rho \to \infty}{\sim} \theta \rho^{-1}$. Unfortunately, (\ref{deriveettail}) is not sufficient to obtain such a result (see for instance \cite{BGT89}).
\\

The paper is organized as follows. In Section 2 we present some notations and general properties about our log-correlated Gaussian field. We prove Theorem \ref{maintheorem} in Section 3 by assuming  Proposition \ref{queuedistrib}. Section 4 is devoted to the study of  the tightness of $\mathbb{M}_t$ and the localization of the path of extremal particles. Assuming Theorem \ref{RESUME}, we prove Proposition \ref{queuedistrib} in Section 5. Finally Theorem \ref{RESUME} is proven in Section 6.

\section{ Preliminaries }
Here we state some elementary results and notations used through the paper. Let us start by a definition.
\begin{Definition}
For any domain $D\subset \r^\d$ and any $f(\cdot) \in \mathcal{C}(D,\r) $, let \nomenclature[g]{$w^{(D)}_{f(\cdot)}(\delta)$}{$:=\underset{y,\, x\in D,  \mid x-y\mid \leq \delta}{\sup} \mid f(x)-f(y)\mid $} \nomenclature[g]{$w_{g_\cdot(\cdot)}(\delta, y,t)$}{$:= \underset{s\leq t,\, x\in D,\, \mid x-y\mid \leq \delta}{\sup} \mid g_s(x)-g_s(y)\mid $}
\begin{equation}
\label{modulus}
w^{(D)}_{f(\cdot)}(\delta):=  \underset{y,\, x\in D,\,  |x-y|\leq \delta}{\sup} |f(x)-f(y)|,\qquad   w^{(D,1/3)}_{f(\cdot)}(\delta) := \underset{y,x\in D,\,  |x-y|\leq \delta}{\sup} \frac{|f(x)-f(y)|}{|x-y|^{1/3}},\qquad \forall \delta>0.
\end{equation}
For any  function $  g_\cdot(\cdot)\in \mathcal{C}(\r^+ \times D,\r)$, let
\begin{equation}
\label{remodulus}
w_{g_\cdot(\cdot)}(\delta,y,t):=  \underset{s\leq t,\,x\in D,\,   |x-y|\leq \delta}{\sup} |g_s(x)-g_s(y)|,\qquad \forall\, \delta>0,\, y\in D,\, t>0.
\end{equation}\nomenclature[g]{$w^{(D,1/3)}_f(\delta)$}{$:= \underset{ \mid x-y\mid \leq \delta}{\sup} \frac{\mid f(x)-f(y)\mid }{\mid x-y\mid^{1/3}}$}  
\end{Definition}
When $D=[0,R]^\d$ ($R>0$), we will use $w_{f(\cdot)}^{(R)}(\delta)$ and $w_{f(\cdot)}^{(R,1/3)}(\delta)$ instead of respectively $w_{f(\cdot)}^{([0,R]^\d)}(\delta)$ and $w_{f(\cdot)}^{([0,R]^\d,1/3)}(\delta)$. Similarly, when $D=B(0,b)$ ($b>0$), we denote $w_{f(\cdot)}^{(0,b)}(\delta):= w_{f(\cdot)}^{(B(0,b))}(\delta)$.
We cite (with our notations) a Lemma due to \cite{DRSV12a} . 
\begin{Lemma}[\cite{DRSV12a}]
\label{rRhodesdec}
Recall that $\mathtt{g}(\cdot):=1- \mathtt{k}(\cdot)$. For any fixed $u\neq x$, the process $(Y_t(u))_{t\geq 0}$ can be decomposed as:
\begin{eqnarray*}
Y_t(u)=P_t^{x}(u)+Z_t^{x}(u) - \zeta_t^x(u),\qquad \forall t>0,
\end{eqnarray*}
where

- $\zeta_t^{x}(u):=\sqrt{2\d} \int_0^{t}\mathtt{g}(\ee^{s}(x-u))ds$, $t>0$,\nomenclature[a1]{$\zeta_t^{x}(u) $}{$:= \sqrt{2\d} \int_0^{t}\mathtt{g}(\ee^{s}(x-u))ds$} \nomenclature[a2]{$ P_t^{x}(u)$}{$:=\int_0^{t}\mathtt{k}(\ee^{s}(x-u))dY_s(x) $} \nomenclature[a3]{$(Z_t^{x}(u))_{t\geq 0} $}{is a centered Gaussian process independent of $(Y_t(x))_{ t\geq 0}$} 

- $P_t^{x}(u):=\int_0^{t}\mathtt{k}(\ee^{s}(x-u))dY_s(x)$ is measurable with respect to the $\sigma$-algebra generated by $(Y_t(x))_{t\geq 0}$,

- the process $(Z_t^{x}(u))_{t\geq 0}$  is a centered Gaussian process independent of $(Y_t(x))_{ t\geq 0}$ with covariance kernel:
\begin{equation}
 \label{covgeneralZ}\E\left( Z_t^x(u)  Z_{t'}^x(v) \right):= \int_0^{t \wedge t'} \left[ \mathtt{k}(\ee^s(u-v))-\mathtt{k}(\ee^s(x-u))\mathtt{k}(\ee^{s}(x-v))\right]ds  ,\qquad \forall\, t,t'>0,\, x,u,v\in \r^\d.
\end{equation}
\end{Lemma}
Observe that (\ref{covgeneralZ}) implies $ (Z_t^x(x+u))_{t\in \r^+,\, u\in \r^\d }\overset{(law)}{=}(Z_t^0(u))_{t\in \r^+,\, u\in \r^\d} $ for any $x \in \r^\d$. Some simple computations lead to
\begin{Lemma}
(i) For any $x,\, u,\,v \in \r^\d$ and $t,\, t'\in \r^+$, we have:
\begin{eqnarray}
\label{covgeneralP}\E\left( P_t^x(u)  P_{t'}^x(v) \right)&:=& \int_0^{t \wedge t'} \mathtt{k}(\ee^s(x-u))\mathtt{k}(\ee^{s}(x-v))ds .
\end{eqnarray}

(ii) For any  $l>0$ the following equality holds:
\begin{equation}
\label{scalingl}
(Y_{s+l}(x)-Y_l(x))_{s\in \r^+, x\in \r^\d}=:\left(Y_s^{(l)}(x)\right)_{s\in \r^+,\,x\in \r^\d}\overset{(d)}{=} (Y_s(x\ee^{l}))_{s\in \r^+,\,x\in \r^\d}.
\end{equation} 
 
 (iii) For any $b>0$, uniformly in $u\in B(0,b)$, $\lim_{t\to\infty} \zeta_{t}(u\ee^{-t})= \sqrt{2\d} \int_{-\infty}^{0}\mathtt{g}(\ee^{v}u)dv:= \zeta(u)$.
\end{Lemma}

Finally we state a Proposition which will be used in the proof of Proposition \ref{renouv}.
\begin{Proposition}
\label{Weakconv}
Let $b,t>0$. For almost every $w\in \Omega$, $(B(0,\ee^b)\ni y \mapsto (Z_t^0(y\ee^{-t})(w))$ belongs to $\mathcal{C}(B(0,\ee^b),\r)$. Moreover when $t$ goes to infinity, the Gaussian process $(B(0,\ee^b)\ni y \mapsto (Z_t^0(y\ee^{-t}))$ converges weakly (according to the topology induced by the uniform norm in $\mathcal{C}(B(0,\ee^b),\r)$) toward the centered Gaussian process $B(0,\ee^b) \ni y\mapsto Z(y)  $ defined by:
\begin{equation}
\label{Zlimit}\E(Z(y) Z(z))=\int_{-\infty}^0 \left[ \mathtt{k}((y-z)\ee^v)-\mathtt{k}(y\ee^v)\mathtt{k}(z\ee^v)\right] dv,\qquad y,\, z\in B(0,\ee^b).
\end{equation}
\end{Proposition}
\noindent{\it Proof of Proposition \ref{Weakconv}.} By standard results on the Gaussian processes, the regularity of the kernel $\mathtt{k}$ implies the continuity of $Z_t^0(\cdot)$ (see for instance \cite{Fer75}). Then it suffices to observe that
\begin{eqnarray*}
\E\left(Z_t^{0}(y\ee^{-t})Z_t^{0}(z\ee^{-t})\right)&=& \int_{-t}^0 \mathtt{k}((y-z)\ee^v)-\mathtt{k}(y\ee^v)\mathtt{k}(z\ee^v) dv
\\
&\underset{t\to\infty}{\to}& \int_{-\infty}^0 \mathtt{k}((y-z)\ee^v)-\mathtt{k}(y\ee^v)\mathtt{k}(z\ee^v) dv.
\end{eqnarray*}
So the finite dimensional laws of $(B(0,\ee^b)\ni y \mapsto (Z_t^0(y\ee^{-t}))_{t\geq 0}$ converge to those of $B(0,\ee^b) \ni y\mapsto Z^y  $. Finally it remains to show the tightness of $(Z_t^0(\cdot))_{t\geq 0}$ which is routine (cf \cite{Bil99}) and we omit the details.
\hfill$\Box$

 \paragraph{Convention:} Throughout the paper, $c,\, c',\, c''$ denote generic constants and may change from paragraph to paragraph. These constants are independent of the parameters $t,\, l,\, R,\, L,\, b,\, M,\, \sigma...$, according to the context of the lemmas and propositions.

\section{Proof of Theorem \ref{maintheorem} assuming  Proposition \ref{queuedistrib}}

Let $l,\, R\geq 0$ satisfying $m:=\frac{\ee^l+1}{(R+1)}\in \N$.  For any $  \{1,...,m\}^\d \ni  {\bf i}= (i_1,...,i_\d)$, let

\nomenclature[c2]{$a_{{\bf i}}$}{$:= (R+1)\left((i_1-1),...,(i_\d-1)\right),\qquad (i_1,...,i_\d)=:{\bf i}\in \{1,...,m \}^\d  $}
- $a_{{\bf i}}:= (R+1)\left((i_1-1),...,(i_d-1)\right)$, which is a point of $[0,\ee^{l}]^\d$,

- $A_{{\bf i}}$  be a subset of $[0,\ee^{l}]^\d$ defined by $A_{{\bf i}}:= a_{\bf i}+ [0,R]^\d$.

As in figure \ref{Tux12} (pp 7), we also define $N_{R,l}:={[0,\ee^{l}]^\d}_{ \big \backslash \underset{{\bf i}\in \{1,...,m\}^\d}{\bigcup} A_{\bf i}}$, which corresponds to ``a buffer zone". Indeed for any $ {\bf i} \neq {\bf j}$, $d(A_{\bf i},A_{\bf j}):=\inf\{|x-y|,\, x\in A_{\bf i},\, y\in A_{\bf j}\}\geq 1 $, then $N_{R,l}$ is the minimal area needed to make sure that the values taken by the process $Y_{t}-Y_l$ inside each $A_{\bf i}$ are independent of its values on all other $A_{\bf j}$ for $j\neq i$.
\nomenclature[c3]{$N_{l,d} $}{$:=[0,\ee^{l}]^\d\backslash_{\underset{{\bf i}\in \{1,...,m \}^\d}{\bigcup} A_{\bf i}},\,\,\,\, m:=\frac{\ee^l+1}{(R+1)},\, {\bf i}\in \{1,...,m\}^\d $} 
\\

The proof of the following three lemmas are postponed at the end of this section.
\begin{Lemma}
\label{lto0}
For  any $z\in \r $, $\epsilon >0$, there exists $R_0\geq 0$ such that for any $l,\, R\geq R_0$ with $\frac{\ee^l+1}{(R+1)}\in \N$, 
\begin{eqnarray}
\label{eq0lto}
\P\left(  |z|M_l([0,1]^\d) \geq \epsilon\ee^{-\sqrt{2\d}z}  \right) + \P\left(   M_l'(\ee^{-l}N_{R,l})\geq \epsilon\ee^{-\sqrt{2\d}z}  \right) \leq \epsilon.
\end{eqnarray}
\end{Lemma}

\begin{Lemma}
\label{lto2}
For any $ \epsilon>0$ there exists $l_0>0$  such that for any $l\geq l_0$,
\begin{equation}
\label{eqlto2}
 \P\left(  w_{Y_l(\cdot)}^{(1,1/3)}(\frac{1}{l}\ee^{-l})\geq \ee^{\frac{l}{3}} \right) =\P\left(\underset{x,y\in [0,\ee^l]^\d,\, |x-y|\leq \frac{1}{l}}{\sup} \frac{\left|Y_l(\frac{x}{\ee^l})-Y_l(\frac{y}{\ee^l}) \right|}{|x-y|^{1/3}}\geq 1\right)\leq \epsilon.
\end{equation}
\end{Lemma}

\begin{Lemma}
\label{1lto}
For any $\epsilon>0$, $ a< \frac{1}{2\sqrt{2\d}}  $  there exists $l_0>0$ such that for any $l\geq l_0$,
\begin{equation}
\label{eq1lto}
\P\left(\forall x\in [0,1]^\d,\, -10\sqrt{2\d} l\leq Y_l(x)\leq -a \log l\right)\geq 1- \epsilon.
\end{equation}
\end{Lemma}
The Lemmas \ref{lto0} and \ref{1lto} are essentially contained in \cite{DRSV12a}, whereas Lemma \ref{lto2} stems from \cite{Fer75}.

Now by admitting Lemma \ref{lto0}, \ref{lto2} and \ref{1lto}, we can give the 

\noindent{\it Proof of Theorem \ref{maintheorem}.} Let $z\in \r$. Fix $\epsilon>0$. Let us choose in the following order

A) a constant $R\geq R_0 $ associated to $z,\,\epsilon$ as in Lemma \ref{lto0},

B) a constant $l_0$ associated to $\epsilon $ as in Lemma \ref{lto2},

C) a constant $l_1\geq l_0$ associated to $R,\, \epsilon $ as in Proposition \ref{queuedistrib},

D) a constant $l_2\geq l_1$ associated to $ \epsilon $ as in Lemma \ref{1lto} with $a= \kappa_\d(= \frac{1}{4\sqrt{2\d}} )$,

E) Finally a constant  $l\geq l_2+\ee^{z}$ such that $\frac{\ee^l+1}{(R+1)}\in \N$.

 According to the previous lemmas the probability of 
\begin{eqnarray}
\nonumber\mathcal{Y}_{R,z}(l):=\{    w_{Y_l(\cdot)}^{(1,1/3)}(\frac{1}{l}\ee^{-l})\leq \ee^{\frac{l}{3}},\,  |z|M_l([0,1]^\d)+ |M_l'(\ee^{-l}N_{R,l})|  \leq \epsilon \ee^{-\sqrt{2d}z},\,\qquad
\\
\label{defYRz}  \forall x\in [0,1]^\d ,\,-10\sqrt{2\d}l\leq Y_l(x)\leq - \kappa_\d \log l\}.  
\end{eqnarray} 
\nomenclature[h1]{$\mathcal{Y}_{R,z}(l)$}{$:=  \{    w_{Y_l(\cdot)}^{(1,\frac{1}{3})}(\frac{1}{l\ee^l})\leq \ee^{\frac{l}{3}},\,  \mid z \mid M_l(\ee^{-l}N_{R,l})+ \mid M_l'(\ee^{-l}N_{R,l})\mid  \leq  \epsilon \ee^{-\sqrt{2d}z},\,  -10l\leq Y_l(x)\leq - a\log l,\, \forall x \}$} is bigger than $ 1-3\epsilon$. For any $t\geq \ee^l$,
\begin{eqnarray}
\nonumber &&\P\left(\M_t([0,1]^\d)\leq a_t-z\right) \geq \P\left(\M_t([0,1]^\d)\leq a_t-z,\, \mathcal{Y}_{R,z}(l)\right)
 \\
\label{eqq123} &&\qquad\qquad\qquad \geq  \E\left(\1_{\{\M_t( \ee^{-l}{\cup} A_{\bf i} )\leq a_t-z\}};\, \mathcal{Y}_{R,z}(l) \right)-\E\left(\1_{\{ \exists x\in \ee^{-l}N_{R,l},\, Y_t(x) \geq a_t-z\}};\, \mathcal{Y}_{R,z}(l) \right).
\end{eqnarray}
Let us bound the second term in (\ref{eqq123}). By the Markov property at time $l$ and the scaling property (\ref{scalingl}) applied to the set $N_{R,l}$, we get that
\begin{eqnarray}
\nonumber &&\E\left(\1_{\{ \exists x\in \ee^{-l}N_{R,l},\, Y_t(x) \geq a_t-z\}};\, \mathcal{Y}_{R,z}(l) \right)
\\
\nonumber&=&\E\Big( \P\left(  \exists x\in \ee^{-l}N_{R,l},\, Y^{(l)}_{t-l}(x) \geq a_t-z-\chi(x)  \right)_{\Big| \chi(\cdot)=Y_l(\cdot)};\mathcal{Y}_{R,z}(l) \Big)
\\
\label{thisis} &=&\E\Big( \P\left(  \exists x\in N_{R,l},\, Y_{t-l}(x) \geq a_t-z-\chi(x)  \right)_{\Big| \chi(\cdot)=Y_l(\frac{\cdot}{\ee^l})};\mathcal{Y}_{R,z}(l) \Big).
\end{eqnarray}
We can find a collection $(y_j)_{j\in J}\in ([0,\ee^l]^\d)^J,\,  \# J<\infty $ satisfying

- for any distinct $j_1,..., j_{\d+2}\in J $, $ \underset{k=1}{\overset{\d+2}{\cap}} (y_{j_k}+[0,1]^\d)=\emptyset$,

- The set $\underset{j\in J}{\cup}(y_j+[0,1]^\d)$ is contained in the closure $\overline{N_{R,L}}$ of $N_{R,L}$. 

Moreover for $t$ sufficiently large, on $\mathcal{Y}_{R,z}(l)$, for any $j\in J$, we have $-Y_l(y_j+\frac{\cdot}{\ee^l}) - z\in \mathcal{C}_1(l, \frac{\kappa_\d}{2} \log l,\log t)$. So on $\mathcal{Y}_{R,z}(l)$, by the invariance by translation and Proposition \ref{queuedistrib} (notice that $a_{t}-a_{t-l}\to 0 $ when $t$ goes to infinity), there exists $T_0$ such that for any $t\geq T_0$, $j\in J$,
\begin{eqnarray*}
 \P\left(  \exists x\in y_j+[0,1]^\d ,\, Y_{t-l}(x) \geq a_t-z-\chi(x)  \right)_{ | \chi(\cdot)=Y_l(\frac{\cdot}{\ee^l})} 
\\
\leq (C^*+1)   \int_{x\in y_j+[0,1]^\d} (-z- Y_l( x \ee^{-l}))\ee^{\sqrt{2\d}( Y_l( x \ee^{-l})+z)}dx.
\end{eqnarray*}
Recall that $-Y_l(y_j+\frac{\cdot}{\ee^l}) - z\in \mathcal{C}_1(l, \frac{\kappa_\d}{2} \log l,\log t)$ implies $-z-Y_l(x\ee^{-l})\geq \frac{\kappa_\d}{2}\log l\geq 0$, $\forall x\in [0,1]^\d$. So the expectation in (\ref{thisis}) is smaller than
\begin{eqnarray}
\nonumber&\leq& (\d+2)(C^*+1)\E\left( \int_{x\in N_{R,l}} (-z- Y_l( x \ee^{-l}))\ee^{\sqrt{2\d}( Y_l( x \ee^{-l})+z)}dx  ;\mathcal{Y}_{R,z}(l) \right)
\\
\label{thiis}&\leq &c\E\Big(   \int_{x\in \ee^{-l}N_{R,l}} (-z-Y_l(x))\ee^{\sqrt{2\d}(Y_l(x)+z)+\d l}dx   ;\mathcal{Y}_{R,z}(l) \Big).
\end{eqnarray}
Last inequality stems from the change of variables $x\ee^{-l}\to x$. We recognize the expression of the additive martingale and  the derivative martingale as in (\ref{1.4}), therefore (\ref{thiis}) is equal to 
\begin{eqnarray*}
 \E\left(  \ee^{\sqrt{2\d}z}[(-z)M_l(\ee^{-l}N_{R,l})  +  M_{l}'(\ee^{-l}N_{R,l})]        ;\mathcal{Y}_{R,z}(l) \right) \leq  \epsilon,
\end{eqnarray*}
by definition of $ \mathcal{Y}_{R,z}(l)$ in (\ref{defYRz}). Finally 
\begin{eqnarray}
\label{eqq1234} \E\left(\1_{\{ \exists x\in \ee^{-l}N_{R,l},\, Y_t(x) \geq a_t-z\}};\, \mathcal{Y}_{R,z}(l) \right)\leq \epsilon.
\end{eqnarray}
Let us go back to (\ref{eqq123}). To treat the first term in (\ref{eqq123}), we start as before, by applying the Markov property at time $l$ and the scaling property (\ref{scalingl}). Then observing that for any ${\bf i}\neq {\bf j}$, $d(A_{\bf i},A_{\bf j})\geq 1$ and using (\ref{eqq123}) and (\ref{eqq1234}), we deduce that
\begin{eqnarray*}
&& \P\left(\M_t([0,1]^\d)\leq a_t-z\right) 
 \\
  && \geq  \E\Big(\P\left(\forall x\in {\cup} A_{\bf i},\ Y_{t-l}(x) \leq a_t-z-\chi(x)\right)_{\Big | \chi(\cdot)=Y_l(\frac{\cdot}{\ee^l})} ;\, \mathcal{Y}_{R,z}(l) \Big)-\epsilon
 \\
 && =  \E\Big(\underset{{\bf i}\in \{1,...,m\}^\d}{\prod}\P\left(\forall x\in  A_{\bf i},\ Y_{t-l}(x) \leq a_t-z-\chi(x)\right)_{\Big | \chi(\cdot)=Y_l(\frac{\cdot}{\ee^l})} ;\, \mathcal{Y}_{R,z}(l) \Big)-\epsilon.
\end{eqnarray*}  
For $t$ sufficiently large, on $\mathcal{Y}_{R,z}(l)$  we have $-Y_l(a_{\bf i}+\frac{\cdot}{\ee^l}) - z\in \mathcal{C}_R(l,\frac{\kappa_\d}{2} \log l,\log t)$, thus by Proposition \ref{queuedistrib}, there exists $T_1\geq T_0$ such that for any $t\geq T_1$, ${\bf i}\in \{1,...,m\}^\d$
\begin{eqnarray*}
\P\Big(\forall x\in  A_{\bf i},\ Y_{t-l}(x) \leq a_t-z-\chi(x)\Big) &\geq&  1- C^*(1+\epsilon)\int_{A_{\bf i}} [-z-\chi (x)]\ee^{\sqrt{2\d} (\chi(x)+z)}dx
\\
&= & 1- C^*(1+\epsilon)\ee^{\sqrt{2\d}z}(-z M_l(\ee^{-l}A_{\bf i} )+M_l'(\ee^{-l} A_{\bf i})  ).
\end{eqnarray*}
Finally we get:
\begin{eqnarray*}
 && \P\left(\M_t([0,1]^\d)\leq a_t-z\right)
  \\
  &\geq & \E\Big(\underset{{\bf i}\in \{1,...,m\}^\d}{\prod} [ 1- C^*(1+\epsilon)\ee^{\sqrt{2\d}z}(M_l'(\ee^{-l} A_{\bf i} ) -z M_l(\ee^{-l} A_{\bf i}) )];\, \mathcal{Y}_{R,z}(l) \Big)-\epsilon.
\end{eqnarray*}
On $\mathcal{Y}_{R,z}(l) $,  $\forall x\in [0,1]^\d,\,  \frac{\kappa_\d}{2} \log l\leq -Y_l(x)$, thus $\forall {\bf i}\in \{1,...,m\}^\d$ we clearly have 
\begin{eqnarray*}
M_l(\ee^{-l} A_{\bf i})+|M_l'(\ee^{-l} A_{\bf i} )| &=& \int_{\ee^{-l}A_i} (-Y_l(x)+1)\ee^{\sqrt{2\d}Y_l(x) +\d l}dx
\\
& \leq &c R^\d\frac{\log l}{l^{\kappa_\d \frac{\sqrt{2\d}}{4}}}\leq \epsilon.
\end{eqnarray*}
So we deduce
\begin{eqnarray*}
   \P\left(\M_t([0,1]^\d)\leq a_t-z\right) &\geq& \E\left(\ee^{-C^* \ee^{\sqrt{2\d}z} (1+c\epsilon) [M_l'(\ee^{-l}\underset{{\bf i}}{\cup}A_{\bf i}))-z M_l(\ee^{-l}\underset{{\bf i}}{\cup}A_{\bf i})]  }  ;\, \mathcal{Y}_{R,z}(l) \right)-\epsilon.
  \end{eqnarray*}
On $\mathcal{Y}_{R,z}(l) $, 
\begin{eqnarray*} 
  |z M_l(\ee^{-l}\underset{{\bf i}}{\cup}A_{\bf i})+ M'(\ee^{-l}\underset{{\bf i}}{\cup}A_{\bf i}))-M_l'([0,1]^\d)|\leq  |z| M_l(\ee^{-l}\underset{{\bf i}}{\cup}A_{\bf i})+ |M_l'(\ee^{-l}N_{R,l})|\leq \epsilon \ee^{-\sqrt{2d}z}  ,
\end{eqnarray*}
so $ \P\left(\M_t([0,1]^\d)\leq a_t-z\right) \geq \E\left(\ee^{-C^* \ee^{\sqrt{2d}z} (1+c\epsilon) M_l'([0,1]^\d)-2\epsilon C^* } \right)-2\epsilon $. By combining this inequality with Theorem A, we obtain the lower bound for Theorem \ref{maintheorem}.  The upper bound of (\ref{eqmaintheorem}) can be derived in the same way.
\hfill$\Box$

\subsection{Proof of Lemmas \ref{lto0}, \ref{1lto} and \ref{lto2}}
\noindent{\it Proof of Lemma \ref{lto0}.} 
By \cite{Kah85}, observe that 
\begin{eqnarray*}
 M_l([0,1]^\d) \underset{l\to \infty}{\to}0,\,\qquad \text{a.s},
\end{eqnarray*}  
which is sufficient to treat the first probability in (\ref{eq0lto}). To treat $M_l'(\ee^{-l}N_{R,l})$ we use the following fact (see \cite{DRSV12a}): For any $ \beta>0$ we can find two processes $ \tilde{Z}^\beta_l(A), Z^\beta_l(A),\, l\geq 0,\, A\in \mathcal{B}([0,1]^\d)$ satisfying

- almost surely there exists $ \beta$ large enough such that $M_l'(A)= \tilde{Z}^\beta_l(A)$, $\forall A\in \mathcal{B}([0,1]^\d)$ (see \cite{DRSV12a} pp 22),

- $\forall A\in \mathcal{B}([0,1]^\d)$, $|Z_t^\beta(A)-\tilde{Z}_t^\beta(A)|\leq \beta M_t^{\sqrt{2\d}}([0,1]^\d)$ (see pp 22, \cite{DRSV12a}),

- for any $l\geq 0$, $\forall A\in \mathcal{B}([0,1]^\d)$, $\E(Z_l^\beta(A))=\beta \lambda(A)$ (see pp 9, \cite{DRSV12a}).

Now let $\epsilon>0$. We fix $\beta >0$ large enough such that $ \P\left( \exists A\in \mathcal{B}([0,1]^\d),\, M_l'(A)\neq \tilde{Z}_l^\beta(A)\right)\leq \epsilon $, $\forall l\geq 0$. We choose $l_0>0$ large enough such that for any $l>l_0$,
\begin{eqnarray*}
\forall A\in \mathcal{B}([0,1]),\quad \P\left(|Z_l^\beta(A)-\tilde{Z}_l^\beta(A)|\geq \frac{\epsilon}{4} \right)\leq \P\left(\beta M_l^{\sqrt{2\d}}([0,1]^\d)\geq \frac{\epsilon}{4} \right)\leq \epsilon.
\end{eqnarray*}
Finally we fix $R$ large enough such that for any $l\geq l_0$, $\lambda(\ee^{-l}N_{R,l})\leq \frac{\epsilon^2}{\beta}$. We deduce that for any $l\geq l_0$,
\begin{eqnarray*}
\P\left(M_l'(\ee^{-l}N_{R,l})\geq \epsilon\right)&\leq& \P\left(M_l'(\ee^{-l}N_{R,l})\geq \epsilon,\,M_l'(\ee^{l}N_{R,l})= \tilde{Z}_l^\beta(\ee^{l}N_{R,l}) \right)+ \epsilon
\\
&\leq& \P\left(\tilde{Z}^\beta_l(\ee^{-l}N_{R,l})\geq \epsilon,\,|Z_l^\beta(\ee^{-l}N_{R,l})-\tilde{Z}_l^\beta(\ee^{-l}N_{R,l})|\leq \frac{\epsilon}{4} \right)+2\epsilon
\\
&\leq& \P\left(Z^\beta_l(\ee^{-l}N_{R,l})\geq \frac{3\epsilon}{4}  \right)+2\epsilon\leq 3\epsilon,
\end{eqnarray*}
where in the last inequality we have used the Markov inequality. 
\hfill $\Box$\\

\noindent{\it Proof of Lemma \ref{1lto}.} From Proposition 19 in \cite{DRSV12a},  for all $a\in [0, \frac{1}{2\sqrt{2\d}} )$,
\begin{eqnarray*}
\underset{t\geq 0}{\sup}(\underset{x\in [0,1]^\d}{\sup}\, Y_t(x)+ a \log (t+1))<\infty,\,\, \text{a.s}.
\end{eqnarray*}
Then by studying $\sup_{x\in [0,1]^\d}(-Y_t(x))$ we get easily Lemma \ref{1lto}.
\hfill$\Box$\\

\noindent{\it Proof of Lemma \ref{lto2}.} The proof is a consequence of  Fernique \cite{Fer75} pp 54. Let $\epsilon>0$ and $l\geq 0$. We consider
\begin{eqnarray*}
\varphi_l(h)&:=&\underset{(x,y)\in ([0,1]^\d)^2,\, |x-y|\leq h}{\sup}\sqrt{\E((Y_l(x)-Y_l(y))^2)}
\\
&=& \underset{(x,y)\in ([0,1]^\d)^2,\, |x-y|\leq h}{\sup}\sqrt{2\int_0^{l}\mathtt{g}(\ee^u[x-y])du}.
\end{eqnarray*}
As $x\mapsto \mathtt{g}(x)$ is $\mathcal{C}^1$ constant equal to $1$ outside $B(0,1)$, symmetric, with $\mathtt{g}(0)=0$ and thus $\mathtt{g}'(0)=0$, there exists $c>0$ such that for any $h>0$, $\varphi_l(h)\leq  ch\ee^l$. We imitate the proof of Theorem 4.2.2 in \cite{Fer75} and in particular use the following assertion (see pp 54 in \cite{Fer75}): {\it`` $\forall p\geq 2,\, b\geq \sqrt{1+4\d\log p},\, m\geq \frac{1}{2h}$ we have,
\begin{eqnarray}
\nonumber \P\left(\underset{x,y\in [0,1]^\d,\, |x-y|\leq h} {\sup}  \left|Y_l(x)-Y_l(y) \right|\geq b[\varphi_l(h)+2\varphi_l(\frac{1}{m})+2(2+\sqrt{2})\int_1^\infty \varphi_l(\frac{1}{m}p^{-u^2})du] \right)\qquad 
 \\
\label{boubou} \leq\left[5m^\d p^{2\d}+m^\d(2mh+1)^\d \right]\int_{b}^\infty \ee^{-\frac{u^2}{2}}du."
\end{eqnarray}}
Using $\varphi_l(h)\leq ch\ee^l$, we get
\begin{equation}
\label{boubou2}
 \P\left(\underset{x,y\in [0,1]^\d,\, |x-y|\leq h} {\sup}  \left|Y_l(x)-Y_l(y) \right|\geq c b \ee^l[h+\frac{1}{m}+\frac{1}{m p\log p}] \right) 
\leq c^\d \left[m^\d p^{2\d}+m^{2\d}h^\d  \right] \ee^{-\frac{b^2}{2}}.
\end{equation}
We set $p=2$ and for any $k\in \N$, $h_k:=\ee^{-k},  \, m_k:=2 \ee^k,\, b_k:=\sqrt{7\d k} $, then we observe that  
\begin{eqnarray}
\nonumber \sum_{k\geq l}^\infty \P\left(\underset{x,y\in [0,1]^\d,\, |x-y|\leq \ee^{-k}} {\sup}  \left|Y_l(x)-Y_l(y) \right|\geq c b_k \ee^{l-k} \right) 
&\leq&  \sum_{k\geq l}^\infty  c\ee^{\d k} \ee^{-\frac{b_k^2}{2}}
\\
\label{sum52l}&\leq&  \sum_{k\geq l}^\infty  \ee^{\d k} \ee^{-\frac{7}{2}\d k} \leq \ee^{-\frac{5}{2}\d l}.
\end{eqnarray}
Furthermore,
\begin{eqnarray*}
&&\P\left(\underset{x,y\in [0,\ee^l]^\d,\, |x-y|\leq \frac{1}{l}}{\sup} \frac{\left|Y_l(\frac{x}{\ee^l})-Y_l(\frac{y}{\ee^l}) \right|}{|x-y|^{\frac{1}{3}}}\geq 1\right) =\P\left(\underset{x,y\in [0,1]^\d,\, |x-y|\leq \frac{1}{l\ee^l}}{\sup} \frac{\left|Y_l(x)-Y_l(y) \right|}{(\ee^l|x-y|)^{\frac{1}{3}}}\geq 1\right)
\\
&&\leq \underset{k\geq l+\log l}{\sum}\P\left(\underset{x,y\in [0,1]^\d,\, \ee^{-(k+1)} \leq  |x-y|\leq \ee^{-k}}{\sup} {\left|Y_l(x)-Y_l(y) \right|}\geq {(\ee^l \ee^{-(k+1)})^{\frac{1}{3}}} \right)
\\
&&\leq \underset{k\geq l+\log l}{\sum}\P\left(\underset{x,y\in [0,1]^\d,\, |x-y|\leq \ee^{-k}}{\sup} {\left|Y_l(x)-Y_l(y) \right|}\geq cb_k \ee^{l-k}  \frac{\ee^{\frac{2}{3}(k-l)}}{c b_k } \right).
\end{eqnarray*} 
If $l$ is large enough,  $k\geq l+\log l$ implies $\ee^{\frac{2}{3}(k-l)}\geq c b_k$, thus by applying (\ref{sum52l}) we obtain 
\begin{eqnarray*}
\P\left(\underset{x,y\in [0,\ee^l]^\d,\, |x-y|\leq \frac{1}{l}}{\sup} \frac{\left|Y_l(\frac{x}{\ee^l})-Y_l(\frac{y}{\ee^l}) \right|}{|x-y|^{\frac{1}{3}}}\geq 1\right) \leq \ee^{-\frac{5}{2}\d l},
\end{eqnarray*}
from which Lemma \ref{lto2} follows. \hfill$\Box$

\section{Tightness of the maximum $\M_t$}
The main aim of this section is the following
\begin{Proposition}[Tightness]
\label{tension}
Recall that $\mathtt{I}_\d(\rho) $ is defined in (\ref{defId}). There exist $c_1,\, c_2>0 $ such that for any $l\geq 2$ we can find $T(l)>0$ so that the following inequality holds
\begin{eqnarray}
\label{eq2tension} {c_1}  \mathtt{I}_\d (\rho) \leq  \P\left( \exists x\in [0,R]^\d ,\, Y_t(x)\geq a_t+\rho(x) \right)  \leq  {c_2}   \mathtt{I}_\d (\rho),
\end{eqnarray} 
provided that $R\in [1,\log l]$, $t\geq T$ and  $\rho(\cdot)\in \mathcal{C}_R(l, \kappa_\d\log l, \log t)$.
\end{Proposition}

To obtain Proposition \ref{tension} we need more information about the path of particles $x$ such that $Y_t(x)\geq a_t + \rho(x)$. First we pay attention to the maximum on the trajectory after $\log l$.
\begin{Lemma}
\label{maxfx}
There exists $c_3>0$ such that for any $l \geq 2$, $R\geq 1$, $\rho(\cdot)\in \mathcal{C}_R(l ,\,10,\, +\infty)$,
\begin{equation}
\label{eqmaxfx}
\P\left( \exists x\in [0,R]^\d,\, \overline{Y}_{[\log l,\infty)}(x)\geq \rho(x)\right)\leq  c_3\int_{[0,R]^\d}((\log l)^\frac{3}{8}+  \rho(x)^{\frac{3}{4}})\ee^{-\sqrt{2\d}\rho(x)}dx.
\end{equation} 
\end{Lemma}
\paragraph{Remark:}  This Lemma is similar to the reasoning pp 43 in \cite{Aid11}. However because of the ``irregularity" of the function $\rho(\cdot) $, here we only control the trajectories after the time $\log l $.
\\

\noindent{\it Proof of Lemma \ref{maxfx}.} Without loss of generality we can assume that $\log l\in \N $. For any $k\in \N,\, k\geq  \log l+1 ,\, z\geq 0 $, we define
\begin{eqnarray}
\label{Delta}
 \Delta_{k,l}^z  &:=&\left\{f: \overline{f}_{[\log l,k-1)} \leq z ,\, \ \overline{f}_{[k-1,k]}\geq  z\right\} ,
\\
\label{Delta1}  \DDelta_{k,l}^z &:=&  \left\{f:\overline{f}_{[\log l,k-1)} \leq z+1, \,\overline{f}_{[k-1,k]}\geq z-1 \right\} .
\end{eqnarray}
We say that $ \Delta_{k,l}^z$ is a strong condition on the path whereas $ \DDelta_{k,l}^z$ is a weak condition on the path.  In particular, $ Y_\cdot(x)\in \Delta_{k,l}^{\rho(x)}$ and $ Y_\cdot(y)\notin \DDelta_{k,l}^{\rho(y)}$ imply that
\begin{eqnarray}
\label{strweak} \underset{s\leq k}{\sup} |Y_s(x)-\rho(x)-Y_s(y)+\rho(y)|\geq 1 .
\end{eqnarray}

Let us start with the following decomposition
\begin{eqnarray*} 
\P\left( \exists x\in [0,R]^\d,\, \overline{Y}_{[\log l,\infty)}(x)\geq \rho(x)\right)\leq \overset{+\infty}{\underset{k= \log l +1}{\sum}}\P\left(\exists x\in [0,R]^\d, Y_\cdot(x)\in \Delta_{k,l}^{\rho(x)} \right).
\end{eqnarray*}
We fix $k\geq  \log l+1  $, and study $\P\left(\exists x\in [0,R]^\d,\, Y_\cdot(x) \in \Delta_{k,l}^{\rho(x)}\right) $. \nomenclature[b2]{$ \Delta_{k,l}^z$}{$ :=\{f:\, \overline{f}_{[\log l,k-1]}\leq z ,\, \overline{f}_{[k-1,k]}\geq z\} $}\nomenclature[b3]{$\DDelta_{k,l}^z$}{$:=\{ f:\,  \overline{f}_{[\log l,k-1]}\leq z+1 ,\, \overline{f}_{[k-1,k]}\geq z-1 \}$}By continuity of $y\mapsto (Y_s(y)-\rho(y))_{s\leq k}$, if $Y_\cdot(x)\in \Delta_{k,l}^{\rho(x)}$ ($x$ satisfies the strong condition) there exist $z\in (0,R)^\d $ and $ r>0$ such that $x\in B(z,r)\subset [0,R]^\d $ and for any $y\in B(z,r)$, $y \in \DDelta_{k,l}^{\rho(y)}$ ($y$ satisfies the weak condition).

Thus on the set $\{ \exists x\in [0,R]^\d,\, Y_\cdot(x)   \in \Delta_{k,l}^{\rho(x)}\}$, we can introduce ${\bf r}>0$ (${\bf r}$ is random) the biggest radius such that

- there exists $z_{\bf r}\in [0,R]^\d$, with $B(z_{\bf r},{\bf r})\subset [0,R]^\d$,

- there exists $x_{\bf r}\in  B(z_{\bf r},{\bf r}) $ such that $Y_\cdot(x_{\bf r})\in  \Delta_{k,l}^{\rho(x_{\bf r})} $,

- for any $y\in B(z_{\bf r},{\bf r})$, $Y_\cdot(y) \in \DDelta_{k,l}^{\rho(y)}$.

Roughly speaking, the (random) radius ${\bf r}>0$ plays a quantitative role to estimate \\ $\P\left(\exists x\in [0,R]^\d,\, Y_\cdot(x)   \in \Delta_{k,l}^{\rho(x)} \right)$. Such a technique will be used several times in the sequel.

We denote by $S$ the volume of the unit ball. On the set $\{ \exists  x\in [0,R]^\d,\, Y_\cdot(x) \in \Delta_{k,l}^{\rho(x)}\} $, by definition of ${\bf r}>0$, for any ${\bf c}>0$, (${\bf c}>0$ will be determined later) we have
\begin{eqnarray*}
1&=&\frac{1}{S{\bf r}^\d}\int_{B(z_{\bf r},{\bf r}) }\1_{\{Y_\cdot(y) \in \DDelta_{k,l}^{\rho(y)}\}} dy
\\
&=& \left( \1_{\{ {\bf r}\geq \frac{\ee^{-(k+{\bf c})}}{ 4}\}}+\underset{p\geq k+{\bf c}}{\sum} \1_{\{\frac{\ee^{-(p+1)}}{4}\leq {\bf r}< \frac{\ee^{-p}}{4} \}}   \right)\frac{1}{S{\bf r}^\d}\int_{B(z_{\bf r},{\bf r}) }\1_{\{Y_\cdot(y) \in \DDelta_{k,l}^{\rho(y)} \}}dy.
\end{eqnarray*}
Taking the expectation, we obtain that
\begin{eqnarray}
&&\nonumber\P\left(\exists x\in [0,R]^\d,\, Y_\cdot(x)\in  \Delta_{k,l}^{\rho(x)}\right)\leq S^{-1}4^\d\ee^{\d (k+{\bf c})}\int_{[0,R]^\d}\P\left( Y_\cdot(y) \in \DDelta_{k,l}^{\rho(y)}  \right)   dy  
\\
\label{malin} && \qquad\qquad\qquad\qquad+\underset{p\geq (k+{\bf c})}{\sum} S^{-1}4^\d \ee^{\d (p+1)}\E\left(\1_{\{  {\bf r}\leq \frac{\ee^{-p}}{4}\}}\int_{B(z_{\bf r},{\bf r}) }\1_{\{Y_\cdot(y) \in \DDelta_{k,l}^{\rho(y)}\}}dy\right).
\end{eqnarray}

Fix $p\geq k+ {\bf c}$. For any $ R\geq 1$, on $\{    {\bf r}\leq \ee^{-p}/4\}$, $B(z_{\bf r},{\bf r})\neq [0,R]^\d$. So there is $\overline{z}\in [0,R]^\d$ with $  |\overline{z}-z_{\bf r}|\leq 2{\bf r}\leq \frac{\ee^{-p}}{2}$ and $Y_\cdot(\overline{z})\notin  \DDelta_{k,l}^{\rho(\overline{z})}$ which implies that $\underset{s\leq k}{\sup} |Y_s(\overline{z})-\rho(\overline{z})-Y_s(x_{\bf r})+\rho(x_{\bf r})|\geq 1$ (recall that $Y_\cdot(x_{\bf r})\in \Delta_{k,l}^{\rho(x_{\bf r})}$). Therefore for any $y\in B(z_{\bf r}, {\bf r}) $, by the triangular inequality we deduce that there exists $u\in [0,R]^\d$, $|u-y|\leq \ee^{-p}$ ($u$ is either $x_{\bf r}$ or $\overline{z}$) such that $\underset{s\leq k}{\sup} |Y_s(u)-\rho(u)-Y_s(y)+\rho(y)|\geq \frac{1}{2}$. To summarize,  
\begin{equation}
\label{bruxelles} \{    {\bf r}\leq \ee^{-p}/4\}\cap \{y\in B(z_{\bf r},{\bf r}) \}\subset \Big\{  \underset{ u\in B(y,\ee^{-p}),\, s\leq k}{\sup} |Y_s(y)-\rho(y)-Y_s(u)+\rho(u)|\geq \frac{1}{2}\Big\}.
 \end{equation}
Furthermore, we remark that
 
 {\bf a)} For any ${\bf c}> \log(8^3)$, as $\rho(\cdot)\in \mathcal{C}_R(l,10 , +\infty)$ and $  \ee^{-p} \leq \frac{1}{l}  $ we deduce that $ \underset{u\in B(y,\ee^{-p}) }{\sup} |\rho(y)-\rho(u)|\leq   (\ee^{-p})^{\frac{1}{3}}<\frac{1}{8}$ (recall that $p\geq \log l+{\bf c}$).

 {\bf b)} For any $y,u\in [0,R]^d$ such that $|y-u|\leq \ee^{-p}  $, as $\mathtt{k}$ is $\mathcal{C}^1$ with compact support (thus Lipschitz), according to Lemma \ref{rRhodesdec} we have
 \begin{eqnarray*}
(Y_s(u))_{s\leq k}&=&( P_s^y(u)+ Z_s^y(u)- \zeta_s^y(u))_{s\leq k}
\\
&=&( P_s^y(u)+ Z_s^y(u)+O(\ee^{k-p}))_{s\leq k}.
\end{eqnarray*}
Now, we choose ${\bf c}> \log(8^3)$ large enough such that for any $p\geq k+{\bf c} $ the $ O(\ee^{k-p})$ is smaller than $\frac{1}{2^7}$ (we stress that such $ {\bf c}$ does not depend on $ k$). Consequently by {\bf a)} and {\bf b)}, for any $p\geq k+ {\bf c}$  the event in the right-hand side of (\ref{bruxelles}) is included in
\begin{eqnarray*}
\{ \underset{u\in B(y,\ee^{-p}) ,\, s\leq k}{\sup} |P_s^y(u)-Y_s(y)|\geq 2^{-3}\} &\cup &\{ \underset{ u\in B(y,\ee^{-p}),\, s\leq k}{\sup} |Z^y_s(u)|\geq 2^{-3}\}
\\
 =\{ w_{P_{\cdot}^y(\cdot)}( \ee^{-p},y,k)  \geq 2^{-3} \}&\cup &\{ w_{Z_{\cdot}^y(\cdot)}( \ee^{-p},y,k)  \geq 2^{-3} \},\qquad ( w_\cdot(\cdot,\cdot,\cdot) \text{ is defined in (\ref{remodulus})}) .
\end{eqnarray*}
We go back to (\ref{malin}), and use the independence between $(Z^y(u))_{u\in [0,R]^d}$ and $Y_\cdot(y)$ to deduce that there are some constants $ c,\, {\bf  c}>0$ (independent of $k$) such that
\begin{eqnarray*}
&&\P\left(\exists\,   x\in [0,R]^\d,\, Y_\cdot(x) \in  \Delta_{k,l}^{\rho(x)} \right)\leq  c
\int_{[0,R]^\d}\P\left( Y_\cdot(y) \in \DDelta_{k,l}^{\rho(y)} \right)    \Big[\ee^{\d k}+
\\
&& \underset{p\geq  k+{\bf c}}{\sum} \ee^{\d p}\P\left(     w_{Z_{\cdot}^y(\cdot)}( \ee^{-p},y,k) \geq 2^{-3} \right)\Big]+ \underset{p\geq  k+{\bf c}}{\sum} \ee^{\d p} \P\left(  Y_\cdot(y) \in \DDelta_{k,l}^{\rho(y)},\,   w_{P_\cdot^y(\cdot)}( \ee^{-p},y,k)   \geq 2^{-3} \right)dy .
\end{eqnarray*}
Referring to the Appendix, by (\ref{eqpartieind}) in Lemma \ref{partieind} we get that
\begin{eqnarray*}
\underset{p\geq  k+{\bf c}}{\sum} \ee^{\d p}\P\left(  w_{Z_\cdot^y(\cdot)}( \ee^{-p},y,k)   \geq 2^{-3} \right) &= &  \underset{p\geq  k+{\bf c}}{\sum} \ee^{\d p}\P\left( \underset{|u|\leq \ee^{-p},\,s\in[0,k]}{\sup}|Z^0_s(u)|  \geq 2^{-3} \right)
\\
&\leq &  \ee^{\d  (k+{\bf c})} \underset{p\geq  k+{\bf c}}{\sum} \ee^{\d(p- (k+{\bf c}))} c_{15} \ee^{-c_{16}2^{-6} \ee^{2(p-k)}}
\\
&\leq & c\ee^{ \d{\bf c}}\ee^{\d k} = c' \ee^{\d k}.
\end{eqnarray*}
Thus we deduce that for any $k\geq \log l+1$, 
\begin{eqnarray}
\label{def12}   &&\P\left(\exists x\in [0,R]^\d,\, Y_\cdot(x) \in \Delta_{k,l}^{\rho(x)} \right)\leq (1)_k+(2)_k, \qquad \text{   with    }
\\
&& \nonumber (1)_k:= c\ee^{\d k}\int_{[0,R]^\d}\P\left( Y_\cdot(y) \in \DDelta_{k,l}^{\rho(y)} \right)   dy ,
\\
&&\nonumber    (2)_k:=c\underset{p\geq k+{\bf c}}{\sum} \ee^{\d p} \int_{[0,R]^\d}\P\left(  Y_\cdot(y) \in \DDelta_{k,l}^{\rho(y)},\,   w_{P_\cdot^y(\cdot)}(\ee^{-p},y,k)  \geq  2^{-3}\right)dy.
\end{eqnarray}
The Lemma \ref{maxfx} will be proved once the following two estimates are shown:
\begin{eqnarray}
\label{boud(1)}\underset{k\geq \log l+1}{\sum}(1)_k &\leq &c \int_{[0,R]^d}\ee^{-\sqrt{2\d} \rho(y)}dy,
\\
\label{boud(2)}\underset{k\geq \log l+1}{\sum}(2)_k  &\leq&c \int_{[0,R]^d}[(\log l)^\frac{3}{8}+ \rho(x)^\frac{3}{4}]\ee^{-\sqrt{2\d} \rho(y)}dy.
\end{eqnarray}
For any $y\in[0,R]^d$, set $T_k(y):=\inf\{ s\geq k-1,\, Y_s(y)\geq \rho(y)-1\}$,  $\tau_k(y):= \inf\{ s\geq k-1,\, B_s\geq \rho(y)-1\}$ and $\tau(y):= \inf\{ s\geq \log l,\, B_s\geq \rho(y)-1\}$.

{\bf Proof of (\ref{boud(1)})}. Fix $y\in[0,R]^d$. By Girsanov's transformation we observe that:
\begin{eqnarray*}
\P\left( Y_\cdot(y) \in \DDelta_{k,l}^{\rho(y)}\right) &=& \P\left( \overline{Y}_{[\log l,k-1)}(y) \leq \rho(y)+1 ,\,  T_k(y)\leq k \right)
=\E\left(\frac{\1_{\{ \overline{B}_{[\log l,k-1)} \leq  \rho(y)+1,\, \tau_k(y) \leq k  \}}}{\ee^{\sqrt{2\d}B_{\tau_k(y)}+\d \tau_k(y)}}\right)
\\
&\leq & c\ee^{-\sqrt{2\d}\rho(y)}\ee^{-\d k} \P\left(\overline{B}_{[\log l,k-1)} \leq  \rho(y)+1,\, \tau_k(y)\leq k \right)
\\
&=&  c\ee^{-\sqrt{2\d}\rho(y)}\ee^{-\d k} \P\left(\overline{B}_{[\log l,k-1)} \leq \rho(y)+1,\, \overline{B}_{[k-1,k]} \geq \rho(y)-1  \right).
\end{eqnarray*}
Then 
\begin{eqnarray}
\label{com2} \underset{k\geq \log l+1}{\sum}(1)_k\leq  c \int_{[0,R]^\d} \ee^{-\sqrt{2\d}\rho(y)} \E\Big( \underset{k\geq \log l+1}{\sum} \1_{\{ \overline{B}_{[\log l,k-1]} \leq  \rho(y)+1 ,\,  \overline{B}_{[k-1,k]} \geq \rho(y)-1  \}} \Big) dy .
\end{eqnarray}
To bound the expectation inside (\ref{com2}), observe that $k<\tau(y) $ implies $\sup_{s\in [\log l,k]}B_s<\rho(y)-1$ and thus $\1_{\{ \overline{B}_{[\log l,k-1]} \leq  \rho(y)+1 ,\,  \overline{B}_{[k-1,k]} \geq \rho(y)-1  \}} =0$. So by the strong Markov property at time $\tau(y)$ we obtain 
\begin{eqnarray*}
&&\E\Big( \underset{k\geq \log l+1}{\sum} \1_{\{ \overline{B}_{[\log l,k-1)}  \leq  \rho(y)+1 ,\,  \overline{B}_{[k-1,k]} \geq \rho(y )-1  \}}\Big)
\\
& &\leq  2 +\E\Big( \sum_{k\geq \tau(y)+2} \1_{\{ \overline{B}_{[\log l, \tau(y)]}\leq \rho(y)\}} \P_{B_{\tau(y)}}\big(  \overline{B}_{k-1-s}\leq \rho(y)+1,\, \overline{B}_{[k-1-s,k-s]}\geq \rho(y)-1\big) \Big)
\\
 &&\leq  2 +\E\Big( \E\big( \underset{k\geq s+2}{\sum} \1_{\{ \overline{B}_{[0, k -s)} \leq 2 ,\,  \overline{B}_{[k-s,k+1-s]} \geq 0 \}}\big)_{\big | \tau(y)=s}\Big).
\end{eqnarray*} 
Let us assume for an instant the following assertion: {\it  there is $c_{12}''>0$ such that for any $t\geq 1$,
\begin{equation}
\label{Pascimiipris}
\P\left(\overline{B}_t\leq 2,\, \overline{B}_{[t,t+1]}\geq 0\right)\leq c_{12}'' t^{-\frac{3}{2}}.
\end{equation}}
Assuming (\ref{Pascimiipris}) and recalling (\ref{com2}), we get that 
\begin{eqnarray*}
&&\underset{k\geq \log l+1}{\sum}(1)_k\leq  c \int_{[0,R]^\d} \ee^{-\sqrt{2\d}\rho(y)}  \big( 2+c_{12}''\underset{k\geq 1}{\sum} k^{-\frac{3}{2}}\big)  dy  \leq c'\int_{[0,R]^\d}\ee^{-\sqrt{2\d}\rho(y)}dy,
 \end{eqnarray*}
which proves (\ref{boud(1)}). It remains to prove (\ref{Pascimiipris}). This is a consequence of (\ref{2.6}):
\begin{eqnarray}
\nonumber \P\left(\overline{B}_t\leq 2,\, \overline{B}_{[t,t+1]}\geq 0\right)&=& \E\Big(\1_{\{\overline{B}_t\leq 2\}} \P_{B_t}(\overline{B}_1\geq 0)\Big)
\\
\nonumber &\leq & \sum_{k=0}^\infty \P\Big( \overline{B}_t\leq 2,\, B_t\in [1-k,2-k]  \Big)\P(\overline{B}_1\geq k-1)
\\
\label{simiiilaireB2}&\overset{(\ref{2.6})}{\leq} & t^{-\frac{3}{2}}c_{12} \sum_{k=0}^\infty 4(1+k)\P(|B_1|\geq k-1) \leq c_{12}'' t^{-\frac{3}{2}}.
\end{eqnarray}

{\bf Proof of (\ref{boud(2)})}. Fix  $y\in[0,R]^d$. The strategy is similar but we have to work on the event $\{ w_{P_\cdot^y(\cdot)}(\ee^{-p},y,k)  \geq  2^{-3}\}$. Let us observe that
\begin{eqnarray*}
&&\P\left(  y \in \DDelta_{k,l}^{\rho(y)},\,   w_{P_\cdot^y(\cdot)}(\ee^{-p},y,k)  \geq  2^{-3}\right)\leq P_A(y,k,p)+ P_B(y,k,p),\qquad\text{with }
\\
&&P_A(y,k,p):=\P\left(  Y_\cdot(y) \in \DDelta_{k,l}^{\rho(y)},\,   w_{P_\cdot^y(\cdot)}(\ee^{-p},y,T_k(y))  \geq  2^{-4}\right),
\\
&&P_B(y,k,p):=\P\Big(  Y_\cdot(y) \in \DDelta_{k,l}^{\rho(y)},\,  \underset{|u|\leq \ee^{-p}}{\sup}\Big|\int_{T_k(y)}^{.} \mathtt{g}(\ee^s u)dY_s(y)\Big|_{[T_k(y),k]} \geq 2^{-4}   \Big).
\end{eqnarray*}
We study first $P_B(y,k,p) $. As $ \{Y_\cdot(y)\in  \DDelta_{k,l}^{\rho(y)} \}= \{\overline{Y}_{[\log l,k-1]}(y) \leq \rho(y)+1 ,\,  T_k(y)\leq k   \} $, by the Markov property at time $T_k(y) $,
\begin{eqnarray}
\nonumber P_B(y,k,p)&=&\P\left(  Y_\cdot(y) \in \DDelta_{k,l}^{\rho(y)} \right)\P\left(   \underset{|u|\leq \ee^{-p}}{\sup}|\int_{T_k(y)}^{.} \mathtt{g}(\ee^s u)dY_s(y)|_{[T_k(y),k]} \geq 2^{-4}   \right)
\\
\label{touti} &\leq & \P\left(  Y_\cdot(y) \in \DDelta_{k,l}^{\rho(y)} \right)\P\left(   \underset{|u|\leq \ee^{-p}}{\sup}|\int_{0}^{.} \mathtt{g}(\ee^s u)dB_s(y)|_{k} \geq 2^{-5}   \right),
\end{eqnarray}
where we have used that $\underset{s\leq k}{\sup}\, \zeta_s(u)= O(\ee^{k-p})\leq 2^{-7} $, $\forall |u|\leq \ee^{-p} $. By (\ref{eqmajorAB3}) in  Lemma \ref{majorAB}, for any $p\geq k+{\bf c} $, we have
\begin{equation}
\label{youpi}   \P\left(\underset{|u|\leq  \ee^{-p}}{\sup}|\int_{0}^{.} \mathtt{g}(\ee^s u)  dB_s|_{k} \geq  2^{-5}   \right)=\P(A_{p,k,2^{-5}})\leq 
 c_{20}\exp({-c_{19}  2^{-10}\ee^{2(p-k)}}),
\end{equation}
Therefore combining (\ref{youpi}) and (\ref{touti}), we get that
\begin{eqnarray}
\nonumber \underset{k\geq \log l+1}{\sum}\underset{p\geq  (k+{\bf c})}{\sum} \int_{[0,R]^\d}\ee^{\d p} P_B(y,k,p) dy& \leq& c\underset{k\geq \log l+1}{\sum} \int_{[0,R]^\d}\ee^{\d k} \P\left( Y_\cdot(y) \in \DDelta_{k,l}^{\rho(y)}\right) dy
\\
\label{P_b}& \leq &  c' \int_{[0,R]^\d} \ee^{-\sqrt{2\d}\rho(y)} dy,
\end{eqnarray}
where we have used (\ref{boud(1)}) in the second inequality.

Now we treat $ P_A(y,k,p) $. By Girsanov's transformation (with density $ \ee^{\sqrt{2\d}B_{T_k(y)}+\d T_k(y)}$ ) we have
\begin{eqnarray}
\nonumber P_A(y,k,p)&=& \P\left(  \overline{Y}_{[\log l,k-1)}(y) \leq \rho(y)+1,\, T_k(y)\leq k,\, w_{P_\cdot^y(\cdot)}(  \ee^{-p},y,T_k(y))  \geq  2^{-4} \right)
\\
\nonumber &\leq & c\ee^{-\sqrt{2\d}\rho(y)- \d k}  \P\Big(\overline{B}_{[\log l,k-1)} \leq \rho(y)+1,\, \tau_k(y)   \leq k,\,  
\\
\label{idoine1} && \qquad\qquad\qquad\qquad   \underset{|u|\leq  \ee^{-p}  }{\sup} |\int_{0}^{.} \mathtt{g}(\ee^s u) dB_s|_{\tau_k(y)}\geq 2^{-4}  \Big).
\end{eqnarray}
By using in turn the Hölder inequality and then (\ref{eqmajorAB3}) in Lemma \ref{majorAB} (observe that $\{    \underset{|u|\leq  { \ee^{-p}}}{\sup} |\int_0^. \mathtt{g}(\ee^s u)dB_s |_{k}\geq 2^{-4} \} =A_{p,k,2^{-4}} $), we get that the probability in (\ref{idoine1}) is smaller than
\begin{eqnarray}
\nonumber &&  \P\left(\overline{B}_{[\log l,k-1)} \leq \rho(y)+1,\, \tau_k(y)\leq k\right)^\frac{3}{4}  \P\left(   \underset{|u|\leq  { \ee^{-p}}}{\sup} |\int_0^. \mathtt{g}(\ee^s u)dB_s |_{k}\geq 2^{-4}  \right)^{\frac{1}{4}}
\\
\label{idoine2} &&\leq  \P\left(\overline{B}_{[\log l,k-1)}\leq  \rho(y)+1,\,  \overline{B}_{[k-1,k]} \geq  \rho(y)-1\right)^\frac{3}{4}  c_{20}^{\frac{1}{4}}\exp(-\frac{c_{19}}{4} 2^{-8}\ee^{2(p-k)}).
\end{eqnarray}
Furthermore using the inequality (\ref{2.6}) (as in (\ref{simiiilaireB2})), we get that
\begin{eqnarray}
\nonumber && \P\left(\overline{B}_{[\log l,k-1)}\leq \rho(y)+1,\,  \overline{B}_{[k-1,k]}\geq \rho(y)-1\right) 
\\
\nonumber &&=  \E\left(  \P_{B_{\log l}-\rho(y)-1}\left( \overline{B}_{k-1-\log l}\leq 0,\,  \overline{B}_{[k-1-\log l,k-\log l]} \geq -2 \right) \right)
\\
\label{idoine3} &&\leq c\frac{\E\left( ( B_{\log l}+\rho(y)+1)\1_{\{ B_{\log l}+\rho(y)+1\geq 0\}}\right)}{(k-\log l)^\frac{3}{2}}\leq c'\frac{\rho(y) +(\log l)^\frac{1}{2}}{(k-\log l)^\frac{3}{2}} .
\end{eqnarray}

%
%
%
Finally as $\sum_{p\geq k+{\bf c}}\ee^{\d(p-k)}\exp(-\frac{c_{19}}{4}  2^{-8}\ee^{2(p-k)})  \leq c  $, gathering (\ref{idoine1}), (\ref{idoine2}) and (\ref{idoine3}) we obtain that
\begin{eqnarray}
\nonumber \underset{k\geq \log l+1}{\sum}\underset{p\geq  (k+{\bf c})}{\sum} \int_{[0,R]^\d}\ee^{\d p} P_A(y,k,p) dy& \leq&  c  \int_{[0,R]^\d}\sum_{k\geq \log l+1} \left[ \frac{\rho(y) +(\log l)^\frac{1}{2}}{(k-\log l)^\frac{3}{2}} \right]^\frac{3}{4} \ee^{-\sqrt{2\d}\rho(y)} dy
\\
\label{P_a}& \leq &   c'\int_{[0,R]^\d} ((\log l)^\frac{3}{8}+ \rho(y)^{\frac{3}{4}})\ee^{-\sqrt{2\d}\rho(y)} dy.
\end{eqnarray}
By combining (\ref{P_a}) with (\ref{P_b}), we deduce that
\begin{eqnarray*}
\nonumber\underset{k\geq \log l+1}{\sum}(2)_k  
&\leq &     c\int_{[0,R]^\d} ((\log l)^\frac{3}{8}+ \rho(y)^{\frac{3}{4}})\ee^{-\sqrt{2\d}\rho(y)} dy,
\end{eqnarray*}
which completes the proof of Lemma \ref{maxfx}.
\hfill$\Box$.
\\

The subsequent Lemma, similar to Lemma 3.3 in \cite{Aid11},  concerns the localization of the trajectory after $\frac{t}{2}$ of an extremal particle at time $t$.
\begin{Lemma}
\label{LemL}
There exist $c_4,\, c_5>0 $ such that for any $l\geq 2$ there is $T(l)>0$ so that the following inequality holds
\begin{eqnarray}
\nonumber
\P\left(\exists x\in [0,R]^\d,\, \overline{Y}_{[\log l, t]}(x)\leq \rho(x),\, \overline{Y}_{[\frac{t}{2},t]}(x) \geq a_t+\rho(x)+L,\, Y_t(x)\geq  a_t+\rho(x)\right)
\\
\label{eqLemL} \leq c_4 \ee^{-c_5 L}\int_{[0,R]^\d} (\sqrt{\log l}+\rho(x)) \ee^{-\sqrt{2\d}\rho(x)}dx ,
\end{eqnarray}
provided that $t\geq T$, $L>0,\, R\geq 1$ and $\rho(\cdot)\in \mathcal{C}_R(l, 10, +\infty)$.
\end{Lemma}
\noindent{\it Proof of Lemma \ref{LemL}.} Instead of (\ref{eqLemL}), it is sufficient to prove that there exist $c_{4*},\, c_{5*}>0 $ such that for any $l\geq 2$ there is $T(l)>0$ such that for any $t\geq T$, $L\geq 1$, $ R\geq 1$ and  $\rho(\cdot)\in \mathcal{C}_R(l, 10, +\infty)$,
\begin{eqnarray}
\nonumber
\P\left(\exists x\in [0,R]^\d,\, \overline{Y}_{[\log l, t]}(x) \leq \rho(x),\, \overline{Y}_{[\frac{t}{2},t]}(x) \in I_t(\rho(x)+L ),\, Y_t(x)\in  I_t(\rho(x))\right)
\\
\label{suffieqLemL}\leq  c_{4*}\ee^{-c_{5*} L} \int_{[0,R]^\d} (  \sqrt{\log l}+ \rho(x)) \ee^{-\sqrt{2\d}\rho(x) }dx .  
\end{eqnarray} 
Indeed let us assume (\ref{suffieqLemL}) and prove (\ref{eqLemL}). We note that the probability in (\ref{suffieqLemL}) is null when $L> -a_t+1$, so we deduce that
\begin{eqnarray*}
&&\P\left(\exists x\in [0,R]^\d,\, \overline{Y}_{[\log l, t]}(x) \leq \rho(x),\, \overline{Y}_{[\frac{t}{2},t]}(x) \geq a_t+\rho(x)+L,\, Y_t(x)\geq  a_t+\rho(x)\right)
\\
&&\leq \sum_{L'=L+1}^{-a_t+1} \sum_{u=1}^{L'}   \P\left(\exists x\in [0,R]^\d,\, \overline{Y}_{[\log l, t]}(x)\leq \rho(x),\, \overline{Y}_{[\frac{t}{2},t]}(x) \in I_t(\rho(x)+L') ,\, Y_t(x)\in  I_t(\rho(x)+u)\right) 
\\
&&\leq \sum_{L'=L+1}^{-a_t+1} \sum_{u=1}^{L'}   \P\left(\exists x\in [0,R]^\d,\, \overline{Y}_{[\log l, t]}(x)\leq \rho(x)+u,\, \overline{Y}_{[\frac{t}{2},t]}(x) \in I_t(\rho(x)+u+L'-u ),\, \right.
\\
&&\qquad \qquad\qquad\qquad\qquad\qquad  \qquad \qquad\qquad\qquad\qquad\qquad   \qquad \qquad\qquad   Y_t(x)\in  I_t(\rho(x)+u)\Big) 
\\
&&\leq \sum_{L'=L+1}^{-a_t+1} \sum_{u=1}^{L'}   c\ee^{-c' (L'-u)} \ee^{-\sqrt{2\d}u} \int_{[0,R]^\d} (\sqrt{\log l}+ \rho(x)+u) \ee^{-\sqrt{2\d}\rho(x)}dx   
\\
&&\leq c \ee^{-c'' L}\int_{[0,R]^\d} (\sqrt{\log l}+\rho(x)) \ee^{-\sqrt{2\d}\rho(x)}dx,
\end{eqnarray*}
which yields (\ref{eqLemL}).

{\bf It remains to prove (\ref{suffieqLemL}).} Let $a>0$, let us introduce (with $I_t^{\bf 1}(z):=[a_t+z-2,a_t+z+1]$),
\begin{eqnarray}
 \qquad \quad \BB_{t,a}^{z,L} &:=&\{ f: \overline{f}_{[\log l, \frac{t}{2}]} \leq z,\, \overline{f}_{[\frac{t}{2},t-a]}\in I_t(z+L ),\, \overline{f}_{[t-a,t]}\leq a_t+z+L ,\,f_t\in I_t(z )\} ,
\\
\label{6.22}
\BBi_{t,a}^{z,L} &:=&  \{f: \overline{f}_{[\log l, \frac{t}{2}]}\leq z+1,\,   \overline{f}_{[\frac{t}{2},t-a]}\in I^{\bf 1}_t(z+L) ,\, \overline{f}_{[t-a,t]}\leq a_t+z+L+1,\, f_t\in I^{\bf 1}_t(z )\}, 
\\
 \Tr_{t,a}^{z,L} &:=& \{f:  \overline{f}_{[\log l, \frac{t}{2}]} \leq z,\,   \overline{f}_{[\frac{t}{2},t-a]}\leq  a_t+z+L-1 ,\, \overline{f}_{[t-a,t]}\in I_t(z+L)\} .
\end{eqnarray}\nomenclature[b4]{$  \BBi_{t,a}^{z,L}$}{$:=  \{f:\, \overline{f}_{[\log l, \frac{t}{2}]}\leq z+1,\, \overline{f}_{[\frac{t}{2},t-a]} \in I^{\bf 1}_t( z+L),\,\overline{f}_{[t-a,t]}\leq a_t+z+L+1 ,\, f_t\in I^{\bf 1}_t(z)\} $}\nomenclature[b5]{$  \BB_{t,a}^{z,L}$}{$:=  \{f:\, \overline{f}_{[\log l, \frac{t}{2}]}\leq z,\, f_{[\frac{t}{2},t-a]} \in I_t( z+L),\,\overline{f}_{[t-a,t]}\leq a_t+z+L ,\, f_t\in  I_t(z)\} $} \nomenclature[b6]{$  \Tr_{t,a}^{z,L}$}{$:=  \{f:\, \overline{f}_{[\log l, \frac{t}{2}]}\leq z,\, \overline{f}_{[\frac{t}{2},t-a]} \leq a_t+ z+L-1,\, \overline{f}_{[t-a,t]}\in  I_t(z+L)\} $}\nomenclature[f5]{$I^{\bf 1}_t(z)$}{$:= [a_t+z-2,a_t+z+1],\quad z\in \r $}
We say that $ \BB_{t,a}^{z,L}$ and $ \Tr_{t,a}^{z,L}$ are strong conditions on the paths whereas $ \BBi_{t,a}^{z,L}$ is a weak one.

If the path of $Y_\cdot(x)$ satisfies all the conditions in the probability of (\ref{suffieqLemL}), either $Y_\cdot(x)\in  \BB_{t,a}^{\rho(x),L}$  or $Y_\cdot(x)\in \Tr_{t,a}^{\rho(x),L}$. So Lemma \ref{LemL} is a consequence of the following assertion: {\it There exists $c>0$ such that for any $l\geq 2$ there is $T(l)>0$ so that the following inequalities hold
\begin{eqnarray}
\label{intLem1}\qquad &&\P\left(\exists x\in [0,R]^\d,\, Y_\cdot(x) \in \BB_{t,a}^{\rho(x),L} \right)\leq c(1+L)a^{-\frac{1}{2}} \int_{[0,R]^d} (\sqrt{\log l}+ \rho(y))\ee^{-\sqrt{2\d}\rho(y)}dy ,
\\
\label{intLem2}\qquad && \P\left(\exists x\in [0,R]^\d,\, Y_\cdot(x) \in \Tr_{t,a}^{\rho(x),L} \right)\leq c(1+a)\ee^{-\sqrt{2\d}L} \int_{[0,R]^d}  (\sqrt{\log l}+ \rho(y))\ee^{-\sqrt{2\d}\rho(y)}dy ,
\end{eqnarray}
provided that $R\geq 1$,  $\rho(\cdot)\in \mathcal{C}_R(l, 10, +\infty)$, $t\geq T$, $L\leq -a_t+1$ and $a\in [0, \frac{t}{3}]$.}

Indeed if we choose $a=\ee^{\frac{\sqrt{2\d}}{2}\frac{2}{3} L}$ in (\ref{intLem1}) and (\ref{intLem2}), then (\ref{suffieqLemL}) follows with $c_{5*}=\sqrt{2\d}/6 $.
\\

In what follows we prove first (\ref{intLem1}) and then (\ref{intLem2}).
\\

\noindent{\bf Proof of (\ref{intLem1}).} As in the proof of Lemma \ref{maxfx}, with the same arguments, on the set $\{ \exists  x\in [0,R]^\d,\, Y_\cdot(x) \in   \BB_{t,a}^{\rho(x),L} \}$, we can define ${\bf r}>0$ be the biggest radius such that

- there exists $z_{\bf r}\in [0,R]^\d$, with $B(z_{\bf r},{\bf r})\subset [0,R]^\d$,

- there exists $x_{\bf r}\in  B(z_{\bf r},{\bf r}) $ such that $Y_\cdot(x_{\bf r})\in  \BB_{t,a}^{\rho(x_{\bf r}),L} $

- for any $y\in B(z_{\bf r},{\bf r})$, $Y_\cdot(y) \in \BBi_{t,a}^{\rho(y),L} $.

So on  $\{\exists x\in [0,R]^\d,\, Y_\cdot(x) \in \BB_{t,a}^{\rho(x),L}\}  $, $\frac{1}{S{\bf r}^\d}\int_{B(z_{\bf r},{\bf r})}\1_{\{Y_\cdot(y)\in  \BBi_{t,a}^{\rho(y),L}\}}dy=1$ (with $S$ the volume of the unit ball). Then for any ${\bf c} > 0$,
\begin{eqnarray}
\nonumber &&\P\left(\exists x\in [0,R]^\d,\, Y_\cdot(x) \in \BB_{t,a}^{\rho(x),L}\right) \leq   S^{-1}4^\d \ee^{\d(t+{\bf c})} \int_{[0,R]^\d }\P\left( Y_\cdot(y)\in \BBi_{t,a}^{\rho(y),L}\right)  dy 
\\
\label{backto} &&    +  \underset{p\geq t+ {\bf c} }{\sum} S^{-1}4^\d\ee^{\d p} \E\left( \1_{\{   {\bf r}\leq \ee^{-p}/4\}}\int_{B(z_{\bf r},{\bf r}) }\1_{\{ Y_\cdot(y)\in \BBi_{t,a}^{\rho(y),L}\}}dy   \right).
\end{eqnarray}
Reproducing the reasoning in the proof of Lemma \ref{maxfx},  we obtain that
\begin{eqnarray*}
 {\bf r}\leq \ee^{-p}/4\, \text{ and }  \,  y\in B(z_{\bf r},{\bf r})  &\Longrightarrow&     \underset{ u\in B(y,\ee^{-p}),\, s\leq t}{\sup} |Y_s(y)-\rho(y)-Y_s(u)+\rho(u)|\geq \frac{1}{2} 
 \\
 & \Longrightarrow &  w_{P_\cdot^y(\cdot)}( \ee^{-p},y,t)  \geq 2^{-3}  \text{ or } w_{Z_\cdot^y(\cdot)}( \ee^{-p},y,t)  \geq 2^{-3},
\end{eqnarray*}
where we have chosen:

- $t\geq l$ large enough such that $\forall p\geq t+{\bf c},\,  \underset{ u\in B(y,\ee^{-p})}{\sup} |\rho(y)-\rho(u)|\leq 2^{-3}  $,

- ${\bf c}$ (which does not depend on $t$) large enough such that for any $p\geq t+{\bf c}$, $\underset{|u|\leq  \ee^{-p},\, s\leq t }{\sup} |\zeta^y_{s}(u)|\leq 2^{-7}$.

Going back to (\ref{backto}) and using the independence between $(Z_\cdot^y(u))_{u\in [0,R]^\d}$ and $Y_\cdot(y)$, we deduce that
\begin{eqnarray*}
&&\P\left(\exists x\in [0,R]^\d,\, Y_\cdot(x) \in \BB_{t,a}^{\rho(x),L}\right) \leq c\int_{[0,R]^d }\P\left(Y_\cdot(y)\in \BBi_{t,a}^{\rho(y),L} \right) \Big[\ee^{\d t }  
\\
&& +  \underset{p\geq t+{\bf c}}{\sum}\ee^{\d p} \P\left(w_{Z_\cdot^y(\cdot)}( \ee^{-p},y,t) \geq  2^{-3}  \right)\Big] +   \underset{p\geq t+{\bf c}}{\sum} \ee^{\d p} \P\left( Y_\cdot(y)\in \BBi_{t,a}^{\rho(y),L},\,   w_{P_\cdot^y(\cdot)}( \ee^{-p},y,t)\geq  2^{-3}  \right)dy .
\end{eqnarray*}
By Lemma \ref{partieind}, 
\begin{eqnarray}
\nonumber \underset{p\geq t+{\bf c}}{\sum}\ee^{\d p}\P\left(w_{Z_\cdot^y(\cdot)}( \ee^{-p},y,t) \geq 2^{-3} \right) \leq  \underset{p\geq t+{\bf c}}{\sum}\ee^{\d p}c_{15}\exp(-c_{16} 2^{-6} \ee^{2(p-t)})\leq  c\ee^{\d t},
\end{eqnarray}
which implies that
\begin{eqnarray*}
&&\P\left(\exists x\in [0,R]^\d,\, Y_\cdot(x) \in \BB_{t,a}^{\rho(x),L}\right) \leq  c\ee^{\d t} \int_{[0,R]^d }\P\left( Y_\cdot(y)\in \BBi_{t,a}^{\rho(y),L}    \right) dy  
\\
&&\qquad \qquad\qquad\qquad \qquad     +   c\int_{[0,R]^d}\underset{p\geq t+{\bf c}}{\sum} \ee^{\d p} \P\left( Y_\cdot(y) \in \BBi_{t,a}^{\rho(y),L},\,   w_{P_\cdot^y(\cdot)}( \ee^{-p},y,t)\geq  2^{-3}   \right)dy .
\end{eqnarray*}
By Girsanov's transformation (with density $ \ee^{\sqrt{2\d}Y_t(y)+\d t}$) we obtain that
\begin{eqnarray}
\nonumber&&\P\left(\exists x\in [0,R]^\d,\, Y_\cdot(x)  \in \BB_{t,a}^{\rho(x),L}\right)    \leq c \int_{[0,R]^\d }\ee^{-\sqrt{2\d}\rho(y)}t^{\frac{3}{2}}  \Big[\P\left( B\in \BBi_{t,a}^{\rho(y),L}    \right)+    
\\
\label{numer1}&&\qquad \qquad\qquad\qquad \qquad   \underset{p\geq t+{\bf c}}{\sum} \ee^{\d(p-t)} \P\Big( B \in \BBi_{t,a}^{\rho(y),L},\,     \underset{|u|\leq  \ee^{-p}  }{\sup} |\int_{0}^{.} \mathtt{g}(\ee^s u) dB_s |_{t}\geq  2^{-3}   \Big)\Big]dy,
\end{eqnarray}
where we have used $\ee^{-\sqrt{2\d}Y_t(y)}\leq t^\frac{3}{2}\ee^{2-\sqrt{2\d}\rho(y)}$ for $Y_\cdot(y)\in \BBi_{t,a}^{\rho(y),L}$. By (\ref{trajA.3}), for any $y\in [0,R]^\d$,  
\begin{eqnarray}
\nonumber t^\frac{3}{2}\P\left( B\in \BBi_{t,a}^{\rho(y),L}    \right)&\leq&  c_{12}\E_{\rho(y)}\left( B_{\log l}\1_{\{B_{\log l}\geq 0\}}\right)  (1+L)a^{-\frac{1}{2}}
\\
\label{tB.16} &\leq&  c(\sqrt{\log l} +\rho(y)) (1+L)a^{-\frac{1}{2}},
\end{eqnarray}
and by (\ref{trajB.16}) (notice that $\{\underset{|u|\leq  \ee^{-p}  }{\sup} |\int_{0}^{.} \mathtt{g}(\ee^s u) dB_s |_{t}\geq  2^{-3}\}=A_{p,t,2^{-3}}  $), for any $u\in [0,R]^\d$,
\begin{eqnarray}
\nonumber &&t^\frac{3}{2} \P\Big( B \in \BBi_{t,a}^{\rho(y),L},\,     \underset{|u|\leq  \ee^{-p}  }{\sup} |\int_{0}^{.} \mathtt{g}(\ee^s u) dB_s |_{t}\geq  2^{-3}  \Big)
\\
 \nonumber & &\leq c_{22}\E_{\rho(y)}\left( B_{\log l}\1_{\{B_{\log l}\geq 0\}}\right)  (1+L)a^{-\frac{1}{2}} \exp( -\frac{c_{19}}{2}2^{-6}\ee^{2(p-t)}) 
\\
\label{tA.1}&&\leq   c(\sqrt{\log l} +\rho(y)) (1+L)a^{-\frac{1}{2}}   \exp(-\frac{c_{19}}{2}2^{-6}\ee^{2(p-t)}).
\end{eqnarray}
Finally using (\ref{tA.1}), (\ref{tB.16}) and (\ref{numer1}) we get (\ref{intLem1}).

{\bf Proof of (\ref{intLem2}). } We introduce:\nomenclature[b7]{$  \Tr_{t,a}^{z,L}(m)$}{$ :=\{f: \overline{f}_{[\log l, \frac{t}{2}]}\leq z,\,  \overline{f}_{[\frac{t}{2}, t-a+m-1]}\leq  a_t+z+L-1,\,  \overline{f}_{[t-a+m-1 , t-a+m]}\in I_t(z+L)\}$} \nomenclature[b8]{$\Tra_{t,a}^{z,L}(m) $}{$ :=\{f: \overline{f}_{[\log l, \frac{t}{2}]}\leq z+1,\,  \overline{f}_{[\frac{t}{2}, t-a+m-1]}\leq  a_t+z+L,\,  \overline{f}_{[t-a+m-1 , t-a+m]}\geq  a_t+z+L-2 \}$}
\begin{equation}
\label{6.26}
 \Tr_{t,a}^{z,L}(m) :=\{f: \overline{f}_{[\log l, \frac{t}{2}]}\leq z,\,  \overline{f}_{[\frac{t}{2}, t-a+m-1]}\leq  a_t+z+L-1,\,  \overline{f}_{[t-a+m-1 , t-a+m]}\in I_t(z+L)\}.
\end{equation}
and the associated weak condition
\begin{equation}
\label{6.27}
  \Tra_{t,a}^{z,L}(m) :=\{f: \overline{f}_{[\log l, \frac{t}{2}]}\leq z+1,\,  \overline{f}_{[\frac{t}{2}, t-a+m-1]}\leq  a_t+z+L,\,  \overline{f}_{[t-a+m-1 , t-a+m]}\geq  a_t+z+L-2 \},
\end{equation}

We decompose the event $\{\exists  x\in [0,R]^\d,\, Y_\cdot(x)\in  \Tr_{t,a}^{\rho(x),L}\} $ in $ \overset{a}{\underset{m=1}{\bigcup}}\{\exists  x\in [0,R]^\d,\, Y_\cdot(x) \in  \Tr_{t,a}^{\rho(x),L}(m) \}$. To prove (\ref{intLem2}) it is sufficient to show that for any  $m\in \{1,...,a\}$,
\begin{equation}
\label{ayout}
\P\left( \exists x\in [0,R]^\d,\, Y_\cdot(x) \in \Tr_{t,a}^{\rho(x),L}(m)\right) \leq c\ee^{-\sqrt{2\d}L}\int_{[0,R]^\d}(\rho(y)+\sqrt{\log l})\ee^{-\sqrt{2\d}\rho(y)}dy.
\end{equation}
To begin with, we reason as in the proof of (\ref{intLem1}). On the set $\{ \exists x \in [0,R]^\d,\, Y_\cdot(x) \in   \Tr_{t,a}^{\rho(x),L}(m) \}$, let ${\bf r}$ be the biggest radius such that

- there exists $z_{\bf r}\in [0,R]^\d$, with $B(z_{\bf r},{\bf r})\subset [0,R]^\d$,

- there exists $x_{\bf r}\in B(z_{\bf r},{\bf r}) $ such that $Y_\cdot(x_{\bf r})\in  \Tr_{t,a}^{\rho(x_{\bf r}),L}(m) $

- for any $y\in B(z_{\bf r},{\bf r})$, we have $Y_\cdot(y) \in \Tra_{t,a}^{\rho(y),L}(m) $.

Then for any ${\bf c}>0$, we have
\begin{eqnarray*}
&& \P\left( \exists x\in [0,R]^\d,\, Y_\cdot(x) \in \Tr_{t,a}^{\rho(x),L}(m)\right) \leq   S^{-1}4^\d \ee^{\d(t-a+m+{\bf c})}\int_{[0,R]^\d) }\P\left( Y_\cdot(y)\in  \Tra_{t,a}^{\rho(y),L}(m) \right) dy 
\\
&&\qquad\qquad   + S^{-1}4^\d \underset{p\geq t-a+m+{\bf c} }{\sum} \ee^{\d p} \E\left(\1_{\{    {\bf r}\leq \frac{\ee^{-p}}{4}\}}\int_{B(z_{\bf r}, {\bf r}) }\1_{\{Y_\cdot(y)\in  \Tra_{t,a}^{\rho(y),L}(m)\}}dy\right).
\end{eqnarray*}
Reproducing the reasoning in the proof of Lemma \ref{maxfx}, we obtain that
\begin{eqnarray*}
 {\bf r}\leq \ee^{-p}/4\,\,\, \text{and} \,\,\,  y\in B(z_{\bf r},{\bf r})  \Longrightarrow     \underset{ u\in B(y,\ee^{-p}),\, s\leq t-a+m }{\sup} |Y_s(y)-\rho(y)-Y_s(u)+\rho(u)|\geq 1/2     
 \\
  \Longrightarrow   w_{P_\cdot^y(\cdot)}( \ee^{-p},y,t-a+m)  \geq 2^{-3} \text{ or } w_{Z_\cdot^y(\cdot)}( \ee^{-p},y,t-a+m)  \geq 2^{-3},
\end{eqnarray*}
where we have chosen:

- $t\geq l+a-m$ large enough such that $\forall p\geq t-a+m+{\bf c},\,  \underset{ u\in B(y,\ee^{-p})}{\sup} |\rho(y)-\rho(u)|\leq 2^{-3}  $,

- ${\bf c}$ (which does not depend on $t-a+m$) large enough such that for any $p\geq t-a+m+{\bf c}$, $\underset{|u|\leq  \ee^{-p},\, s\leq t-a+m }{\sup} |\zeta^y_{s}(u)|\leq 2^{-7}$.

We use the independence between $(Z_\cdot^y(u))_{u\in [0,R]^\d} $ and $Y_\cdot(y)$ to deduce that
\begin{eqnarray}
\nonumber &&\P\left(\exists x\in  [0,R]^\d,\, Y_\cdot(x) \in \Tr_{t,a}^{\rho(x),L}(m)\right) \leq (\ref{defaaa*}) +(\ref{defbbb}), \qquad\text{with, }
\\
 &\label{defaaa*}&:=  c\int_{[0,R]^\d }\P\left(Y_\cdot(y)\in \Tra_{t,a}^{\rho(y),L}(m) \right)  \Big[\ee^{\d (t-a+m) }  +  \underset{p\geq t-a+m+{\bf c}}{\sum}\ee^{\d p} 
 \\
 \nonumber &&\qquad\qquad \qquad  \qquad\qquad \qquad  \qquad\qquad \qquad   \P\left(w_{Z_\cdot^y(\cdot)}( \ee^{-p},y,t-a+m) \geq 2^{-3} \right)\Big]dy ,
\\
  \label{defbbb}&&  :=\int_{[0,R]^\d}\underset{p\geq t-a+m+{\bf c}}{\sum} \ee^{\d p} \P\left( Y_\cdot(y) \in \Tra_{t,a}^{\rho(y),L}(m),\,   w_{P_\cdot^y(\cdot)}( \ee^{-p},y,t-a+m)\geq 2^{-3} \right)dy .
\end{eqnarray}
By Lemma \ref{partieind}, 
\begin{eqnarray}
\nonumber \underset{p\geq t-a+m+{\bf c}}{\sum}\ee^{\d p} \P\left(w_{Z_\cdot^y(\cdot)}( \ee^{-p},y,t-a+m) \geq 2^{-3} \right)&\leq& \underset{p\geq t-a+m+{\bf c}}{\sum}\ee^{\d p} c_{15}\exp( -c_{16}2^{-6} \ee^{2(p-(t-a+m)})
\\
&\leq& c\ee^{\d(t-a+m)}.
\end{eqnarray}
Therefore we get
\begin{equation}
\label{defaaa} (\ref{defaaa*})\leq c\ee^{\d (t-a+m)} \int_{[0,R]^\d }\P\left( Y_\cdot(y)\in \Tra_{t,a}^{\rho(y),L}(m)    \right) dy  .
\end{equation}
For any $y\in[0,R]^d$ set $T(y):= \inf\{ s\geq t-a+m-1,\, Y_s(y)\geq a_t+\rho(y)+L-2\}$.

{\bf Study of the right hand term in (\ref{defaaa}) .} Fix $p\geq t-a+m+{\bf c},\,  y\in [0,R]^\d$. Observe that
\begin{eqnarray*}
Y_\cdot(y)\in \Tra_{t,a}^{\rho(x),L}(m)  \Longleftrightarrow \overline{Y}_{[\log l,\frac{t}{2}]}(y)\leq \rho(y)+1,\, \overline{Y}_{[\frac{t}{2},t-a+m-1]}(y)\leq a_t+\rho(y)+L,\,\qquad
\\
 T(y)\leq t-a+m.
\end{eqnarray*}
Thus by Girsanov's transformation with density $ \ee^{\sqrt{2\d}Y_{T(y)}(y)+\d T(y)}$, we have
\begin{eqnarray}
\label{hatchou21} \P\left( Y_\cdot(y)\in \Tra_{t,a}^{\rho(y),L}(m)    \right)\leq  c  t^{\frac{3}{2}} \ee^{-\sqrt{2\d}(\rho(y)+L) -\d(t-a+m)} \P\left( B\in \Tra_{t,a}^{\rho(y),L}(m)  \right). 
\end{eqnarray}
According to (\ref{trajA.4}),
\begin{eqnarray}
\nonumber
t^{\frac{3}{2}} \P\left( B\in \Tra_{t,a}^{\rho(y),L}(m)    \right)&\leq & c_{12}\E_{\rho(y)}\left( B_{\log l}\1_{\{ B_{\log l}\geq 0\}}\right)
\\
\label{hatchou}&\leq& c  {(\sqrt{\log l}+\rho(y))}.
\end{eqnarray}
Finally with (\ref{hatchou}) and (\ref{hatchou21}) for any $y\in [0,R]^\d $, we have
\begin{equation}
\label{hatchou3} \P\left( Y_\cdot(y)\in \Tra_{t,a}^{\rho(y),L}(m)    \right)\leq 
c \ee^{-\sqrt{2\d}(\rho(y)+L) -\d(t-a+m)}{(\sqrt{\log l}+\rho(y))} .
\end{equation}
With (\ref{defaaa}), we immediately deduce that
\begin{eqnarray}
\label{boundA}  (\ref{defaaa*})&\leq& c  \int_{[0,R]^\d}{(\sqrt{\log l} +\rho(y))}\ee^{-\sqrt{2\d}(\rho(y)+L)}   dy.
\end{eqnarray}

{\bf Study of $(\ref{defbbb}) $.} Fix $p\geq t-a+m+{\bf c},\, y\in [0,R]^\d $. We use the independence between $(Y_s(y))_{s\leq T(y)}$ and $(Y^{(T(y))}_s(y))_{s\geq 0}$ to get that
\begin{eqnarray*}
&&\P\left( Y_\cdot(y) \in \Tra_{t,a}^{\rho(y),L}(m),\,   w_{P_\cdot^y(\cdot)}( \ee^{-p},y,t-a+m)\geq 2^{-3}  \right)\leq (\ref{defbbb})_1+(\ref{defbbb})_2,\qquad \text{ with }
\\
&&(\ref{defbbb})_1:=  \P\left( Y_\cdot(y) \in \Tra_{t,a}^{\rho(y),L}(m)\right)\P\Big(    \underset{|u|\leq \ee^{-p}  }{\sup} |\int_{T(y)}^{.} \mathtt{g}(\ee^s u) dY_s(y) |_{[T(y),t-a+m]} \geq 2^{-4}  \Big),
\\
&&(\ref{defbbb})_2:= P\left( Y_\cdot(y) \in \Tra_{t,a}^{\rho(y),L}(m)  ,\, w_{P_\cdot^y(\cdot)}(\ee^{-p},y,T(y))\geq 2^{-4} \right).
\end{eqnarray*}
According to Lemma \ref{majorAB} and (\ref{hatchou3}), we have
\begin{equation}
\label{crrt} (\ref{defbbb})_1\leq   c\ee^{-\d(t-a+m)}  \ee^{-\sqrt{2\d}(\rho(y)+L)}  (\sqrt{\log l}+\rho(y)) c_{20} \ee^{-{c_{19}}2^{-8}\ee^{2(p-(t-a+m))}}  .
\end{equation}
Concerning $ (\ref{defbbb})_2 $, we apply  Girsanov's transformation with density  $ \ee^{\sqrt{2\d}Y_{T(y)}(y)+\d T(y)}$ and (\ref{trajB.17}) (notice that $\{    \underset{|u|\leq  \ee^{-p}  }{\sup} \big|\int_{0}^{.} \mathtt{g}(\ee^s u) dB_s \big|_{t-a+m}\geq 2^{-4}  \}= A_{p,t-a+m, 2^{-4}} $) to get that 
\begin{eqnarray}
\nonumber (\ref{defbbb})_2\leq  c  \ee^{-\sqrt{2\d}(\rho(y)+L)-\d(t-a+m)} t^{\frac{3}{2}} \P\left(B\in \Tra_{t,a}^{\rho(y),L}(m) ,\,   \underset{|u|\leq  \ee^{-p}  }{\sup} \big|\int_{0}^{.} \mathtt{g}(\ee^s u) dB_s \big|_{t-a+m}\geq 2^{-4} \right)
\\
\label{crrt2}\leq c\ee^{-\sqrt{2\d}(\rho(y)+L)-\d(t-a+m)}   (\sqrt{\log l}+\rho(y)) c_{22}  \exp({-\frac{c_{19}}{2}2^{-8}\ee^{2(p-(t-a+m))}})  .
\end{eqnarray}
Combining (\ref{crrt2}) and (\ref{crrt}) we deduce that
\begin{eqnarray}
\nonumber(\ref{defbbb})&\leq&   c    \underset{p\geq t-a+m+{\bf c}}{\sum} \ee^{\d (p-(t-a+m)}\exp({- \frac{c_{19}}{2}2^{-8}\ee^{2(p-(t-a+m))}})   \int_{[0,R]^d}(\sqrt{\log l}+\rho(y))\ee^{-\sqrt{2\d}(\rho(y)+L)}    dy  
\\
\label{boundB} &&\qquad \leq  c'\ee^{-\sqrt{2\d}L} \int_{[0,R]^d}(\sqrt{\log l}+\rho(y))\ee^{-\sqrt{2\d}\rho(y)}    dy .
\end{eqnarray}
With (\ref{boundA}) and (\ref{boundB}), we get inequality (\ref{intLem2}), and therefore the proof of Lemma \ref{LemL} is completed. \hfill$\Box$
\\

Now we will tackle the proof of Proposition \ref{tension}. For any $ L,\alpha,\, t\geq 0  $, we introduce  \nomenclature[b9]{$  \DDi_t^{\alpha,L}  $}{$ := \{f: \overline{f}_{t}\leq \alpha + 1,\, \overline{f}_{[\frac{t}{2},t]} \leq a_t+\alpha+L+1,\, f_t\geq a_t+\alpha-2\}  $}   \nomenclature[b92]{$  \Loz_{t,l}^{\alpha,L}(m)  $}{$ :=  \{f:\overline{f}_{\log l}\geq \alpha,\,  \overline{f}_{[\log l,\frac{t}{2}]}\leq \alpha,\, \overline{f}_{[\frac{t}{2},t]} \leq a_t+\alpha+L,\, f_t\in I_t(\alpha+m)\}  $}
\begin{eqnarray}
\qquad  \DD_t^{\alpha,L}&:=&\{f: \overline{f}_{t}\leq \alpha,\, \overline{f}_{[\frac{t}{2},t]} \leq a_t+ \alpha+L,\, f_t\geq a_t+\alpha-1\},
\\
\label{defclude}\qquad  \DDi_t^{\alpha,L}&:=&\{ \overline{f}_{t}\leq \alpha + 1,\, \overline{f}_{[\frac{t}{2},t]} \leq a_t+\alpha+L+1,\, f_t\geq a_t+\alpha-2\},
\\
\label{Dloza} \Loz_{t,l}^{\alpha,L}(m)&:=& \{f:\overline{f}_{\log l}\geq \alpha,\,  \overline{f}_{[\log l,t]}\leq \alpha,\, \overline{f}_{[\frac{t}{2},t]} \leq a_t+\alpha+L,\, f_t\in I_t(\alpha+m)\}.
\end{eqnarray}

%

The following proposition implies Proposition \ref{tension}.
 \begin{Proposition} 
\label{tension*}
There exist $c_6,\, c_7>0$ such that:

\noindent(i) For any  $l>2$ there exists $T(l)>0$ so that the following inequality holds
\begin{equation}
\label{eq1tension}  \P\left( \exists x\in [0,R]^\d ,\, Y_t(x)\geq a_t+\rho(x) \right)  
\leq  {c_6}  \int_{[0,R]^\d} (\sqrt{\log l} + \rho (x) ) \ee^{-\sqrt{2\d}\rho(x)}dx,
\end{equation}
provided that $t\geq T$, $R\geq 1$, $\rho(\cdot)\in \mathcal{C}_R(l,10, +\infty )$.
\\

\noindent(ii)  For any  $\epsilon>0$ we can find $L,\, l_0(L)>0$ such that for any  $l> l_0$, there exists $T(l)>1$ so that the following inequality holds
\begin{equation}
\label{eq3tension}    \P\left( \exists x\in [0,R]^\d ,\, Y_t(x)\geq a_t+\rho(x),\, Y_\cdot(x)\notin \DD_t^{\rho(x),L} \right)  
\leq  \epsilon \mathtt{I}_\d(\rho),
\end{equation}
provided that $t\geq T$, $R\geq 1$, $\rho(\cdot) \in \mathcal{C}_R(l, \kappa_\d \log l, +\infty )$.
\\

\noindent(iii) For any  $l>1$ there exists $T(l)>0$ so that the following inequality holds
\begin{equation}
\label{eq2tension*}      \P\left( \exists x\in [0,R]^\d ,\, Y_t(x)\in I_t(\rho(x)) \right)\geq c_7 \mathtt{I}_\d(\rho),
\end{equation}
provided that $t\geq T$ $, R\in[1,\log l]$, $\rho(\cdot) \in \mathcal{C}_R(l,\kappa_\d \log l, \log t )$.
\end{Proposition}

Observe that(\ref{eq1tension}) gives the upper bound of Proposition \ref{tension}, (\ref{eq3tension}) ensures that for  $L,l $ large enough, with an overwhelming probability all the extremal particles $x$ satisfy  $Y_\cdot(x)\in  \DD_t^{\rho(x),L} $. Finally (\ref{eq2tension*}) is the lower bound of Proposition \ref{tension}, which will be essential to prove Proposition \ref{queuedistrib} (see (\ref{utillower})).
\\

\noindent{\it Proof of Proposition \ref{tension*}.} (\ref{eq1tension}) and (\ref{eq3tension}) can be deduced from the following two assertions:

{\it  -There exists $c_{6*}>0$ such that for any  $L,\, l>1$ there is $T>0$ so that the following inequality holds
\begin{equation}
\label{valseme1}\P\left( \exists x\in [0,R]^\d,\, Y_\cdot(x)\in  \DD_t^{\rho(x),L} \right)\leq c_{6*}(1+L)^2\int_{[0,R]^d}\rho(y)\ee^{-\sqrt{2\d}\rho(y)} dy,
\end{equation}
provided that $t\geq T$, $R>1$, $\rho(\cdot) \in \mathcal{C}_R(l, 10, +\infty )$. 
\\

-There exists $c_{6**}>0$ such that for any  $L,\, l>1$ there is $T>0$ so that the following inequality holds
\begin{equation}
\label{valseme2} \P\left(\exists x\in [0,R]^\d,\, Y_\cdot(x)\in  \Loz_{t,l}^{\rho(x),L}(m) \right)\leq  c_{6**}(1+L)\int_{[0,R]^\d} \E_{\rho(y)}\left(B_{\log l}^+ \1_{\{\underline{B}_{\log l}\leq 1\}}\right) \ee^{-\sqrt{2\d}(\rho(y)+m)} dy,
\end{equation}
provided that $t\geq T$, $R>1$ $m>0$, $\rho(\cdot) \in \mathcal{C}_R(l,10, +\infty )$.}
 \\

 {\bf Proof of (\ref{eq1tension}) and (\ref{eq3tension}) assuming (\ref{valseme1}) and (\ref{valseme2}).} We will decompose the event $\{\exists x\in [0,R]^\d,\, Y_t(x)\geq a_t+\rho(x)\}$. For any $L>0$ there exist four possible cases:
 
i) There exists $x\in [0,R]^\d$ such that $\overline{Y}_{[\log l,\infty]}(x)\geq \rho(x) $. So we define:
\begin{eqnarray*}
A&=&\{\exists x\in [0,R]^\d,\, \overline{Y}_{[\log l,\infty]}(x)\geq \rho(x)\},
\end{eqnarray*}
ii) If $A$ is not achieved we can consider the case when $\exists x\in [0,R]^\d$ such that $ \overline{Y}_{[\log l, t]}(x)\leq \rho(x)$ , $ Y_t(x)\geq  a_t+\rho(x)$ and $ \overline{Y}_{[ \frac{t}{2},t ]}(x) \geq a_t+\rho(x)+L$. So we define:
\begin{eqnarray*}
B_L&=&\{\exists x\in [0,R]^\d,\, \overline{Y}_{[\log l, t]}(x)\leq \rho(x),\, \overline{Y}_{[ \frac{t}{2},t ]}(x) \geq a_t+\rho(x)+L,\, Y_t(x)\geq  a_t+\rho(x)\},
\end{eqnarray*}
iii) If $A$ and $B_L$ are not achieved, we consider the case when $\exists x\in [0,R]^\d$ such that $  \overline{Y}_{[\log l, t]}(x)\leq \rho(x)$,  $ \overline{Y}_{[ \frac{t}{2},t ]}(x) \leq a_t+\rho(x)+L,\, Y_t(x)\geq  a_t+\rho(x)$ and $\overline{Y}_{\log l}(x)\geq \rho(x) $. So we define:
\begin{eqnarray*}
C_L&=& {\underset{m\geq  1}{\bigcup}}\{ \exists x\in[0,R]^\d,\, Y_\cdot(x)\in  \Loz_{t,l}^{\rho(x),L}(m)\},
\end{eqnarray*}
iv) Finally if $A$, $B_L$ and $C_L$ are not achieved, it remains the case when $\exists x\in [0,R]^\d$ such that $ \overline{Y}_{t}(x)\leq \rho(x)$,  $ \overline{Y}_{[ \frac{t}{2},t ]}(x) \leq a_t+\rho(x)+L,\, Y_t(x)\geq  a_t+\rho(x)$. So we define:
\begin{eqnarray*}
D_L&=&\{\exists   x\in[0,R]^\d,\,  Y_t(x)\geq a_t+\rho(x),\,  Y_\cdot(x) \DD_t^{\rho(x),L}\}  .
\end{eqnarray*}

Let $\epsilon>0 $. Recalling (\ref{eqLemL}), we fix $L\geq 1 $ large enough such that $c_4\ee^{-c_5L}\leq \epsilon $. Then we choose $l_0(L)\geq 1 $ large enough such that for any $l>l_0$ there is $T(l)>1$ such that for any $R\geq 1$, $\rho(\cdot) \in \mathcal{C}_R(l, \kappa_\d \log l, +\infty )$, $t\geq T$:

\noindent- From Lemma \ref{maxfx},
\begin{eqnarray}
\nonumber \P(A)&\leq& c_3\int_{[0,R]^\d}((\log l)^\frac{3}{8}+ \rho(y)^\frac{3}{4}) \ee^{-\sqrt{2\d}\rho(y)} dy
\\
\label{v21}&\leq & \epsilon  \int_{[0,R]^\d} \rho(y) \ee^{-\sqrt{2\d}\rho(y)} dy.
\end{eqnarray}
- From Lemma \ref{LemL},
\begin{eqnarray}
\nonumber \P(B_L)&\leq& c_4\ee^{-c_5 L}\int_{[0,R]^\d}( \sqrt{\log l} +\rho(x)) \ee^{-\sqrt{2\d}\rho(y)}dy
\\
\label{v22} &\leq &\epsilon  \int_{[0,R]^\d} \rho(y) \ee^{-\sqrt{2\d}\rho(y)} dy.
\end{eqnarray}
- From (\ref{valseme2}), 
\begin{eqnarray}
\nonumber \P(C_L)& \leq& \sum_{m\geq 0} \P\left( \exists x\in  \Loz_{t,l}^{\rho(x),L}(m)\right) 
\\
\nonumber &\leq & c_{6**}(1+L)\underset{m\geq 1}{\sum} \ee^{-\sqrt{2\d} m} \int_{[0,R]^\d} \E_{\rho(y)}\left(B_{\log l}^+ \1_{\{\underline{B}_{\log l}\leq 1\}}\right) \ee^{-\sqrt{2\d}\rho(y)} dy
\\
\label{valseme2bis}&\leq & \epsilon  \int_{[0,R]^\d} \rho(y) \ee^{-\sqrt{2\d}\rho(y)} dy.
\end{eqnarray}
In the last inequality we have used $\rho(y)\geq \kappa_\d\log l$ (as $\rho(\cdot)\in \mathcal{C}_R(l,\kappa_\d \log l,+\infty)$) which implies $\E_{\rho(y)}\left(B_{\log l}^+ \1_{\{\underline{B}_{\log l}\leq 1\}}\right) = \E(  (\rho(y)+ B_{\log l})_+\1_{\{\underline{B}_{\log l}\leq 1-\rho(y) \}} ) \leq o_l(1)\rho(y) $. 

Combining (\ref{v21}), (\ref{v22}) and (\ref{valseme2bis}) we get (\ref{eq3tension}). To obtain (\ref{eq1tension}) we use (\ref{v21}), (\ref{v22}) and (\ref{valseme2bis}) with $\epsilon=1 $, and add from (\ref{valseme1}),
\begin{eqnarray}
\label{valseme1bis}\nonumber \P\left( D_L\right) &\leq&    c_{6*}(1+L)^2\int_{[0,R]^d}\rho(y)\ee^{-\sqrt{2\d}\rho(y)} dy .
\end{eqnarray}
Thus it yields (\ref{eq1tension}) with $c_6=c_{6*}(1+L)^2+ 3$. \hfill$\Box $
\\

%
%
%
%
%
%
%
%
%
%
%
%
%
%
%

\noindent{\bf Proof of (\ref{valseme1}) and (\ref{valseme2}).} The studies of $ \P\left( \exists x\in [0,R]^\d,\, Y_\cdot(x)\in  \DD_t^{\rho(x),L} \right)$ and   \\ $ \P\left( \exists x\in [0,R]^\d,\, Y_\cdot(x)\in  \Loz_t^{\rho(x),L}(m) \right)$ are quite redundant with that of $ \P\left( \exists x\in [0,R]^\d,\, Y_\cdot(x) \in \BB_{t,a}^{\rho(x),L} \right)$ in Lemma \ref{LemL}. Then we just mention the main steps:

A) Introduce the weak condition  $   f\in \Loza_{t,l}^{\alpha,L}(m) $ ($m,\,\alpha,\, L\geq 0,\, t\geq l\geq 0 $) which is defined by \nomenclature[b94]{$ \Loza_{t,l}^{\alpha,L}(m) $}{$:=\{ f:\, \overline{f}_{\log l}\geq \alpha-1,\,  \overline{f}_{[\log l, \frac{t}{2}]}\leq \alpha + 1,\, \overline{f}_{[\frac{t}{2},t]} \leq a_t+\alpha+L+1  ,\, f_t\in I_t^{\bf 1}(\alpha+m)  \}  $}
\begin{eqnarray}
\label{Dloza2}&&  \overline{f}_{\log l}\geq \alpha-1,\,  \overline{f}_{[\log l, \frac{t}{2}]}\leq \alpha + 1,\, \overline{f}_{[\frac{t}{2},t]} \leq a_t+\alpha+L+1  ,\, f_t\in I_t^{\bf 1}(\alpha+m),
\end{eqnarray}
(Recall that $ I_t^{\bf 1}(\alpha):= [a_t+\alpha-2,a_t+\alpha+1]  $). Then in the both cases:

B) Introduce the radius $ {\bf r}>0$.

C) Make the common reasoning about the modulus of continuity of $y\mapsto (Y_s(y)-\rho(y))_{s\leq t} $

D) Decompose $Y_\cdot(u)$ by using Lemma \ref{rRhodesdec}, then precise correctly the constant ${\bf c} $ to treat the deterministic part in the modulus of continuity $w$.

E) Apply Lemma \ref{partieind} to treat the probability of $\P\left(w_{Z^y_\cdot(\cdot)}(\ee^{-p},y,t) \geq 2^{-3}\right) $.

F) Apply Girsanov's transformation with density  $ \ee^{\sqrt{2\d}Y_t(y)+\d t}$.

At the end of these steps we can affirm that (as in (\ref{numer1})):  {\it For any  $l>1$ there exists $T>0$ such that for any $R\geq 1$ $m>0$, $\rho(\cdot)\in \mathcal{C}_R(l,10, +\infty )$, $t\geq T$,  
\begin{eqnarray}
\nonumber&&\P\left(\exists x\in [0,R]^\d,\, Y_\cdot(x) \in    \DD_t^{\rho(x),L} \right)    \leq c \int_{[0,R]^\d }\ee^{-\sqrt{2\d}\rho(y)}t^{\frac{3}{2}}  \left[\P\left( B\in    \DDi_t^{\rho(y),L} \right)+   \right.
\\
\label{numer1Sq}&&\qquad \qquad\qquad\qquad \qquad   \underset{p\geq t+{\bf c}}{\sum} \ee^{\d(p-t)} \P\Big( B \in    \DDi_t^{\rho(y),L},\,     \underset{|u|\leq  \ee^{-p}  }{\sup} \big|\int_{0}^{.} \mathtt{g}(\ee^s u) dB_s \big|_{t}\geq 2^{-3}  \Big)\Big]dy ,
\end{eqnarray}
and 
\begin{eqnarray}
\nonumber&&\P\left(\exists x\in [0,R]^\d,\, Y_\cdot(x) \in  \Loz_{t,l}^{\rho(x),L}(m) \right)    \leq c \int_{[0,R]^\d }\ee^{-\sqrt{2\d}(\rho(y)+m)}t^{\frac{3}{2}}  \left[\P\left( B\in  \Loza_{t,l}^{\rho(y),L}(m)   \right)+   \right.
\\
\label{numer1Loz}&&\qquad \qquad\qquad\qquad \qquad   \underset{p\geq t+{\bf c}}{\sum} \ee^{\d(p-t)} \P\Big( B \in  \Loza_{t,l}^{\rho(y),L}(m),\,     \underset{|u|\leq \ee^{-p}  }{\sup} \big|\int_{0}^{.} \mathtt{g}(\ee^s u) dB_s \big|_{t}\geq 2^{-3}  \Big)\Big]dy .
\end{eqnarray}}
Furthermore by (\ref{trajA.5bidbis}), (\ref{trajB.18bidbis}) and (\ref{trajA.6}), (\ref{trajB.19}), noticing that $  \{ \underset{|u|\leq  \ee^{-p}  }{\sup} \big|\int_{0}^{.} \mathtt{g}(\ee^s u) dB_s \big|_{t}\geq 2^{-3} \}=A_{p,t,2^{-3}}$, we have
\begin{eqnarray*}
{t^\frac{3}{2}}\P\left( B\in \DDi_{t}^{\rho(y),L}    \right)&\leq& c_{12}\rho(y)(1+L)^2 ,
\\
{t^\frac{3}{2}}\P\Big( B \in \DDi_{t}^{\rho(y),L},\,     \underset{|u|\leq  \ee^{-p}  }{\sup} \big|\int_{0}^{.} \mathtt{g}(\ee^s u) dB_s \big|_{t}\geq 2^{-3}  \Big) &\leq &  c_{22} \rho(y)(1+L)^2\exp({-\frac{c_{19}}{2}2^{-6}\ee^{2(p-t)}}),
\end{eqnarray*}
and
\begin{eqnarray*}
{t^\frac{3}{2}}\P\left( B\in \Loza_{t,l}^{\rho(y),L}(m)    \right)&\leq&  c_{12} {(1+L)} \E_{\rho(y)}\left(B_{\log l}^+ \1_{\{\underline{B}_{\log l}\leq 1\}}\right),
\\
{t^\frac{3}{2}}\P\Big( B \in \Loza_{t,l}^{\rho(y),L}(m),\,     \underset{|u|\leq  \ee^{-p}  }{\sup} \big|\int_{0}^{.} \mathtt{g}(\ee^s u) dB_s \big|_{t}\geq 2^{-3}  \Big) &\leq &  c_{22} {(1+L)} \E_{\rho(y)}\left(B_{\log l}^+ \1_{\{\underline{B}_{\log l}\leq 1\}}\right) 
\\
&&\qquad\times\exp({- \frac{c_{19}}{2}2^{-6}\ee^{2(p-t)}}).
\end{eqnarray*}
Finally assertions (\ref{valseme1}) and (\ref{valseme2}) follow easily from (\ref{numer1Sq}) and (\ref{numer1Loz}) and the four previous inequalities. \hfill $\Box$

\noindent{ \it Proof of (\ref{eq2tension*}).} The proof relies on a second moment argument. We need some notations: 

- Let \nomenclature[i1]{$ e_s=e_s^{(t)}$}{$:= s^{\frac{1}{12}} $ when $ 0\leq s\leq \frac{t}{2}$ and $  (t-s)^\frac{1}{12}$ when $ \frac{t}{2}\leq s\leq t$}
\begin{eqnarray}
\label{defek} e_s=e_s^{(t)}:=\left\{   \begin{array}{ll}
             s^{\frac{1}{12}} & \qquad \mathrm{if}\quad      0\leq s\leq \frac{t}{2}, \\
            (t-s)^\frac{1}{12} &\qquad \mathrm{if} \quad \frac{t}{2}\leq s\leq t . \\
         \end{array}\right. 
\end{eqnarray}

- For any $x\in [0,R]^\d$ let $(A_k(x))_{k\geq 1}$ be the partition of $-x+[0,R]^\d$ defined by
\begin{equation}
\label{Aet}
 A_1(x):=[0,R]^\d\backslash_{ B(x,1)} ,\quad    A_k(x):=  B(x,\ee^{2-k})\backslash_{B(x,\ee^{1-k})}\cap[0,R]^\d,\quad \forall k\geq 2.
\end{equation} \nomenclature[c4]{$A_k(x) $}{$:=  B(x,\ee^{2-k})\backslash_{B(x,\ee^{1-k})}\cap[0,R]^\d  ,\, k\geq 2 $}\nomenclature[c5]{$ A_1(x) $}{$:= [0,R]^\d\backslash_{ B(x,1)}$}
%
In order to have good bounds in our second moment
argument, we will restrict to 'good' particles.

- Let $ D,\, L>0, k\in \{1,..., \lfloor t\rfloor\}  $, we say that $x\in [0,R]^\d $ is $  L-$\text{good}$_k$ if    \nomenclature[c6]{$L-\text{good}_k$ and  $L-\text{good}$}{ }
\begin{eqnarray}
\nonumber &&\underset{y\in A_k(x)}{\sup}|Y_k(x)-Y_k(y)|\leq e_k+\frac{D}{2}     \text{   and  }
 \\
\label{Aket} &&\qquad\qquad Y_k(x)\leq \left\{ \begin{array}{ll}  \log l^{\frac{2}{3}},\,& \text{if  } k\in \{ 1,...,5\lfloor \log l\rfloor-1 \},
 \\
  \rho(x)-4e_k + D, \, &\text{if  } k\in \{ 5\lfloor \log l \rfloor,...,\lfloor \frac{1}{2}t\rfloor-1 \},
\\
 a_t +\rho(x)+L-4e_k+D, \,&\text{if  } k\in  \{\lfloor\frac{1}{2} t\rfloor,...,\lfloor t\rfloor\}.
\end{array}\right.
\end{eqnarray}
We say that $ x $ is $L-$good  particle (we write $x$ $L-$good   or simply $x$ good if $L=0$) if $x$ is $L-$good$_k$ for any $k\in\{1,..., \lfloor t\rfloor\}  $. Notice that the ``$\frac{2}{3}$" in $\log l^{\frac{2}{3}}$ is arbitrary and any value between $\frac{1}{2}$ and $1$ could be used.  \nomenclature[i2]{$h_{good} $}{$:= \underset{x\in \boxplus_t}{\sum} \1_{\{ x \text{ good} \}}$}\nomenclature[c7]{$ \boxplus_t$}{: a regular subdivision of $ [0,R]^\d $}

- Let $\boxplus_t:=\left\{ \ee^{-t}(i_1,...,i_\d),\,\text{with }  i_j\in \{0,...,\lfloor R\ee^t\rfloor \},\, \forall j\in [1,\d]\right\}$ be a regular subdivision of $[0,R]^\d$. We notice that for any $l>0$, there exists $T>0$ such that for any $t>T,\, \rho(\cdot)\in \mathcal{C}_R(l,\kappa_\d \log l, \infty) $,
\begin{eqnarray}
\label{Riemann}\left| \ee^{-\d t} \sum_{x\in \boxplus_t} \rho(x)\ee^{-\sqrt{2\d} \rho(x)}-   \mathtt{I}_\d(\rho)\right|\leq \frac{1}{2}\mathtt{I}_\d(\rho).
\end{eqnarray}
We also notice that there exists $c>0$, such that for any $k\in \{2,...,\lfloor t\rfloor\}  $,
\begin{equation}
\label{cardinal}
\ee^{-\d(t-k)}\sum_{y\in \boxplus_t,\, y\in A_k(x)}1 \leq c.
\end{equation}
- Finally let 
\begin{equation}
\label{hgood}h_{good}:=\#\{ x\in \boxplus_t:\,\, Y_\cdot(x)\in \DD_{t }^{\rho(x),0},\,\, x\, \, \text{good}_k\,  \forall k\in [2,\lfloor t\rfloor]\}.
\end{equation}

Now we can tackle the proof of (\ref{eq2tension*}). We fix $L=0$. By Corollary \ref{C.3}, there exists $c,\, c'>0$ and $D>0$ large ($D$ from (\ref{Aket})) such that for any $t  \geq 1$ and $\rho(x) \in [\kappa_\d \log l ,\log t]$,  
\begin{equation}
\label{bibt}
c'\ee^{-\d t} \rho(x)\ee^{-\sqrt{2\d}\rho(x)} \leq \P\left( Y_\cdot(x)\in  \DD_{t}^{\rho(x),0}, \,x\,    \text{good}_k\, \forall k\in [2,\lfloor t\rfloor]\right)\leq  \P\left( Y_\cdot(x)\in  \DD_{t}^{\rho(x),0}\right) \leq c \ee^{-\d t} \rho(x)\ee^{-\sqrt{2\d}\rho(x)} .
\end{equation}
So by combining with (\ref{Riemann}) we get that
\begin{equation}
\label{moment1h}
\frac{c'}{2}\mathtt{I}_\d(\rho) \leq  \E\left(h_{good}\right)\leq c(1+\frac{1}{2}) \mathtt{I}_\d(\rho).
\end{equation}
We look at the second moment of $h_{good}$. We recall that for any $x\in [0,R]^\d$, $\#\{ y\in \boxplus_t:\,\, |x-y|\leq \ee^{2-t}\} \leq 2^\d\ee^{2\d}$. Recall also that $|x-y|\geq 1$ implies that the process $Y_.(x)$ and $Y_.(y)$ are independent. So we deduce that
\begin{eqnarray}
\nonumber \E\left( h_{good}^2\right)&\leq&  2^\d\ee^{2\d}\E\left( h_{good}\right) +  \E\Big(\underset{x,y \in \boxplus_t,\, |x-y|\geq \ee^{2-t}}{\sum}\1_{\{Y_\cdot(x)\in \DD_{t}^{\rho(x),0},\, Y_\cdot(y)\in \DD_{t}^{\rho(y),0} ,\, x,y\,\,  \text{good}_k \forall k\in [2,\lfloor t\rfloor]  \}} \Big) 
\\
\nonumber&\leq&  2^\d\ee^{2\d} \E\left( h_{good}\right) +  \underset{x,y \in \boxplus_t,\, |x-y|\geq 1}{\sum} \P\left(Y_\cdot(x)\in \DD_{t}^{\rho(x),0}\right) \P\left( Y_\cdot(y)\in \DD_{t}^{\rho(y),0} \right) + 
\\
\label{h2good1} &&\qquad\qquad\qquad\qquad \E\Big(\underset{x \in \boxplus_t}{\sum} \sum_{y\in \boxplus_t,\, \ee^{2-t}\leq |y-x|\leq 1}  \1_{\{Y_\cdot(x)\in \DD_{t}^{\rho(x),0},\, Y_\cdot(y)\in \DD_{t}^{\rho(y),0},\,  x \,\, \text{good}_k \forall k\in [2,\lfloor t\rfloor]  \}} \Big).
\end{eqnarray}
By inequality (\ref{bibt}),
\begin{eqnarray*}
\underset{x,y \in \boxplus_t,\, |x-y|\geq 1}{\sum} \P\left(Y_\cdot(x)\in \DD_{t}^{\rho(x),0}\right)\P\left(Y_\cdot(y)\in \DD_{t}^{\rho(y),0}\right)  &\leq& \ee^{-2\d t}\underset{x,y \in \boxplus_t,\, |x-y|\geq 1}{\sum} \rho(x)\rho(y)\ee^{-\sqrt{2\d}[\rho(x)+\rho(y)]}
\\
& \leq &c\mathtt{I}_\d(\rho)^2\leq  c \mathtt{I}_\d(\rho),
\end{eqnarray*}
(observe that $R\leq \log l$ and $ \rho(\cdot) \in \mathcal{C}_R(l,\kappa_\d \log l,\infty)$ imply $\mathtt{I}_\d(\rho)\leq 1$). Going back to (\ref{h2good1}) we get     
\begin{eqnarray}
\nonumber \E\left( h_{good}^2\right) &\leq &  c\E\left( h_{good}\right) +  \sum_{k=2}^{t-1} \E\left(  \sum_{x,y\in \boxplus_t} \1_{\{  Y_\cdot(x)\in \DD_{t}^{\rho(x),0},\, Y_\cdot(y)\in \DD_{t}^{\rho(y),0},\,  x  \text{ good}_k,\,  y\in A_k(x)   \}} \right)
\\
\label{h2good}&:= &  c\E\left( h_{good}\right) +   \sum_{k=2}^{t-1}(\ref{h2good})_k   .
\end{eqnarray}
Let us study $(\ref{h2good})_k$. For any $2\leq k\leq t-1$, $x\in [0,R]^\d$ the process $(Y_s^{(k)}(y))_{s\leq t-k,\, y\in A_k(x)}$ is independent of the sigma-field
\begin{eqnarray*}
\mathcal{G}_k(x):=\sigma\left(Y_s(y),\, s\leq k,\, y\in A_k(x),\,\, Y_s(x),\, s\in \r^+\right).
\end{eqnarray*}
By the Markov property at time $k$, $(\ref{h2good})_k$ is equal to
\begin{eqnarray}
\nonumber && \sum_{x,y\in \boxplus_t,\, y\in A_k(x)} \E\left(\1_{\{Y_\cdot(x)\in  \DD^{\rho(x),0}_{t},\, x\, \text{good}_k \}} \P_{Y_k(y)}\big[ \overline{Y}_{\frac{t}{2}-k}(y)\leq \rho(y),\, \overline{Y}_{[\frac{t}{2}-k,t-k]}(y)\leq a_t+\rho(y),\, \right. \qquad 
\\
\nonumber &&\qquad\qquad\qquad\qquad\qquad\qquad\qquad\qquad \qquad\qquad\qquad\qquad    Y_{t-k}(y)\in I_t(\rho(y)) \big] \Big), 
\end{eqnarray}
Now by using the Girsanov's transformation (with density $\ee^{\sqrt{2\d} Y_{t-k}(y)+\d(t-k)} $, recall also that $\ee^{-\sqrt{2\d}Y_{t-k}(y)}\leq t^\frac{3}{2}\ee^{-\sqrt{2\d}\rho(y)}$ when  $Y_{t-k}(y)\in I_t(\rho(y))$), we deduce that $(\ref{h2good})_k$ is smaller than
\begin{eqnarray}
\label{h2goodk1} && \sum_{x,y\in \boxplus_t,\, y\in A_k(x)} \E\Big(\1_{\{Y_\cdot(x)\in  \DD^{\rho(x),0}_{t},\, x\, \text{good}_k \}} c\ee^{-\sqrt{2\d}(\rho(y)-Y_k(y))}\ee^{-\d(t-k)} t^\frac{3}{2} (\ref{h2goodk1})_{y,t,k}\Big),\quad \text{with}
\\
\nonumber &&\qquad \qquad \qquad\qquad (\ref{h2goodk1})_{y,t,k}:= \P_{Y_k(y)-\rho(y)}\big[ \overline{B}_{\frac{t}{2}-k}\leq 0,\, \overline{B}_{[\frac{t}{2}-k,t-k]}\leq a_t,\,B_{t-k}\in I_t(0) \big],
\end{eqnarray}
when $k\leq t/2$ and 
\begin{eqnarray}
\label{h2goodk2}  c\sum_{x,y\in \boxplus_t,\, y\in A_k(x)} \E\left(\1_{\{Y_\cdot(x)\in  \DD^{\rho(x),0}_{t},\, x\, \text{good}_k \}}   \ee^{-\sqrt{2\d}(\rho(y)-Y_k(y))}\ee^{-\d(t-k)} t^\frac{3}{2} \right),
\end{eqnarray}
when $t-1>k\geq t/2$. To treat $(\ref{h2good})_k$ we distinguish four cases.

\paragraph{  a)  $k\leq 5\log l$.} Let $x,y\in \boxplus_t$ with $y\in A_k(x)$, by (\ref{trajA.5bidbis}) in Lemma \ref{trajbrowseul} we have
\begin{eqnarray*}
&&   c\ee^{-\sqrt{2\d}(\rho(y)-Y_k(y))}\ee^{-\d(t-k)} t^\frac{3}{2}  (\ref{h2goodk1})_{y,t,k}
\\
&&\leq c'\ee^{-\d(t-k)}   (\rho(y)-Y_k(y))\ee^{-\sqrt{2\d}(\rho(y)-Y_k(y))}\1_{\{ \rho(y)-Y_k(y)\geq 0\}}
\\
&&\leq c''\ee^{-\d(t-k)}  \ee^{- (\rho(y)-Y_k(y))} .
\end{eqnarray*}
In addition if $x$ is $\text{good}_k$, we can ensure that 
$\rho(y)-Y_k(y)\geq \rho(y)-e_k-\frac{D}{2}-Y_k(x) \geq \frac{\kappa_\d}{2}\log l-\frac{D}{2}$. Finally by combining (\ref{h2goodk1}), (\ref{bibt}) and (\ref{cardinal}), we deduce that
\begin{eqnarray}
\nonumber(\ref{h2good})_k&\leq&  c\frac{1}{l^\frac{\kappa_\d}{2}}\ee^{-\d(t-k)+\frac{D}{2}}  \sum_{x,y\in \boxplus_t,\, y\in A_k(x)} \E\left(\1_{\{ Y_\cdot(x)\in  \DD^{\rho(x),0}_{t},\, x\, \text{good}_k \}} \right)
\\
\nonumber&\leq&   c'\frac{1}{l^\frac{\kappa_\d}{2}}\ee^{D}  \ee^{-\d t} \sum_{x \in \boxplus_t }\rho(x)\ee^{-\sqrt{2\d}\rho(x)}\ee^{-\d(t-k)}\sum_{y\in \boxplus_t,\, y\in A_k(x)}1
\\
\label{a=ine}&\leq &c''\frac{1}{l^\frac{\kappa_\d}{2}}\ee^{ D}\mathtt{I}_\d(\rho).
\end{eqnarray}

\paragraph{  b)  $5 \log l \leq k\leq t/4$.} Let $x,y\in \boxplus_t$ with $y\in A_k(x)$, by (\ref{trajA.5bidbis}) in Lemma \ref{trajbrowseul} to $(\ref{h2goodk1})_{y,t,k}$, (strictly speaking, there is a $t$ instead $t-k$ in (\ref{trajA.5bidbis}), but this does not really matter because of $\frac{t}{2}-k\geq \frac{t}{4}$) we have
\begin{eqnarray*}
&&\ee^{-\sqrt{2\d}(\rho(y)-Y_k(y))}\ee^{-\d(t-k)} t^\frac{3}{2} (\ref{h2goodk1})_{y,t,k}
\\
&&\leq c'\ee^{-\d(t-k)}   (\rho(y)-Y_k(y))\ee^{-\sqrt{2\d}(\rho(y)-Y_k(y))}\1_{\{ \rho(y)-Y_k(y)\geq 0\}}
\\
&&\leq c''\ee^{-\d(t-k)}  \ee^{- (\rho(y)-Y_k(y))} .
\end{eqnarray*}
In addition if $x$ is $\text{good}_k$, we can ensure that 
$\rho(y)-Y_k(y)\geq 3e_k+\rho(y)-\rho(x)-D \geq 3e_k-(D+1)$ (recall that $|x-y|\leq \ee^{-5\log l}$  implies $|\rho(x)-\rho(y)|\leq 1$ as $\rho(\cdot) \in \mathcal{C}_R(l,\kappa_\d \log l, +\infty)$). Finally by combining (\ref{h2goodk1}), (\ref{bibt})  and (\ref{cardinal}), we deduce that
\begin{eqnarray}
\nonumber(\ref{h2good})_k&\leq&  c\ee^{-\d(t-k)-3e_k+D} \sum_{x,y\in \boxplus_t,\, y\in A_k(x)} \E\left(\1_{\{Y_\cdot(x)\in  \DD^{\rho(x),0}_{t},\, x\, \text{good}_k \}} \right)
\\
\nonumber&\leq&  c' \ee^{-3e_k+D}  \ee^{-\d t} \sum_{x \in \boxplus_t  }\rho(x)\ee^{-\sqrt{2\d}\rho(x)}\ee^{-\d(t-k)}\sum_{y\in \boxplus_t,\, y\in A_k(x)}1
\\
\label{b=ine}&\leq &c''\ee^{-3e_k+D}\mathtt{I}_\d(\rho).
\end{eqnarray}

\paragraph{  c)  $t/4 \leq k\leq t/2$.} Let $x,y\in \boxplus_t$ with $y\in A_k(x)$. In addition if $x$ is $\text{good}_k$, we can ensure that 
$\rho(y)-Y_k(u)\geq 3e_k+\rho(y)-\rho(x)-D  \geq 3e_k-(D+1) $ (recall that $|x-y|\leq \ee^{-\frac{t}{4}}$  implies $|\rho(x)-\rho(y)|\leq 1$ as $\rho(\cdot) \in \mathcal{C}_R(l,\kappa_\d \log l, +\infty)$). Finally by combining (\ref{h2goodk2}), (\ref{bibt}) and (\ref{cardinal}), we deduce that
\begin{eqnarray}
\nonumber(\ref{h2good})_k&\leq&  c t^\frac{3}{2}\ee^{-\d(t-k)}\ee^{\sqrt{2\d}(-3e_k+D)} \sum_{x,y\in \boxplus_t,\, y\in A_k(x)} \E\left(\1_{\{Y_\cdot(x)\in  \DD^{\rho(x),0}_{t},\, x\, \text{good}_k \}} \right)
\\
\nonumber&\leq&  c' t^\frac{3}{2}\ee^{\sqrt{2\d}(-3e_k+D)}  \ee^{-\d t} \sum_{x \in \boxplus_t  }\rho(x)\ee^{-\sqrt{2\d}\rho(x)}\ee^{-\d(t-k)}\sum_{y\in \boxplus_t,\, y\in A_k(x)}1
\\
\label{c=ine}&\leq &c'' \ee^{\sqrt{2\d}D}\ee^{ -\frac{\sqrt{2\d}}{2}3e_k } \mathtt{I}_\d(\rho).
\end{eqnarray}

\paragraph{  d)  $t/2 \leq k\leq t-1$.} Let $x,y\in \boxplus_t$ with $y\in A_k(x)$. In addition if $x$ is $\text{good}_k$, then we can ensure that 
$\rho(y)-Y_k(y)\geq 3e_k+\rho(y)-\rho(x)-D-a_t \geq 3e_k-(D+1)-a_t$. Finally by combining (\ref{h2goodk2}), (\ref{bibt}) and (\ref{cardinal}), we deduce that
\begin{eqnarray}
\nonumber(\ref{h2good})_k&\leq&  c  \ee^{-\d(t-k)}\ee^{\sqrt{2\d}(-3e_k+D)} \sum_{x,y\in \boxplus_t,\, y\in A_k(x)} \E\left(\1_{\{Y_\cdot(x)\in  \DD^{\rho(x),0}_{t},\, x\, \text{good}_k \}} \right)
\\
\nonumber&\leq&  c'  \ee^{\sqrt{2\d}(-3e_k+D)}  \ee^{-\d t} \sum_{x \in \boxplus_t  }\rho(x)\ee^{-\sqrt{2\d}\rho(x)}\ee^{-\d(t-k)}\sum_{y\in \boxplus_t,\, y\in A_k(x)}1
\\
\label{d=ine}&\leq &c'' \ee^{\sqrt{2\d}D}\ee^{ - {\sqrt{2\d}} 3e_k } \mathtt{I}_\d(\rho).
\end{eqnarray}

The terms with $e_k$ allow us to control the $\sum_{k=5\log l+1}^{t-1} (\ref{h2good})_k $. Indeed by combining (\ref{a=ine}), (\ref{b=ine}), (\ref{c=ine}) and (\ref{d=ine}) with (\ref{h2good}) we get:
\begin{equation}
\E(h^2_{good})\leq c'\mathtt{I}_\d(\rho)+ \left(c''\ee^D\sum_{k=2}^{5\log l} \frac{1}{l^\frac{\kappa_\d}{2}}+ c''\ee^{\sqrt{2\d}D} \sum_{k=5\log l+1}^{t-1} \ee^{-\frac{\sqrt{2\d}}{2} 3e_k}\right)\mathtt{I}_\d(\rho) \leq c \mathtt{I}_\d(\rho).
\end{equation}
By the Paley-Zygmund inequality, we have $ \P\left(h_{good}\geq 1\right)\geq \frac{\E\left(h_{good}\right)^2}{\E\left(h_{good}^2\right)}\geq c\mathtt{I}_\d(\rho)$. We conclude because of $h_{good}\geq 1$ implies $\exists x\in [0,R]^d,\, Y_t(x)\geq a_t+ \rho(x)-1$. \hfill$\Box$

%
%

\section{Tail of distribution of the maximum $\M_t$}
Our aim is to prove the Proposition \ref{queuedistrib}. We recall (\ref{uUu}) and (\ref{dDdD}).\\

\noindent{\it Proof of Proposition \ref{queuedistrib}.} Let $R$ and $\epsilon>0$. We want to estimate for $\rho(\cdot)  \in\mathcal{C}_R(l,\kappa_\d \log l,\log t)$,  $\P\left(\exists x\in [0,R]^\d,\, Y_t(x)\geq a_t+\rho(x)\right)$. We introduce some notations:
\begin{equation}
\label{def1}
M_{t,\rho} :=\underset{y\in [0,R]^\d}{\sup} (Y_t(y)-\rho(y)),\qquad \o_{t,\rho}:=\{y\in [0,R]^\d,\, Y_t(y)\geq a_t+\rho(y)-1\},
\end{equation}
\begin{equation}
\label{def2}
M_{t,\rho}(x, b):= \underset{y\in B(x,\ee^{b-t})}{\sup}(Y_t(y)-\rho(y)),\qquad \o_{t,\rho}(x,b):=\{y\in B(x,\ee^{b-t}),\, Y_t(y)\geq a_t+\rho(y)-1\},
\end{equation}
\begin{equation}
\label{def3}\text{and}\qquad\qquad\qquad\qquad \R_t:= [\ee^{-t/2},R-\ee^{-t/2}]^\d.\qquad\qquad \qquad\qquad
\end{equation}

\nomenclature[i3]{$ M_{t,\rho} $}{$:=\underset{y\in [0,R]^\d}{\sup}(Y_t(y)-\rho(y))$} 
\nomenclature[i4]{$ M_{t,\rho}(x,b)$}{$:=\underset{y\in B{(x,\ee^{b-t})}}{\sup}(Y_t(y)-\rho(y))$}
\nomenclature[c8]{$\R_t $}{$:=[\ee^{-\frac{t}{2}},R-\ee^{-\frac{t}{2}}]^\d $}
\nomenclature[c9]{$ \o_{t,\rho}(x,b)$}{$:=\{y\in B(x,\ee^{b-t}),\, Y_t(y)\geq a_t+\rho(y)-1\} $}
\nomenclature[c91]{$\mathfrak{O}_{t,\rho}  $}{$:=\{y\in [0,R]^\d,\, Y_t(y)\geq a_t+\rho(y)-1\} $}

For any $t>0$, because of the continuity of the function $x\mapsto Y_t(x)-\rho(x)$, the random variables $\lambda (\o_{t,\rho} ) $ and $\lambda (\o_{t,\rho}(x,b)) $ are strictly positive respectively on $\{ M_{t,\rho}\geq a_t\} $ and $ \{M_{t,\rho}(x,b)\geq a_t\} $. Therefore for any $L\geq 1$, 
\begin{eqnarray*}
\P\left( \exists x\in [0,R]^\d,\, Y_t(x)\geq a_t+\rho(x)\right)= \P\left(M_{t,\rho}\geq a_t\right)&=&\E\left(\int_{[0,R]^\d}\frac{\1_{\{ x\in \o_{t,\rho}\}}\1_{\{M_{t,\rho}\geq a_t \}}}{ \lambda (\o_{t,\rho} )}    dx\right)
\\
&=&(1)_L+(2)_L+(3),
\end{eqnarray*}
with  
\begin{eqnarray*}
(1)_L&:=&\E\left(\int_{\R_t}  \frac{\1_{\{ x\in \o_{t,\rho},\, Y_\cdot(x)\in \DD_t^{\rho(x),L}\}}\1_{\{M_{t,\rho}\geq a_t \}}}{ \lambda (\o_{t,\rho} )}dx\right),
\\
(2)_L&:=&\E\left(\int_{\R_t}  \frac{\1_{\{ x\in \o_{t,\rho},\, Y_\cdot(x)\notin \DD_t^{\rho(x),L}\}}\1_{\{M_{t,\rho}\geq a_t \}}}{ \lambda (\o_{t,\rho} )}     dx\right),
\\
(3)&:=&\E\left(\int_{[0,R]^\d-\R_t}\frac{\1_{\{ x\in \o_{t,\rho}\}}\1_{\{M_{t,\rho}\geq a_t \}}}{ \lambda (\o_{t,\rho} )}    dx\right).
\end{eqnarray*}
We shall show that $(2)_L$ and $ (3)$ are negligible, only  $(1)_L$ contributes in (\ref{eqqueuedistrib}).

Recall (\ref{dDdD}), clearly $(2)_L\leq \P\left(\exists x\in \o_{t,\rho}\cap[0,R]^\d ,\, Y_\cdot(x) \notin \DD_t^{\rho(x),L}\right)$. Via Proposition \ref{tension*} , there exist $L$ and $ l_0(L)\geq 0$ such that for any $l \geq l_0$ there exists $T\geq 0$  such that for any $\rho(\cdot) \in \mathcal{C}_{R}(l,\kappa_\d \log l,\log t)$,   
$$(2)_L\leq \epsilon \mathtt{I}_\d(\rho) .$$ 
Concerning $(3)$, decomposing $[0,R]^\d\backslash_{\R_t}$ in, at most, $ 2^d\ee^{(\d-1)\frac{t}{2}}$ cube of volume $\ee^{-\d\frac{t}{2}} $, and by the invariance by translation of $(Y_s(\cdot))_{s\geq 0} $, we have
\begin{eqnarray}
\label{depart(3)} (3)\leq \P\left(\exists x\in [0,R]^\d\backslash_{\R_t} ,\, x\in \o_{t,\rho}\right) \leq 2^\d \ee^{(\d-1)\frac{t}{2}} \P\left( \exists x\in [0,\ee^{-\frac{t}{2}}]^\d,\, Y_t(x)\geq a_t\right).
\end{eqnarray} 
Furthermore on the event $\{ \exists x\in [0,\ee^{-\frac{t}{2}}  ]^\d,\, Y_t(x)  \geq a_t \}$, we introduce ${\bf r}>0$ (${\bf r}$ is random) the biggest radius, in a similar way as in the proof of Lemma \ref{maxfx}, such that

- there exists $z_{\bf r}\in [0,\ee^{-\frac{t}{2}} ]^\d$, with $B(z_{\bf r},{\bf r})\subset [0,\ee^{-\frac{t}{2}}  ]^\d$,

- there exists $x_{\bf r}\in  B(z_{\bf r},{\bf r}) $ such that $  Y_t(x_{\bf r})  \geq a_t  $,

- for any $y\in B(z_{\bf r},{\bf r})$, $  Y_t(x)  \geq a_t-1 $.

Thus on $\{ \exists  x\in [0, \ee^{-\frac{t}{2}}  ]^\d,\,   Y_t(x)  \geq a_t \} $, by definition of ${\bf r}>0$, for any $t,\, {\bf c}>0$, (${\bf c}>0$ will be determined later) we have 
\begin{eqnarray*}
1&=&\frac{1}{S{\bf r}^\d}\int_{B(z_{\bf r},{\bf r}) }\1_{\{   Y_t(x)  \geq a_t-1\}} dy
\\
&=& \left( \1_{\{ {\bf r}\geq \frac{\ee^{- (t+{\bf c}) }}{ 4}\}}+\underset{p\geq t+{\bf c}}{\sum} \1_{\{\frac{\ee^{-(p+1)}}{4}\leq {\bf r}< \frac{\ee^{-p}}{4} \}}   \right)\frac{1}{S{\bf r}^\d}\int_{B(z_{\bf r},{\bf r}) }\1_{\{  Y_t(y)  \geq a_t-1 \}}dy.
\end{eqnarray*}
By taking the expectation we obtain that
\begin{eqnarray}
&&\nonumber \P\left( \exists x\in [0,\ee^{-\frac{t}{2}}]^\d,\, Y_t(x)\geq a_t\right)  \leq S^{-1}4^\d\ee^{\d (t+{\bf c})}\int_{[0,\ee^{-\frac{t}{2}} ]^\d}\P\left(   Y_t(y)  \geq a_t-1 \right)   dy  
\\
\label{malin(3)} && \qquad\qquad\qquad\qquad+\underset{p\geq  t +{\bf c} }{\sum} S^{-1}4^\d \ee^{\d (p+1)}\E\left(\1_{\{  {\bf r}\leq \frac{\ee^{-p}}{4}\}}\int_{B(z_{\bf r},{\bf r}) }\1_{\{  Y_t(y)  \geq a_t-1\}}dy\right).
\end{eqnarray}

Fix $p\geq t +{\bf c} $. On $\{    {\bf r}\leq \ee^{-p}/4\}$, $B(z_{\bf r},{\bf r})\neq [0,\ee^{-\frac{t}{2}} ]^\d$. So there exists $\overline{z}\in [0,\ee^{-\frac{t}{2}}  ]^\d,\,  |\overline{z}-z_{\bf r}|\leq 2{\bf r}\leq \frac{\ee^{-p}}{2}$ with $Y_t(\overline{z})\leq a_t-1  $ which implies that $  |Y_t(\overline{z})-Y_t(x_{\bf r})|\geq 1$. Thus for any $y\in B(z_{\bf r}, {\bf r}) $, by the triangular inequality we deduce that there exists $u\in [0,\ee^{-\frac{t}{2}}  ]^\d$, $|u-y|\leq \ee^{-p}$ ($u$ is either $x_{\bf r}$ or $\overline{z}$) such that $\underset{s\leq t }{\sup} |Y_s(u)-Y_s(y) |\geq \frac{1}{2}$. To summarize,  
\begin{equation}
\label{bruxelles(3)} \{    {\bf r}\leq \ee^{-p}/4\}\cap \{y\in B(z_{\bf r},{\bf r}) \}\subset \Big\{ \underset{ u\in B(y,\ee^{-p}),\, s\leq t}{\sup} |Y_s(y)-Y_s(u)|\geq \frac{1}{2}\Big\}.
 \end{equation}
According to Lemma \ref{rRhodesdec}, for any $y,u\in [0,\ee^{-\frac{t}{2}} ]^d$ such that $|y-u|\leq \ee^{-p}  $,
 \begin{eqnarray*}
(Y_s(u))_{s\leq t}&=&( P_s^y(u)+ Z_s^y(u)- \zeta_s^y(u))_{s\leq t}
\\
&=&( P_s^y(u)+ Z_s^y(u)+O(\ee^{k-p}))_{s\leq t}.
\end{eqnarray*}
Now, we choose ${\bf c}>0$ large enough such that for any $p\geq t+{\bf c} $ the $ O(\ee^{t-p})$ is smaller than $\frac{1}{2^7}$ (we stress that such $ {\bf c}$ does not depend on $ k$). Consequently for any $p\geq t+ {\bf c}$  the event in the right-hand side of (\ref{bruxelles(3)}) is included in
\begin{eqnarray*}
\{ \underset{ u\in B(y,\ee^{-p}),\, s\leq t}{\sup} |P_s^y(u)-Y_s(y)|\geq 2^{-3}\} &\cup &\{ \underset{ u\in B(y,\ee^{-p}),\, s\leq t}{\sup} |Z^y_s(u)|\geq 2^{-3}\}
\\
 =\{ w_{P_{\cdot}^y(\cdot)}( \ee^{-p},y,t)  \geq 2^{-3} \}&\cup &\{ w_{Z_{\cdot}^y(\cdot)}( \ee^{-p},y,t)  \geq 2^{-3} \},\qquad ( w_\cdot(\cdot,\cdot,\cdot) \text{ is defined in (\ref{remodulus})}) .
\end{eqnarray*}
We go back to (\ref{malin(3)}), and use the independence between $(Z_\cdot^y(u))_{u\in [0,R]^d}$ and $Y_\cdot(y)$ to deduce that there exist some constants $c,\, {\bf c}>0$ (independent of $k$) such that
\begin{eqnarray}
\label{goback(3)} && \qquad \qquad \qquad \P\left( \exists x\in [0,\ee^{-\frac{t}{2}}]^\d,\, Y_t(x)\geq a_t\right)  \leq  c
\int_{[0,\ee^{-\frac{t}{2}} ]^\d}\P\left( Y_t(y) \geq a_t-1 \right)    \Big[\ee^{\d t}+
\\
\nonumber && \underset{p\geq  t+{\bf c}}{\sum} \ee^{\d p}\P\left(     w_{Z_{\cdot}^y(\cdot)}( \ee^{-p},y,t) \geq 2^{-3} \right)\Big]+ \underset{p\geq  t+{\bf c}}{\sum} \ee^{\d p} \P\left(  Y_t(y) \geq a_t-1 ,\,   w_{P_\cdot^y(\cdot)}( \ee^{-p},y,t)   \geq 2^{-3} \right)dy .
\end{eqnarray}
Referring to the Appendix, by (\ref{eqpartieind}) in Lemma \ref{partieind} we get
\begin{eqnarray*}
\underset{p\geq  t+{\bf c}}{\sum} \ee^{\d p}\P\left(  w_{Z_\cdot^y(\cdot)}( \ee^{-p},y,t)   \geq 2^{-3} \right) &= &  \underset{p\geq  t+{\bf c}}{\sum} \ee^{\d p}\P\left( \underset{|u|\leq \ee^{-p},\,s\in[0,t]}{\sup}|Z^0_s(u)|  \geq 2^{-3} \right)
\\
&\leq &  \ee^{\d  (t+{\bf c})} \underset{p\geq  t+{\bf c}}{\sum} \ee^{\d(p- (t+{\bf c}))} c_{15} \ee^{-c_{16}2^{-6} \ee^{2(p-t)}}
\\
&\leq & c\ee^{ \d {\bf c}}\ee^{\d t} = c' \ee^{\d t}.
\end{eqnarray*}
Fix $y\in[0,\ee^{-\frac{t}{2}}  ]^d$. By Girsanov's transformation we observe that
\begin{equation}
\label{youpi(3)1}\P\left( Y_t(y) \geq a_t-1 \right)
=\E\left( {\1_{\{ B_t\geq a_t-1 \}}}{\ee^{-\sqrt{2\d}B_{t}-\d t}}\right)\leq  c\ee^{-\d t} \ee^{-\sqrt{2\d} a_t},
\end{equation}
and
\begin{eqnarray}
\nonumber \P\left(  Y_t(y) \geq a_t-1 ,\,   w_{P_\cdot^y(\cdot)}( \ee^{-p},y,t)   \geq 2^{-3} \right)  &=&  \E\Big( \ee^{-\sqrt{2\d} B_{t}-\d t} \1_{\{ B_t \geq a_t-1 ,\,    \underset{|u|\leq  \ee^{-p}}{\sup}\big|\int_{0}^{.} \mathtt{g}(\ee^s u)  dB_s\big|_{t} \geq  2^{-3}\}}   \Big)
\\
\label{youpi(3)2}&\leq &\ee^{-\d t} \ee^{-\sqrt{2\d} a_t}\P\Big( \underset{|u|\leq  \ee^{-p}}{\sup}\big|\int_{0}^{.} \mathtt{g}(\ee^s u)  dB_s\big|_{t} \geq  2^{-3} \Big) .
\end{eqnarray}
By (\ref{eqmajorAB3}) in Lemma \ref{majorAB}, for any $p\geq t+{\bf c} $, we have
\begin{equation}
\label{youpi(3)}   \P\left(\underset{|u|\leq  \ee^{-p}}{\sup}\big|\int_{0}^{.} \mathtt{g}(\ee^s u)  dB_s\big|_{t} \geq  2^{-3}   \right)= \P(A_{p,t,2^{-3}})\leq 
 c_{20}\exp({-c_{19}  2^{-6}\ee^{2(p-t)}}).
\end{equation}
Go back to (\ref{goback(3)}) combining (\ref{youpi(3)1}), (\ref{youpi(3)2}) and (\ref{youpi(3)}), we obtain that
\begin{eqnarray*}
\P\left( \exists x\in [0,\ee^{-\frac{t}{2}}]^\d,\, Y_t(x)\geq a_t\right)  &\leq&  c \ee^{-\sqrt{2\d} a_t} \int_{[0,\ee^{-\frac{t}{2}}]^\d} [1+\sum_{p\geq t+ {\bf c} }  \ee^{\d p}  c_{20}\exp({-c_{19}  2^{-6}\ee^{2(p-t)}}) ]dy
\\
&\leq& c\ee^{-\sqrt{2\d} a_t} \ee^{-\frac{\d t}{2}}.
\end{eqnarray*}
Finally with (\ref{depart(3)}) for $t>0$ large enough, it stems that
\begin{eqnarray*}   (3) \leq 2^\d\ee^{(\d-1)\frac{t}{2}}c\ee^{-\sqrt{2\d} a_t} \ee^{-\frac{\d t}{2}}\leq \epsilon \mathtt{I}_\d.
\end{eqnarray*} 
Therefore we can fix $L>0$, such that there exist $l>0$ and $T>0$ satisfying: $\forall\, t\geq T,\, \rho(\cdot)  \in \mathcal{C}_R(l,\kappa_\d \log l,\log t)$, 
\begin{eqnarray}
\label{5.1?}|\P(M_{t,\rho}\geq a_t)-(1)_L|\leq \epsilon \mathtt{I}_\d(\rho).
\end{eqnarray}
The previous inequality just express that with an overwhelming probability for any $x\in [0,R]^\d  $,  $ Y_t(x) \geq a_t +\rho(x) $ is equivalent to $ Y_\cdot(x)\in \DD_t^{\rho(x),L} $. We will take advantage of this fact to know the spatial distribution of extremal particles.
\\

For any $t>b\geq 0$, let us introduce:
\begin{equation}
\label{defXi}
\Xi_{\rho,t}(b,x)=\{ \exists y\in [0,R]^\d,\, |y-x|\geq \ee^{b-t},\, Y_t(y)\geq a_t+\rho(y)-1\},
\end{equation} \nomenclature[h2]{$\Xi_{\rho,t}(b,x) $}{$:=\{ \exists y\in [0,R]^\d,\, \mid y-x \mid \geq \ee^{b-t},\, Y_t(y)\geq a_t+\rho(y)-1\} $} On $   \Xi_{\rho,t}(b,x)^c $, $\frac{\1_{\{M_{t,\rho}\geq a_t\}}}{\lambda(\o_{t,\rho})}=\frac{\1_{\{M_{t,\rho}(x,b)\geq a_t\}}}{\lambda(\o_{t,\rho}(x,b))}$ therefore we obtain
\begin{equation}
\label{inter1} (1)_{L}=(1)_{L,b}+(2)_{L,b}-(3)_{L,b},\qquad \forall b\geq 0,
\end{equation}
with 
\begin{eqnarray}
\label{(1)}(1)_{L,b}&:=&\E\left(\int_{\R_t}\frac{\1_{\{Y_\cdot(x) \in \DD_t^{\rho(x),L}\}}\1_{\{M_{t,\rho}(x,b)\geq a_t \}} }{ \lambda (\o_{t,\rho}(x,b) )}dx\right),
\\
\label{(2)}(2)_{L,b}&:=& \E\left(\int_{\R_t} \1_{\Xi_{\rho,t}(b,x)}  \frac{\1_{\{  Y_\cdot(x)\in \DD_t^{\rho(x),L}\}}  \1_{\{M_{t,\rho}\geq a_t \}}}{ \lambda (\o_{t,\rho})} dx\right),
\\
\label{(3)}(3)_{L,b}&:=&\E\left(\int_{\R_t} \1_{\Xi_{\rho,t}(b,x)}  \frac{\1_{\{  Y_\cdot(x)\in \DD_t^{\rho(x),L}\}}\1_{\{M_{t,\rho}(x,b)\geq a_t  \}}}{ \lambda (\o_{t,\rho}(x,b) )} dx\right).
\end{eqnarray}
We shall show, via two lemmas, that $(2)_{L,b}$ and $(3)_{L,b}$ are negligible. 
\begin{Lemma}
\label{avantXi}
There exists $c_8>0$ such that for any $L,\,  t\geq b>1$,\nomenclature[i5]{ $ {\bf r}_t(x)$ }{$ :=\sup\{ r>0,\, w_{Y_\cdot(\cdot)}(r,x ,t)\leq \frac{1}{4}\}  \wedge  \ee^{-t} $}
\begin{eqnarray}
(2)_{L,b}+(3)_{L,b}\leq c_8\int_{[0,R]^\d}\E\left(\frac{\1_{\{Y_\cdot(x) \in \DDi_t^{\rho(x),L}\}}}{{\bf r}_t(x)^\d}\1_{\{ \Xi_{\rho,t}(b-\log 2,x)\}}\right)dx, \text{ with}
\\
\label{defrx}{\bf r}_t(x):=\sup\{ r>0,\, w_{Y_\cdot(\cdot)}(r,x ,t)\leq \frac{1}{4}\}  \wedge  \ee^{-t}.
\end{eqnarray}
\end{Lemma} 
\noindent{\it Proof of Lemma \ref{avantXi}.} Fix $x^*\in \R_t$ and observe that 
\begin{eqnarray*}
  \int_{B{(x^*, \frac{1}{4}\ee^{b-t})}}\1_{\Xi_{\rho,t}(b,x)}  \left[\frac{\1_{\{ Y_\cdot(x) \in  \DD_t^{\rho(x),L} ,\,  M_{t,\rho}(x,b)\geq a_t\}}}{ \lambda (\o_{t,\rho}(x,b) )} + \frac{\1_{\{ Y_\cdot(x)  \in \DD_t^{\rho(x),L},\, M_{t,\rho}\geq a_t\}}}{ \lambda (\o_{t,\rho} )}\right] dx 
  \\
  \leq 2\1_{\{\exists  x\in B{(x^*, \frac{1}{4}\ee^{b-t})},\, Y_\cdot(x)\in \DD_t^{\rho(x),L}  \}}\1_{\{ \Xi_{\rho,t}(b-\log \frac{4}{3},x^*) \}}.
\end{eqnarray*}
By continuity of $y\mapsto (Y_s(y)-\rho(y))_{s\leq t}$, if $x\in B(x^*,\frac{1}{4}\ee^{b-t})$ such that $Y_\cdot(x)\in \DD_t^{\rho(x),L}$ ($x$ satisfies the strong condition) then there exist $r>0$ and $  x_r\in B(x^*,\frac{1}{4}\ee^{b-t}) $ such that: $x \in B(x_r,r)$; $  B(x_r,r)\subset B(x^*,\frac{1}{4}\ee^{b-t}) $; and for any  $y \in B(x_r,r)$, $Y_\cdot(y) \in \DDi_t^{\rho(y),L}$ ($y$ satisfies the weak condition).

\begin{figure}
\centering
\caption{}
\label{Tux2}
\includegraphics[interpolate=true,width=9.5cm,height=9.5cm]{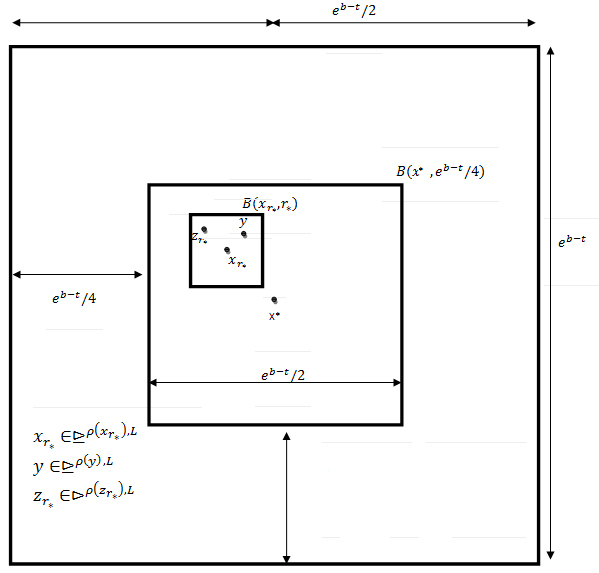} 
\end{figure}

Thus on the set $\{ \exists x\in  B(x^*, \frac{1}{4}\ee^{b-t} ),\, Y_\cdot(x)\in  \DD_t^{\rho(x),L}\}$, there exists ${\bf r}_*>0$ (see figure \ref{Tux2} pp 38) which is the biggest radius such that: 

- there exists $ x_{\bf r_*}\in [0,R]^\d$ with $B(x_{\bf r_*}, {\bf r_*})\subset B(x^*,\frac{1}{4}\ee^{b-t})$,

- there exists $ z_{\bf r_*}\in B(x_{\bf r_*}, {\bf r_*})$ with $Y_\cdot(z_{\bf r_*})\in \DD_{t}^{\rho(z_{\bf r_*}),L}$,

- for any $ y\in B(x_{\bf r_*}, {\bf r_*})$, $Y_\cdot(y)\in \DDi_{t}^{\rho(y),L}$.

By definition, $\frac{1}{S{\bf r}_*^\d}\int_{B(x_{{\bf r}_*},{\bf r}_*) }\1_{\{Y_\cdot(y)\in \DDi_t^{\rho(y),L}  \}} dy=1$ on $\{ \exists x\in  B(x^*, \frac{1}{4}\ee^{b-t}),\, Y_\cdot(x)\in \DD_t^{\rho(x),L}  \}$, so we can affirm that 
\begin{eqnarray*}
&&\int_{B{(x^*, \frac{1}{4}\ee^{b-t})}}\1_{\Xi_{\rho,t}(b,x)}  \left[\frac{\1_{\{ Y_\cdot(x)  \in \DD_t^{\rho(x),L},\, M_{t,\rho}(x,b)\geq a_t \}}}{ \lambda (\o_{t,\rho}(x,b) )} + \frac{\1_{\{ Y_\cdot(x)  \in \DD_t^{\rho(x),L},\, M_{t,\rho}\geq a_t \}}}{ \lambda (\o_{t,\rho} )}\right]   dx
\\
& &\leq   \frac{2}{S{{\bf r}}_*^\d}\int_{B(x_{\bf r_*},{\bf r_*}) }\1_{\{ Y_\cdot(y)\in \DDi_t^{\rho(y),L} \}}\1_{\Xi_{\rho,t}(b-\log {\frac{4}{3}},x^*) } dy
\\
&&\leq  \int_{B(x_{{\bf r}_*},{\bf r}_*) }\frac{c}{{\bf r}_*^\d}\1_{\{ Y_\cdot(y)\in \DDi_t^{\rho(y),L} \}}\1_{\Xi_{\rho,t}(b-\log {2},y) } dy.
\end{eqnarray*} 
Furthermore on $\{y\in B(x_{{\bf r}_*},{\bf r}_*) \}$ there are two possible options:

- Either there exists  $\overline{x}_{{\bf r}_*}\in [0,R]^\d$ such that $|z_{{\bf r}_*}-\overline{x}_{{\bf r}_*}|\leq 2 {\bf r}_*$ and $Y_\cdot(\overline{x}_{{\bf r}_*})\notin  \DD_t^{\rho(\overline{x}_{{\bf r}_*}),L}$ which implies that $\sup_{s\leq t} \big|Y_s(\overline{x}_{{\bf r}_*})-Y_s(z_{{\bf r}_*})+\rho(\overline{x}_{{\bf r}_*})-\rho(z_{{\bf r}_*})\big|\geq 1$. As $  \rho(\cdot) \in \mathcal{C}_R(l,\kappa_\d \log l,\log t)$ (implying $\sup_{z \in B(y,4{\bf r}_* )}|\rho(y)-\rho(x)|  \leq \frac{1}{4}$)  by the triangular inequality, for any $y\in B(x_{{\bf r}_*},{\bf r}_*) $ we have finally $  w_{Y_\cdot(\cdot) }(4{\bf r}_*,y,t)\geq \frac{1}{4}$.

- Or $B(x_{{\bf r}_*}, {\bf r}_*)=B(x^*,\frac{1}{4}\ee^{b-t})$ and thus ${\bf r_*}=\frac{1}{4}\ee^{b-t}$,

Nevertheless in the both cases we have on $ \{ y\in B(x_{{\bf r}_*},{\bf r}_*) \}  $
\begin{eqnarray*}
4{\bf r}_*\geq {\bf r}_t(y):= \sup\{ r>0,\, w_{Y_{\cdot}(\cdot)}(r,y,t)\leq \frac{1}{4}\}\wedge \ee^{-t}  .
\end{eqnarray*}
We deduce that
\begin{eqnarray*}
&&\int_{B{(x^*, \frac{1}{4}\ee^{b-t})}}\1_{\Xi_{\rho,t}(b,x)} \left[\frac{\1_{\{ Y_\cdot(x)  \in \DD_t^{\rho(x),L},\,  M_{t,\rho}(x,b)\geq a_t \}}}{ \lambda (\o_{t,\rho}(x,b) )} + \frac{\1_{\{ Y_\cdot(x) \in  \DD_t^{\rho(x),L},\,  M_{t,\rho}\geq a_t \}}}{ \lambda (\o_{t,\rho} )}\right]   dx 
\\
&&\leq   4^\d c \int_{B(x_{\bf r_*},{\bf r_*}) }\frac{\1_{\{Y_\cdot(y)\in \DDi_t^{\rho(y),L} \}}}{{\bf r}_t^\d(y)}\1_{\Xi_{\rho,t}(b-\log 2,y) } dy
\\
&&\leq   c'\int_{ B{(x^*, \frac{1}{4}\ee^{b-t})}}\frac{\1_{\{Y_\cdot(y)\in \DDi_t^{\rho(y),L} \}}}{{\bf r}_t^\d(y)}\1_{\Xi_{\rho,t}(b- \log 2,y) } dy.
\end{eqnarray*}
This inequality is true for any $x^*\in \R_t$, moreover we can find $\N\ni  m\leq c\ee^{d(t-b)}$, and $(x_i)_{i\leq m}$ a collection of $\R_t$ such that:

(i) $\R_t \subset \underset{1\leq i\leq m}{\bigcup}B(x_i,\frac{\ee^{b-t}}{4})\subset [0,R]^\d  $,

(ii) for any $(i_1,...,i_{d+2})\in \{1,... , m\}$ distinct, $\underset{j=1}{\overset{d+2}{\bigcap}}B(x_i,\frac{\ee^{b-t}}{4})=\emptyset$. 

Finally there exists $c>0$ independent of $L,\,l$ or $t\geq b> 1 $ such that
\begin{eqnarray*}
\int_{\R_t} \1_{\Xi_{\rho,t}(b,x)}  \left[\frac{\1_{\{ Y_\cdot(x)  \in \DD_t^{\rho(x),L},\, M_{t,\rho}(x,b)\geq a_t \}}}{ \lambda (\o_{t,\rho}(x,b) )} + \frac{\1_{\{ Y_\cdot(x) \in   \DD_t^{\rho(x),L},\, M_{t,\rho}\geq a_t \}}}{ \lambda (\o_{t,\rho} )}\right] dx\leq
\\
 c \int_{[0,R]^\d }\frac{\1_{\{Y_\cdot(y)\in \DDi_t^{\rho(y),L}\}}}{{\bf r}_t^\d(y)}\1_{\Xi_{\rho,t}(b- \log{2},y) } dy,\qquad \text{a.s}.
\end{eqnarray*}
Lemma \ref{avantXi} follows easily. \hfill$\Box$

The proof of the following lemma is postponed at the end of this section:
\begin{Lemma}
\label{htoui}
 Let $R$, $L$ be two constants fixed. For any $\epsilon>0$ we can find $b_0,\, l_0 >1$ large enough such that for any $l\geq l_0$, $b> b_0$,  $\exists T>0$ so that the following inequality holds 
\begin{equation}
\label{eqhtoui}
(4)_{L,b}:= \int_{[0,R]^\d}\E\left(\frac{\1_{\{ Y_\cdot(x) \in \DDi_t^{\rho(x),L}\}}}{{\bf r_t}(x)^\d}\1_{ \Xi_{\rho,t}(b,x)}\right)dx \leq \epsilon   \mathtt{I}_{\d}(\rho),
\end{equation}
provided that $ t\geq T$, $\rho(\cdot) \in \mathcal{C}_R(l,\kappa_\d \log l, \log t)$.
\end{Lemma}
\paragraph{Remark:}  This lemma gives a description in ``cluster" for the repartition of the extremal particles in $[0,R]^\d $. About this question, see also \cite{AZi12} for a slightly different model. 
\\


Assuming this Lemma, combining (\ref{inter1}) and (\ref{eqhtoui}) we can fix $b(L),\, l_0(L) >0$, such that for any   $l>l_0$ there exists $T>0$ such that $\forall\, t\geq T,\, \rho(\cdot) \in \mathcal{C}_R(l,\kappa_\d \log l,\log t)$, 
\begin{equation}
\label{inter2}
\left|\P(M_{t,\rho}\geq a_t)-(1)_{L,b}\right| \leq  2\epsilon \mathtt{I}_{\d}(\rho).
\end{equation}
Therefore we can  restrict our study to $(1)_{L,b}$. The Markov property at time $t_b=t-b$ and the invariance by translation of our model give
\begin{eqnarray}
\nonumber (1)_{L,b}&=&\E\left(\int_{\R_t}\frac{\1_{\{  Y_\cdot(x) \in \DD_t^{\rho(x),L} \}} \1_{\{ M_{t,\rho}(x,b)\geq a_t\}}  }{ \lambda ( \o_{t,\rho}(x,b))}dx\right)
\\
\label{tito}&=&   \int_{\R_t} \E\left(\1_{\{ \overline{Y}_{t_b}(x)\leq \rho(x),\, \overline{Y}_{[\frac{t}{2},t_b]}(x)\leq a_t+\rho(x)+L\}} (\ref{tito})^{b,L}_{x,t}\right)dx,
\end{eqnarray}
%
with 
\begin{eqnarray}
 \nonumber (\ref{tito})^{b,L}_{x,t}&:=&\E\left( \frac{\1_{\{ \overline{Y}^{(t_b)}_b(0)+z_{(\ref{tito})}\leq 0,\, Y^{(t_b)}_b(0)\geq -L-1,\, \exists y\in B(0,\ee^{b-t}),\, Y^{(t_b)}_{b}(y)+ z_{(\ref{tito})}\geq -L- g_{(\ref{tito})} (y) \}}}{\lambda_{B(0,\ee^{b-t}) }(\{y: Y^{(t_b)}_b(y)+ z_{(\ref{tito})}\geq -L-1-g_{(\ref{tito})} (y)\} )}\right) , 
\\ 
\nonumber\text{with },\,\,   \quad   && g_{(\ref{tito})} (y)=Y_{t_b}(x+y)-Y_{t_b}(x)-(\rho(x+y)-\rho(x)),
\\ 
\nonumber && z_{(\ref{tito})} = Y_{t_b}(x)-a_t -\rho(x)-L .
\end{eqnarray}
\nomenclature[i51]{$\rho_x(.)$}{$:=\rho(x+.)-\rho(x)$}  In the following we will denote
\begin{equation}
\forall x\in [0,R]^\d,\, \rho_x(.):=\rho(x+.)-\rho(x) .
 \end{equation} 
According to the scaling property $\left(Y_s^{(t_{b})}(y)\right)_{s\leq b,\,y\in B(0,\ee^{b-t})}\overset{(d)}{=} (Y_s(y\ee^{t-b}))_{s\leq b,\,y\in B(0,\ee^{b-t})}$, thus we can rewrite $  (\ref{tito})_{x,t}^{b,L} $ as
\begin{eqnarray*}
\ee^{\d t_b}\E_{z_{(\ref{tito})}}\left( \frac{\1_{\{ \overline{Y}_b(0)\leq 0,\, Y_b(0)\geq -L-1\}} \1_{\{\exists y\in B(0,1),\, Y_{b}(y)\geq -L- g_{(\ref{tito})}(y\ee^{b-t}) \}}}{\lambda_{B(0,1) }(\{y: Y_b(y)\geq -L-1-g_{(\ref{tito})}(y\ee^{b-t})\} )}\right).
\end{eqnarray*}
In addition Lemma \ref{rRhodesdec}  and the Girsanov's transformation lead to 
\begin{eqnarray}
\nonumber (1)_{L,b}&=&\int_{\R_t} \E\left(\ee^{\sqrt{2\d}Y_{t_b}(x)+\d t_b}\1_{\{ \overline{Y}_{t_b}(x)\leq \rho(x),\, \overline{Y}_{[\frac{t}{2},t_b]}(x)\leq a_t+\rho(x)+L\}}   \ee^{-\sqrt{2\d}Y_{t_b}(x)-\d t_b}\ee^{\d t_b}  (\ref{tito})^{b,L}_{x,t} \right)dx
\\
\label{egal(1)L} &=&\int_{\R_t} {\ee^{-\sqrt{2\d}\rho(x)}}{t^\frac{3}{2}}\E_{-\rho(x)}\left(\1_{\{ \overline{B}_{t_b}\leq 0,\, \overline{B}_{[\frac{t}{2},t_b]}\leq a_t+L\}}F_{L,b}\left(B_{t_b}-a_t -L,\mathfrak{G}_{t,b}^{\rho_x}\right)\right)dx,
\end{eqnarray}
with 

- as before $B$ a standard Brownian motion, 

- for $ g\in \mathcal{C}(B{(0,\ee^b)},\r),\, z\in \r  $ ,\nomenclature[i6]{$ F_{L,b}(z,g)$}{: defined in (\ref{defFLb})}  \nomenclature[i7]{$ \mathfrak{G}_{t,b}^{\Psi} $}{: defined in (\ref{defG})}
\begin{equation}
\label{defFLb}F_{L,b}(z,g):=\ee^{-\sqrt{2\d}(z+L)}\E_{z}\left( \frac{\1_{\{ \overline{Y}_b(0)\leq 0,\, Y_b(0)\geq -L-1\}} \1_{\{\exists y\in B(0,1),\, Y_{b}(y)\geq -L- g(y\ee^{b}) \}}}{\lambda_{B(0,1) }(\{y: Y_b(y)\geq -L-1-g(y\ee^b)\} )}\right),  
\end{equation}

 - for any $\Psi \in \mathcal{C}_R(B(0,\ee^b),\r)$, 
 \begin{equation}
\label{defG} \mathfrak{G}_{t,b}^\Psi:B{(0,\ee^{b})}\ni y\mapsto   -\int_{0}^{t_b}\mathtt{g}(\ee^{s-t} y)dB_{s} - \zeta_t(y\ee^{-t})+ Z_{t_b}^{0}(y\ee^{-t})-\Psi(y\ee^{-t}).
 \end{equation}
For $\Psi=0$ we denote $\mathfrak{G}_{t,b}^0=\mathfrak{G}_{t,b}$. In passing we take the opportunity to define for any $\sigma \in [0,t_b]$, \nomenclature[i8]{$ \mathfrak{G}_{t,b,\sigma}  $}{: defined in (\ref{defGsigma})}
\begin{equation}
\label{defGsigma}
\mathfrak{G}_{t,b,\sigma} :B{(0,\ee^{b})}\ni y\mapsto   -\int_{t_b-\sigma}^{t_b}\mathtt{g}(\ee^{s-t} y)dB_{s} -  \zeta_{t_b}(y\ee^{-t})+ Z_{t_b}^{0}(y\ee^{-t}).
\end{equation}

Notice that $Z^0_{t_b}(\cdot)  $ is a centred Gaussian process, independent of $B$, which have the covariances as in (\ref{covgeneralZ}). Furthermore by Proposition \ref{Weakconv} for  any $b>0$, the  Gaussian process $B{(0,\ee^b)}\ni y\mapsto Z_{t_b}^0({y\ee^{-t}})-\zeta_{t_b}(y\ee^{-t})$, converges in law to $y\mapsto Z(y\ee^{-b})-\zeta(y\ee^{-b})$.
\\

 Now we want to get (via a renewal theorem): for any $L,b>0$, uniformly in $x\in \R_t $,
$$\E_{-\rho(x)}\left(\1_{\{ \overline{B}_{t_b}\leq 0,\, \overline{B}_{[\frac{t}{2},t_b]}\leq a_t+L\}}F_{L,b}\left(B_{t_b}-a_t-L,\mathfrak{G}_{t,b}^{\rho_x}\right)\right)\sim \frac{C^*}{t^\frac{3}{2} }\rho(x).$$
We stress that $C^*$ must not depend on $x$ or $\rho$. To obtain this result, we need yet two steps:

-Study the regularity of $F_{L,b}$ (Lemma \ref{controlF})

- Use this regularity to apply a renewal theorem (Theorem \ref{RESUME})

%

 \begin{Definition}
 \label{defbregular}
 A continuous function $F:\r\times \mathcal{C}(B(0,e^b),\r)\to \r^+$ is ''{\bf $b$ regular}" if there exist two functions $h:\r\to \r_+$ and $F^*: \mathcal{C}(B(0,e^b))\to \r^+$ satisfying 
 
 (i) 
 \begin{equation}
 \label{et2rajout}
 \underset{z\in \r }{\sup}\, h(z)<+\infty ,\quad  \text{and } h(z)\underset{z\to -\infty}{=}O(e^z).
 \end{equation}
 
 (ii) There exists $c>0$ such that for any $ \delta\in (0,1)$, $g\in \mathcal{C}(B{(0,e^b)},\r)$ with $w^{(0,1)}_{g(\cdot e^b)}(\delta)\leq \frac{1}{4} $, 
 \begin{equation}
 \label{et1}
  F^*(g)\leq c{\delta^{-10}}.
 \end{equation}
 
 (iii) For any $z\in \r,\, g\in \mathcal{C}(B(0,e^b),\r)$,  $F(z,g)\leq h(z)  F^*(g) $.

 (iv) There exists $c>0$ such that for any $z\in \r,\, g_1,g_2\in \mathcal{C}(B(0,e^b),\r)$ with  $||g_1-g_2||_\infty\leq \frac{1}{8} $,   
 \begin{equation}
 \label{et3bibi}
 |F(z,g_1)-F(z,g_2)|\leq c ||g_1-g_2||_\infty^\frac{1}{4} h(z) F^*(g_1).
 \end{equation}
 \end{Definition}

 \begin{Definition}
  For any $M\geq 0$ and $F$ a function $b$ {\bf regular}, we define
 \begin{equation}
 \label{defMreload}F^{(M)}(x,g):=(F(x,g)\wedge M)  \1_{\{x\geq -M\}}.
 \end{equation}
 \end{Definition}

The proof of the following two results are postponed to the next sections. 
\begin{Lemma}[Control of $F_{L,b}$]
\label{controlF} For any $L,b>0$ the function $F_{L,b}$ defined in (\ref{defFLb}) is $b$ {\bf regular}.
\end{Lemma}

\nomenclature[i91]{$F_{L,b,M}(x,g) $}{$:= (F_{L,b}(x,g)\wedge M)  \1_{\{x\geq -M\}} $} \nomenclature[i92]{$ \overline{F}_{L,b,M} $}{$:= F_{L,b}-F_{L,b,M} $}

For any $\gamma\in \r$, let $  T_{\gamma}:= \inf\{s\geq 0,\, B_s=\gamma \}$. Let $(\mathtt{R}_s)_{s\geq0}$ be a three dimensional Bessel process starting from $0$.
 \begin{theorem}
 \label{RESUME}
 Let $b>0$ and $F:\r\times \mathcal{C}(B(0,e^b),\r)\to \r^+$ be a function {\bf $b$ regular}. For any $\epsilon>0$,  there exists $M,\sigma, l,  T >0$ such that for any $t\geq T$, $\rho(\cdot) \in \mathcal{C}_R(l,\kappa_\d \log l,\log t) $, $\mathtt{z}\in [1,\log t)^{30}$,
 \begin{equation}
 \label{eqRESUME}
 \left| \int_{\R_t} t^\frac{3}{2}e^{-\sqrt{2\d}\rho(x)} \E_{-\rho(x)} \left(\1_{\{  \overline{B}_{t_b}\leq 0,\,  \overline{B}_{[\frac{t}{2},t_b]}\leq -\mathtt{z} \}}  F\left(B_{t_b}+\mathtt{z},\mathfrak{G}^{\rho_x}_{t,b}\right)\right)dx- C_{M,\sigma}(F) I_\d(\rho)\right|\leq \epsilon I_\d(\rho).
 \end{equation}
 with
 \begin{eqnarray}
 \nonumber C_{M,\sigma}(F):= \sqrt{\frac{2}{\pi}}\int_{0}^M\int_0^u\E\Big(  F^{(M)}\Big(-u, y\mapsto Z (ye^{-b})-\zeta(ye^{-b}) \qquad\qquad\qquad\qquad\qquad
 \\
 \label{defCLBMS}  -\int_{0}^{T_{-\gamma} \wedge \sigma }(1-\mathtt{k}(\ee^{-s}ye^{-b})) dB_s - \int_{T_{-\gamma}\wedge \sigma}^\sigma (1-\mathtt{k}(\ee^{-s}ye^{-b}) )d \mathtt{R}_{s-T_{-\gamma}}\Big)\Big)d\gamma du.
 \end{eqnarray}
 \end{theorem}
%
%
 Assuming Lemma \ref{controlF} and Theorem \ref{RESUME}, we are in position to {\bf end the proof of the Theorem \ref{queuedistrib}}. Indeed combining (\ref{inter2}), (\ref{egal(1)L}), Lemma \ref{controlF} and Theorem \ref{RESUME} (applied with $\mathtt{z}=-a_t-L$ and $F=F_{L,b}$) we deduce that:  {\it $\forall \epsilon>0$ there exist $L,\, b,\,M,\, \sigma>0$ such that for $l,T>0$ large enough we have : {\bf for any $t\geq T,\, \rho(\cdot)\in \mathcal{C}_R(l,\kappa_\d \log l,\, \log t)$,
\begin{equation}
\label{eqinterfinal}
| \P(M_{t,\rho}\geq a_t)-C_{M,\sigma}(F_{L,b})\mathtt{I}_\d(\rho)|\leq \epsilon \mathtt{I}_\d(\rho).
\end{equation}}}

In addition by Proposition \ref{tension}: {\it There exist $c_1,\, c_2>0$ and  $l,\, T>0$ large such  that: {\bf  for any $t\geq T,\, \rho(\cdot) \in \mathcal{C}_R(l,\kappa_\d\log l,\log t)$, 
\begin{equation}
\label{utillower} c_1 \mathtt{I}_\d(\rho)  \leq \P(M_{t,\rho}\geq a_t)   \leq c_2{\mathtt{I}_\d(\rho)}.
\end{equation}}}
For any $ n\in \N^*$, let $(L,b,M,\sigma)_n$ such that (\ref{eqinterfinal}) is true with $\epsilon=\frac{1}{n}$. Clearly $C_n:=C_{M_n,\sigma_n}(F_{L_n,b_n})\in [\frac{c_1}{2}, 2c_2]$ for any $n\in \N$ large enough. Let $\phi:\N\to \N$ strictly increasing such that $C_{\phi(n)}\to C^*\in [c_1/2,2c_2]$. 

Now we fix $\epsilon>0$. Let $N_0>0$ such that for any $n\geq N_0$, $|C_{\phi(n)}-C^*|\leq \epsilon$. Then we choose $N_1> N_0$ 
such that $n\geq N_1$ implies $\frac{1}{\phi(n)}\leq \epsilon$. {\bf Finally there exist (according to (\ref{eqinterfinal})) $l(N_1),T(N_1)>0$ such that for any $t\geq T,\, \rho (\cdot)\in \mathcal{C}_R(l,\kappa_\d \log l,\log t)$,
\begin{eqnarray*}
|  \P(M_{t,\rho}\geq a_t)-C^* \mathtt{I}_\d(\rho)|\leq \epsilon \mathtt{I}_\d(\rho).
\end{eqnarray*}}
This completes the proof of Theorem \ref{queuedistrib}.
\hfill$\Box$
\\

In the next two subsections we shall prove Lemmas \ref{htoui} and \ref{controlF}, then in Section 6 we will prove Theorem \ref{RESUME}.

\subsection{Proof of Lemma \ref{htoui}}
This important Lemma gives the cluster representation for the extremal particles. The notion of "good particles", defined in (\ref{Aket}) and studied in section D is essential for its proof.
\\

\noindent{\it  Proof of Lemma \ref{htoui}.} Let $R,\, \epsilon,\, L>0$. Recall the definition of $\Xi_{\rho,t}(b,x)$ in (\ref{defXi}), we want to show that there exist  $b_0,\, l_0 $ such that for any   $b>b_0,\, l>l_0,\, \exists T>0$ so that the following inequality holds
\begin{eqnarray*}
(4)_{L,b}:= \int_{[0,R]^\d}\E\Big(\frac{\1_{\{ Y_\cdot(x)\in \DDi_t^{\rho(x),L}\}}}{{\bf r_t}(x)^\d}\1_{ \Xi_{\rho,t}(b,x)}\Big)dx \leq \epsilon   \int_{[0,R]^\d}\rho(x) \ee^{-\sqrt{2\d}\rho(x)}dx ,
\end{eqnarray*}
provided that $t>T$, $\rho(\cdot) \in \mathcal{C}_R(l,\kappa_\d \log l,\log t)$.

 Let $A^*:=B{(0,\ee^2)}\backslash_{ B(0,\ee)}$. Recall  (\ref{defek}), (\ref{Aet}) and (\ref{Aket}) for the definitions of respectively  $(e_s)_{s\leq t}$, , $A_k(u)$,  $u$  is  $L-$good$_k$ and $u$  is  $L-$good. By Lemma \ref{C.1} we choose $D(L,\epsilon), l_0(L,\epsilon)$ large enough such that for any $l\geq l_0,\, \exists T>0$ so that the following inequality holds
\begin{eqnarray*}
 \int_{[0,R]^\d}\E\left(\frac{\1_{\{ Y_\cdot(x)\in \DDi_t^{\rho(x),L}\}}}{{\bf r_t}(x)^\d}\1_{\{  x\text{ not}  L-\text{good} \}}\right)dx \leq  \epsilon \mathtt{I}_\d(\rho),
\end{eqnarray*}
provided that $t\geq T$, $\rho(\cdot)\in \mathcal{C}_R(l, \kappa_\d\log l,\log t)$. So we can restrict our study to 
\begin{eqnarray*}
\circledast_{L,b}:=   \int_{[0,R]^\d}\E\left(\frac{\1_{\{ Y_\cdot(x)\in \DDi_t^{\rho(x),L}\}}}{{\bf r_t}(x)^\d}\1_{ \Xi_{\rho,t}(b,x)}\1_{\{ x \, L-\text{good} \}}\right)dx.
\end{eqnarray*}
 Without loss of generality we can always assume that $t-b\in \N $, so the subsets $(A_k(y))_{1\leq k\leq t-b+1}$  form a partition of $\{ u\in [0,R]^\d,\, |y-u|\geq \ee^{b-t}\}$, therefore
\begin{eqnarray*}
\circledast_{L,b}&\leq &\underset{k=1}{\overset{  t_b+1 }{\sum}}  \int_{[0,R]^\d}\E\left(\frac{\1_{\{ Y_\cdot(x)\in \DDi_t^{\rho(x),L},\, x \, L-\text{good}\}}}{{\bf r_t}(x)^\d}  \1_{\{ \exists u\in A_k(x),\, Y_t(u)\geq a_t+\rho(u)-1\}}\right)dx .
\end{eqnarray*}
As $\mathtt{k}$ is continuous with support included in $B(0,1) $, for any $k\leq t-b+1$,  the process $( Y^{(k)}_{s}(u) )_{{\tiny\begin{matrix} s\leq t-k
\\
u\in A_k(x)
\end{matrix} }} $ is independent of 
 \begin{eqnarray*}
\mathcal{G}_k(x):=\sigma\left(Y_{s}(u),\, s\leq k,\, u\in A_k(y),\,\,\, Y_s(y),\, s\in \r_+,\, y\in B(x,\ee^{-t})   \right)
 \end{eqnarray*}
According to the definition (\ref{defrx}), clearly ${\bf r}_t(x) $ is measurable with respect to $\mathcal{G}_k(x)$. Then by the Markov property at time $k$,
\begin{eqnarray*}
\circledast_{L,b} &\leq& \underset{k=1}{\overset{  t_b+1 }{\sum}} \int_{[0,R]^\d}\E\left(\frac{\1_{\{ Y_\cdot(x)\in \DDi_t^{\rho(x),L},\, x \, \text{is } L-\text{good}_k\}}}{{\bf r_t}(x)^\d}\times \right.
\\
&&\qquad\qquad\qquad\qquad\quad \left. \P\left( \exists u\in  A_k(x) ,\, Y_{t-k}^{(k)}(u)\geq a_t+\rho(u)-g(u)\right)_{g(\cdot)=Y_k(\cdot)+1} \right)dx 
\\
&=:&\underset{k=1}{\overset{t_b+1}{\sum}} \circledast_{L,b}(k).
\end{eqnarray*}

We remark that for any $ {\bf c}>0 $ 
\begin{equation*}
 \int_{[0,R]^\d} \E\Big(\frac{\1_{\{ Y_\cdot(x)\in \DDi^{\rho(x),L}_t\}}}{{\bf r}_t(x)^\d}\Big)dx \leq  \int_{[0,R]^\d} c\ee^{\d (t+{\bf c})} \E\left(  \1_{\{Y_\cdot(x)\in \DDi^{\rho(x),L}_t\}}\right) + \sum_{p\geq t+{\bf c} }  \ee^{\d p} \E\left(\1_{\{{\bf r}_t(x)\leq \ee^{-p}/4 \}} \1_{\{ Y_\cdot(x)\in \DDi^{\rho(x),L}_t\}}\right)dx.
\end{equation*}
By the arguments C), D), E) and F) pp 28 used to prove (\ref{valseme1}), we can affirm that
\begin{equation}
\label{tensloc}
\int_{[0,R]^\d}\E\left(\frac{\1_{\{ Y_\cdot(x)\in \DDi^{\rho(x),L}_t\}}}{{\bf r}_t(x)}\right) dx \leq c(1+L)^2\int_{[0,R]^\d}\rho(x)\ee^{-\sqrt{2\d}\rho(x)} dx,
\end{equation}
Using (\ref{tensloc}) we will bound $\circledast_{L,b}(k)$ by distinguishing three cases:

\paragraph{      A) If $k\leq 5\log l$.} As $x$ is good$_k$, $ \rho(\cdot)\in \mathcal{C}_R(l,\kappa_\d \log l,\log t) $ and $u\in A_k(x)$,  
\begin{eqnarray*}
\rho(u)-\underset{v\in A_k(x)}{\sup} Y_k(v) -1 &\geq& -[e_k+\frac{D}{2}]-Y_k(x)+\rho(u)-1
\\
&\geq& \rho(u) -(\log l)^\frac{2}{3}-(5 \log l)^\frac{1}{12}-\frac{D}{2}-1
\\
&\geq & \rho(u) -(\log l)^\frac{1}{12}-\frac{D}{2}-1\geq \frac{\kappa_\d}{2}\log l,
\end{eqnarray*}
once $l\geq \ee^{D^2}$. By using in turn the scaling property (\ref{scalingl}) and then the invariance by translation we get that
\begin{eqnarray*}
\circledast_{L,b}(k)&\leq&  \int_{[0,R]^\d}\E\left(\frac{\1_{\{ Y_\cdot(x)\in \DDi_t^{\rho(x),L} \}}}{{\bf r_t}(x)^\d} \P\left( \exists u\in  A_k(x) ,\, Y_{t-k}^{(k)}(u)\geq  a_t+ \frac{\kappa_\d}{2}\log l\right) \right)dx
\\
&\leq & \int_{[0,R]^\d}\E\left(\frac{\1_{\{ Y_\cdot(x) \in \DDi_t^{\rho(x),L} \}}}{{\bf r_t}(x)^\d} \P\left(  \exists v \in A^*,\,  Y_{t-k}(v)  \geq a_t+  \frac{\kappa_\d}{2}\log l \right)  \right)dx .
\end{eqnarray*}
By applying Proposition \ref{tension*} (with the constant function $x\mapsto \frac{\kappa_\d}{2} \log l\in \mathcal{C}_{A^*}(1,10, +\infty) $) we have for some $\beta>0$,
\begin{eqnarray*}
 \P\left(  \exists v \in A^*,\,  Y_{t-k}(v)  \geq a_t+  \frac{\kappa_\d}{2}\log l \right) &\leq& c_6 \int_{A^*} \frac{\kappa_\d}{2}\log l\ee^{-\sqrt{2\d} \frac{\kappa_\d}{2}\log l}du
\\
&\leq &c l^{-\beta},\qquad \text{(notice that $\lambda(A^*)\leq 1$)}.
\end{eqnarray*}
Finally  $\circledast_{L,b}(k)\leq c  l^{-\beta} \int_{[0,R]^\d}\E\big(\frac{\1_{\{ Y_\cdot(x)\in \DDi_t^{\rho(x),L} \}}}{{\bf r_t}(x)^\d} \big)dx  $ then, by applying (\ref{tensloc}) we get that
\begin{eqnarray}
\nonumber \sum_{k=0}^{5\log l}  l^{-\beta}\circledast_{L,b}(k) &\leq& c  l^{-\beta}\sum_{k=0}^{5\log l}  \int_{[0,R]^\d}\E\left(\frac{\1_{\{ Y_\cdot(x) \in\DDi_t^{\rho(x),L} \}}}{{\bf r_t}(x)^\d} \right)dx.
\\
\label{eqq0} &\leq&  c' (1+L)^2 l^{-\beta}\mathtt{I}_\d(\rho),
\end{eqnarray}
which is smaller than $\epsilon \mathtt{I}_\d(\rho)$ for $l$ large enough.

\paragraph{     B) If $5\log l \leq k\leq \frac{t}{2}$.} As $x$ is good$_k$, and $u\in A_k(x)$,  
\begin{eqnarray*}
\rho(u)-\underset{v\in A_k(x)}{\sup} Y_k(v)-1&\geq& -[e_k+\frac{D}{2}]-Y_k(x)+\rho(u)-1
\\
&\geq & -[e_k+\frac{D}{2}]+ (-\rho(x)+4e_k-D) +\rho(u)-1
\\
&\geq & 3e_k -\frac{3}{2}D -2 .
\end{eqnarray*}
For the last inequality, recall that $\rho(\cdot)\in \mathcal{C}_R(l,\kappa_\d \log l,\log t)$,  $|u-x|\leq c\ee^{-5 \log l}$ if $k\geq 5\log l$  and thus $|\rho(u)-\rho (x)|\leq 1$. In addition with the scaling property (\ref{scalingl}) then the invariance by translation we get that
\begin{eqnarray*}
\circledast_{L,b}(k)&\leq&  \int_{[0,R]^\d}\E\left(\frac{\1_{\{ Y_\cdot(x)\in \DDi_t^{\rho(x),L} \}}}{{\bf r_t}(x)^\d} \P\left( \exists u\in  A_k(x) ,\, Y_{t-k}^{(k)}(u)\geq  a_t+  3e_k -\frac{3}{2}D -2 \right)  \right)dx
\\
&\leq & \int_{[0,R]^\d}\E\left(\frac{\1_{\{ Y_\cdot(x)\in \DDi_t^{\rho(x),L} \}}}{{\bf r_t}(x)^\d} \P\left(  \exists v \in A^*,\,  Y_{t-k}(v)  \geq a_t+  3e_k -\frac{3}{2}D -2  \right)  \right)dx .
\end{eqnarray*}
 By applying Proposition \ref{tension*} then (\ref{tensloc}) we get that
\begin{eqnarray*}
 \circledast_{L,b}(k)&\leq& c\int_{[0,R]^\d}\E\left(\frac{\1_{\{ Y_\cdot(x) \in \DDi_t^{\rho(x),L}\}}}{{\bf r_t}(x)^\d}\right)\ee^{-3e_k +\frac{3}{2}D }  dx    \leq c'(1+L)^2  \mathtt{I}_\d(\rho)  \ee^{-3e_k +\frac{3}{2}D }.
\end{eqnarray*}
Finally, there exists $c>0$ such that  
\begin{eqnarray}
\nonumber\underset{k=5\log l}{\overset{\lfloor \frac{t}{2}\rfloor}{\sum}} \circledast_{L,b}(k) \leq  c \ee^{\frac{3}{2}D}\underset{k=5\log l}{\overset{\lfloor \frac{t}{2}\rfloor}{\sum}} (1+L)^2 \mathtt{I}_\d(\rho)  \ee^{-3e_k }&\leq & c'(1+L)^2\ee^{\frac{3}{2}D}\ee^{-3e_{\log l}}\mathtt{I}_\d(\rho)
\\
\label{eqq1}&\leq &\epsilon \mathtt{I}_\d(\rho).
\end{eqnarray}
once $l$ is large enough.

\paragraph{      C) If  $t-b+1\geq k\geq \frac{t}{2}$.} As $x$ is good$_k$ and $u\in A_k(x)$, 
\begin{eqnarray*}
a_t+\rho(u)-\underset{v\in A_k(x)}{\sup} Y_k(v)-1&\geq& -[e_k+\frac{D}{2}]-Y_k(x)+a_t+\rho(u)-1
\\
&\geq & 3e_k -\frac{3D}{2}-L + \rho(u)-\rho(x)-1  \geq 3e_k -3D/2 -L-2.
\end{eqnarray*}
For the last inequality, recall that $\rho(\cdot )\in \mathcal{C}_R(l, \kappa_\d \log l,\log t)$,  $|u-x|\leq c\ee^{-\frac{t}{2}}$ if $k\geq t/2$  and thus $|\rho(u)-\rho (x)|\leq 1$. According to Lemma \ref{maxfx} (with the constant function $x\mapsto 3e_k -3D/2 -L-2\in \mathcal{C}_R(1,10, +\infty) $) one has
\begin{eqnarray*}
 \P\left( \exists u\in  A_k(x) ,\, Y_t^{(k)}(u)\geq a_t+\rho(u)-g(u)\right)_{\Big| g(u)=Y_k(u)}&\leq & \P\left( \exists u\in  A^* ,\, Y_t(u)\geq   3e_k -3D/2-L -2\right)
\\
&\leq& c_3 \ee^{- [3e_k -3D/2-L -2]}.
 \end{eqnarray*}
Finally with in addition (\ref{tensloc}) we get that
\begin{eqnarray*}
&&\circledast_{L,b}(k)=
\\
 &&  \int_{[0,R]^\d}\E\left(\frac{\1_{\{ Y_\cdot(x)\in \DDi_t^{\rho(x),L},\, x \, \text{is } L-\text{good}_k\}}}{{\bf r_t}(x)^\d} \P\left( \exists u\in  A_k(x) ,\, Y_{t-k}^{(k)}(u)\geq a_t+\rho(u)-g(u)\right)_{\Big| g(u)=Y_k(u)} \right)dx 
\\
&\leq & c\ee^{-[3e_k -3D/2-L]} \int_{[0,R]^\d}\E\left(\frac{\1_{\{ Y_\cdot(x) \in \DDi_t^{\rho(x),L}\}}}{{\bf r_t}(x)^\d}\right)  dx  
\\
&\leq&  c'\ee^{-[3e_k -3D/2-L]} (1+L)^2 \int_{[0,R]^d}\rho(x)\ee^{-\sqrt{2\d}\rho(x)}dx,
\end{eqnarray*}
and thus 
\begin{eqnarray*}
\underset{k=\frac{t}{2}+1}{\overset{t_b+1}{\sum}} \circledast_{L,b}(k)&\leq  &\underset{k=\frac{t}{2}+1}{\overset{t_b+1}{\sum}} c  (1+L)^2 \ee^{-[3e_k -3D/2 -1]} \mathtt{I}_\d(\rho) 
\\
&\leq &  c'(1+L)^2 \ee^{ (\frac{3}{2}D+L)}\ee^{- b^\frac{1}{12}} \mathtt{I}_\d(\rho).
\end{eqnarray*}
This yields that there exists $b_0(D,L)\geq 1$, such that for any $b\geq b_0$, and any $t\geq 1$, $\rho(\cdot) \in \mathcal{C}_R(l,\kappa_\d \log l,\log t)$ we have
\begin{equation}
\label{eqq2}
\underset{k=\frac{t}{2}+1}{\overset{t_b+1}{\sum}} \circledast_{L,b}(k)\leq \epsilon \mathtt{I}_\d(\rho).
\end{equation}
From (\ref{eqq0}), (\ref{eqq1}) and (\ref{eqq2}) we get Lemma \ref{htoui}.\hfill$\Box$

\subsection{Proof of Lemma  \ref{controlF}}
 
Fix $L,b>1$. We shall prove that $F_{L,b}$ is $b$ {\bf regular} with 
\begin{eqnarray}
\label{defhLb}
h(z)=h_{L,b}(z)&:=&  \ee^{-\sqrt{2\d}(z+L)}\P_{z+L+1}\left( Y_b(0)\geq 0\right)^{\frac{1}{2}} ,
\\
\label{defF}F^*(z)=F^*_{b}(g)&:=& \underset{z\in \r}{\sup}\, \E_z\left(\frac{\1_{\{\exists y\in B(0,1),\, Y_b(y)\geq  -g(y\ee^b)  \}} }{[\lambda_{B(0,1)}\left(\{ y,\, Y_b(y)\geq -g(y\ee^b)  -\frac{1}{2}\}\right)]^8}\right)^{\frac{1}{4}} .
\end{eqnarray}
 \nomenclature[i93]{$h_{L,b}(z) $}{$ :=\ee^{-\sqrt{2\d}(z+L)}\P_{z+L+1}\left( Y_b(0)\geq 0\right)^{\frac{1}{2}} $}\nomenclature[i94]{$F^*_{b}(g) $}{$ :=  \underset{z\in \r}{\sup}\, \E_z\left(\frac{\1_{\{\exists y\in B(0,1),\, Y_b(y)\geq  g(y\ee^b)  \}} }{[\lambda_{B(0,1)}\left(\{ y,\, Y_b(y)\geq g(y\ee^b)  -\frac{1}{2}\}\right)]^8}\right)^{\frac{1}{4}} $}
\noindent{\it Proof of Lemma \ref{controlF}.} We will show that $h_{L,b}$ and $F^*_b$ satisfy (i), (ii), (iii) and (iv) of Definition \ref{defbregular}.
\\

-To check (i), observe that there exists $c>0 $ such that 
\begin{equation}
\label{et2}
\underset{z\in \r}{\sup}\,h_{L,b}(z)\leq c  ,\quad  \text{and }\quad h_{L,b}(z)\leq  \ee^{z},\, \text{if }   z< -(2b+L+1).
\end{equation}

 -Now we shall prove (ii): Let $g\in \mathcal{C}(B(0,\ee^b),\r)$ such that $w^{(0,1)}_{g(\cdot \ee^b)}(\delta)\leq \frac{1}{4}$. We define $\Lambda=\lambda_{B{(0,1)}}(\{y, Y_b(y)\geq -g(y\ee^b)-\frac{1}{2}\}  ) $.
 On the set $\{ \exists y\in B(0,1),\, Y_b(y)\geq -g(y\ee^{b}) \}$, we introduce ${\bf r} $ the biggest radius such that {\it $\exists x_{\bf r}$ with $B(x_{\bf r}, {\bf r})\subset B(0,1)$; $\exists z_{\bf r}\in B(x_{\bf r}, {\bf r})$ with $Y_b(z_{\bf r})\geq -g(z_{\bf r}\ee^b)$;  $\forall y\in B(x_{\bf r}, {\bf r})$, $Y_b(y)\geq -g(y\ee^b)-\frac{1}{2}$.} By (\ref{defF}),
 \begin{eqnarray*}
 F^*_{b}(g)^4 &=&  \underset{z\in \r}{\sup}\,\E_z\left(\frac{\1_{\{\exists y\in B(0,1),\, Y_b(y)\geq -g(y\ee^{b}) \}}}{\Lambda^8}\right) 
    \\
  &&  \leq  S^{-8}(\ee^{b\d}/\delta)^8+ \underset{k=\ee^{b}/\delta}{\overset{\infty}{\sum}}  S^{-8}(k+1)^{8}  \underset{z\in \r}{\sup}\,\E_{z}\left(   \1_{\{\exists y\in B(0,1),\, Y_{b}(y)\geq -g(y\ee^b) \}} \1_{\{\frac{S}{(k+1)^\d}\leq \Lambda \leq \frac{S}{k^\d}\}}\right),
 \end{eqnarray*}
 with $S$ the volume of the unit ball. Clearly, $\Lambda\leq S (\frac{1}{k})^\d$ implies $  {\bf r}\leq \frac{1}{k}$, moreover on  $\{ {\bf r} \leq \frac{1}{k}<\delta \}$, 
 \begin{eqnarray*}
 1=\1_{\{\underset{ \underset{|x-y|\leq \frac{1}{k}}{x,y\in B(0,1)}}{\sup}|Y_b(x)-Y_b(y)|\geq \frac{1}{2}-w^{(0,1)}_{g(.\ee^b)}(\delta) \}}\leq \1_{\{\underset{ \underset{|x-y|\leq \frac{1}{k}}{x,y\in B(0,1)}}{\sup}|Y_b(x)-Y_b(y)|\geq \frac{1}{4}\}}.
 \end{eqnarray*}
 Furthermore by (\ref{boubou2}) (with $h=\frac{1}{k},\, m=2k,\, p=2,\, l=b$ and $x=c\ee^{-b}k $) we have 
 $$\underset{z\in \r}{\sup}\,\P_z\Big( \underset{\underset{|x-y|\leq \frac{1}{k}}{x,y\in B(0,1)}}{\sup}|Y_b(x)-Y_b(y)|\geq \frac{1}{4}  \Big)=\P_0\Big( \underset{\underset{|x-y|\leq \frac{1}{k}}{x,y\in B(0,1)}}{\sup}|Y_b(x)-Y_b(y)|\geq \frac{1}{4}  \Big)\leq c'\ee^{-\frac{1}{c''}\ee^{-b} k}.$$ Finally
 $
 F^*_{b}(g)^4\leq S^{-8}\ee^{8b\d}/\delta^8 + \underset{k=1+\ee^b/\delta}{\overset{\infty }{\sum}}  S^{-8}(k+1)^{8}c\ee^{-\frac{1}{c''}\ee^{-b} k} \leq  \ee^{4c_I b}\delta^{- 8}$, which suffices to prove of (ii).
 \\

-Check (iii) stems easily from the definition of $F_{L,b}$ in (\ref{defFLb}) and the Cauchy-Schwarz inequality.
 \\ 
 
-It remains to prove (iv). Let $g_1,\, g_2$ two continuous functions from $B(0,\ee^b)\to \r$ such that $||g_1-g_2||_\infty= \delta<\frac{1}{8}$. Let us define (only for this proof) $\forall g\in \mathcal{C}(B(0,\ee^b),\r)$ and $\gamma\in \r$:
\begin{eqnarray*}
&&M(g):=\underset{y\in B(0,1)}{\sup}(Y_b(y)+g(y\ee^b)), \qquad \Lambda_g(\gamma) :=\lambda_{B(0,1)}(\{y,\, Y_b(y)\geq -g(y\ee^b)+\gamma  \}).
\end{eqnarray*}
With these two notations we have:
 \begin{equation}
 \label{bb2}|F_{L,b}(z,g_1)-F_{L,b}(z,g_2)|\leq \ee^{-\sqrt{2\d}(z+L)} \E_{z+L+1}\left(\1_{\{  Y_b(0)\geq 0\}} \left|  \frac{\1_{\{ M({g_1})\geq 1 \}}}{\Lambda_{g_1}(0)} - \frac{\1_{\{  M({g_2})\geq 1 \}}}{\Lambda_{g_2}(0)}\right|\right).
\end{equation}
By the triangular inequality observe that
\begin{eqnarray*}
 \left|  \frac{\1_{\{ M({g_1})\geq 1 \}}}{\Lambda_{g_1}(0)} - \frac{\1_{\{  M({g_2})\geq 1 \}}}{\Lambda_{g_2}(0)}\right| &\leq& \left|\frac{\1_{\{M(g_1)\geq 1\}}-\1_{\{M(g_2)\geq 1\}}}{\Lambda_{g_1}(0)} \right| +\left|\1_{\{M(g_2)\geq 1\}}\Big(\frac{1}{\Lambda_{g_1}(0)}-\frac{1}{\Lambda_{g_2}(0)}\Big) \right|
 \\
&\leq & \frac{\1_{\{M(g_1)\in [1-\delta,1+\delta]\}}}{\Lambda_{g_1}(0)}+ \frac{\Lambda_{g_1}(-\delta)-\Lambda_{g_1}(\delta)}{\Lambda_{g_1}(0) \Lambda_{g_2}(0) }
\end{eqnarray*}
where we have used $||g_1-g_2||_\infty=\delta$. Furthermore from Theorem 3.1 in \cite{PTr79}, as $Var(Y_b(y))=b\geq 1,\,\forall y\in B(0,1)  $,  we can affirm that {\it there exists $c>0 $ such that for any $\delta\in(0,1) $, $g\in  \mathcal{C}(B(0,\ee^b),\r)$,
\begin{equation}
\label{PTa}
\underset{z\in \r}{\sup}\, \P\left(M(g)\in [z-\delta,z+\delta]\right)\leq c \delta.
\end{equation}}
Going back to (\ref{bb2}), we have
\begin{eqnarray*}
&& \ee^{\sqrt{2\d}(z+L)}|F_{L,b}(z,g_1)-F_{L,b}(z,g_2)|
 \\
&&\leq  \E_{z+L+1}\left(\frac{\1_{\{ Y_b(0)\geq 0,\, M(g_1)\in [1-\delta,1+\delta]\}}}{\Lambda_{g_1}(0)}  \right) +  \E_{z+L+1}\left(\1_{\{ Y_b(0)\geq 0,\, M(g_2)\geq 1 \}}  \frac{\Lambda_{g_1}(-\delta)-\Lambda_{g_1}(\delta)}{\Lambda_{g_1}(0) \Lambda_{g_2}(0)} \right)
\\
&& := (A)+(B).
\end{eqnarray*}
By applying twice the Cauchy-Schwarz inequality then (\ref{PTa}) to $(A)$ we get that
\begin{eqnarray*}
(A)&\leq & \P_{z+L+1}\left( Y_b(0)\geq 0\right)^{\frac{1}{2}}\times\E_{z+L+1}\left(\frac{\1_{\{M(g_1)\geq 1  -\delta \}} }{\Lambda_{g_1}^4(\delta)}\right)^{\frac{1}{4}} \times \P_{z+L+1}\left( M(g_1)\in [1-\delta ,1+\delta]\right)^{\frac{1}{4}} 
\\
&\leq& c \P_{z+L+1}\left( Y_b(0)\geq 0\right)^{\frac{1}{2}}\times\E_{z+L+1+\delta}\left(\frac{\1_{\{ M(g_1)\geq 1  \}} }{\Lambda_{g_1}^4(2\delta)}\right)^{\frac{1}{4}}\delta^\frac{1}{4}.
\end{eqnarray*}
Similarly, observing that  $\min(\Lambda_{g_1}(0),\Lambda_{g_2}(0))\geq \Lambda_{g_1}(\frac{1}{4})$, we deduce that
\begin{eqnarray*}
(B)&=&\int_{B(0,1)}\E_{z+L+1}\left(\frac{\1_{\{ Y_b(0)\geq 0,\,  M(g_2)\geq 1 \}}}{\Lambda_{g_1}(0)\Lambda_{g_2}(0)}  \1_{\{Y_b(x)+ g_{1}(x\ee^b) \in [-\delta ,\delta]\}} \right)dx
\\
&&\leq \P_{z+L+1}\left( Y_b(0)\geq 0\right)^{\frac{1}{2}}                     \E_{z+L+1}\left(\frac{\1_{\{ M(g_1)\geq 1-\delta\}} }{[\Lambda_{g_1}(\frac{1}{4})]^8}\right)^{\frac{1}{4}} \int_{B(0,1)}\P_{z+L+1+ g_{1}(x\ee^b) }\left(  Y_b(x)\in [-\delta,\delta]  \right)^\frac{1}{4}dx
\\
&&\leq c\P_{z+L+1}\left( Y_b(0)\geq 0\right)^{\frac{1}{2}}                     \E_{z+L+1+\delta}\left(\frac{\1_{\{\exists y\in B(0,1),\, Y_b(y)\geq  -g_{1}(y\ee^b)+1  \}} }{[\Lambda_{g_1}(\frac{1}{4}+\delta)]^8}\right)^{\frac{1}{4}} \delta^{\frac{1}{4}}.
\end{eqnarray*}
From the bound on (A) and (B) we deduce that 
\begin{equation*}
|F_{L,b}(z,g_1)-F_{L,b}(z,g_2)|\leq \ee^{-\sqrt{2\d}(z+L)} 2c\P_{z+L+1}(Y_b(0)\geq 0)^\frac{1}{2} ||g_1-g_2 ||_\infty^{\frac{1}{4}} F^*_b(g_1),
\end{equation*} 
which proves (iv). 
\hfill$\Box$
\\

\section{Proof of Theorem  \ref{RESUME}}

For any $\sigma \in [0,t_b]$, $
\mathfrak{G}_{t,b,\sigma} :B{(0,\ee^{b})}\ni y\mapsto   -\int_{t_b-\sigma}^{t_b}\mathtt{g}(\ee^{s-t} y)dB_{s} -  \zeta_{t_b}(y\ee^{-t})+ Z_{t_b}^{0}(y\ee^{-t})$. The Theorem \ref{RESUME} is a combination of the two following lemmas.

  \begin{Lemma}
 \label{sansrho}
  Let $b>0$ and $F:\r\times \mathcal{C}(B(0,e^b),\r)\to \r^+$ be a function {\bf $b$ regular}. For any $\epsilon>0$,  there exist $ l,\, T >0$ such that for any $t\geq T$, $\rho(\cdot) \in \mathcal{C}_R(l,\kappa_\d \log l,\log t) $, $\mathtt{z}\leq (\log t)^{30}$,
 \begin{eqnarray}
 \nonumber
 \Big| \int_{\R_t} e^{-\sqrt{2\d}\rho(x)} \E_{-\rho(x)} \Big[\1_{\{ \overline{B}_{t_b }\leq 0,\, \overline{B}_{[\frac{t}{2},t_b] } \leq -\mathtt{z} \}}\Big(F(B_{t_b}+\mathtt{z},\mathfrak{G}^{\rho_x}_{t,b} )-F(B_{t_b}+\mathtt{z},\mathfrak{G}_{t,b})\Big)\Big] dx\Big|
 \\
 \label{eqsansrho} \leq \epsilon \int_{[0,R]^\d}   \rho(x)  e^{-\sqrt{2\d}\rho(x)} dx.
 \end{eqnarray}
  \end{Lemma}
  Recall the definition (\ref{defMreload}).

\begin{Lemma}
\label{passageaueps}
(i) Let $b>0$ and $F:\r\times \mathcal{C}(B(0,e^b),\r)\to \r^+$ be a function {\bf $b$ regular}. For any $\epsilon>0$,  there exist $M,\,\sigma,\, T >0$ such that for any $t\geq T$, $\alpha\in[1,\log t]$, $\mathtt{z}\leq (\log t)^{30}$,
\begin{equation}
\label{eqpassageaueps}
\frac{t^{\frac{3}{2}}}{\alpha }\Big|\E_{-\alpha} \Big[\1_{\{   \overline{B}_{t_b }\leq 0,\, \overline{B}_{[\frac{t}{2},t_b] } \leq -\mathtt{z}   \}}\Big(F^{(M)}(B_{t_b}+\mathtt{z},\mathfrak{G}_{t,b,\sigma} )-F(B_{t_b}+\mathtt{z},\mathfrak{G}_{t,b})\Big) \Big] \Big|\leq \epsilon.
\end{equation}

(ii) Let $b>0$ and $F:\r\times \mathcal{C}(B(0,e^b),\r)\to \r^+$ be a function {\bf $b$ regular}. Fix $M,\sigma > 0$. There exists $C_{M,\sigma}(F)>0$ such that for any $\epsilon > 0$, there exists $T\geq 0$ such that for any  $t\geq T$, $\alpha \in [1,\log t]$, $\mathtt{z}\leq (\log t)^{30}$,
\begin{equation}
\label{eqpassageaueps2}
\Big|\frac{t^{\frac{3}{2}}}{\alpha } \E_{-\alpha}\Big[\1_{\{   \overline{B}_{t_b }\leq 0,\, \overline{B}_{[\frac{t}{2},t_b] } \leq -\mathtt{z}\}}F^{(M)}\big(B_{t_b}+\mathtt{z},\mathfrak{G}_{t,b,\sigma}\big)\Big]-C_{M,\sigma}(F)\Big|\leq \epsilon.
\end{equation}
\end{Lemma}
Displays (\ref{eqsansrho}) and (\ref{eqpassageaueps}) may to replace $\mathfrak{G}_{t,b}^{\rho_x}$ by ${\mathfrak{G}_{t,b,\sigma}}$ in the argument of $F$. Then thanks to the properties of ${\mathfrak{G}_{t,b,\sigma}} $ we can prove the renewal result (\ref{eqpassageaueps2}). Theorem \ref{RESUME} is obtained by replacing $\alpha$ by $\rho(x)$, then integrating on $[0,R]^\d$, the displays (\ref{eqpassageaueps}) and (\ref{eqpassageaueps2}).

Before to tackle the proof of Lemma \ref{sansrho} we need a control on the function $F^*$ and $h$ associated to the {\bf $b$ regularity} of $F$:


\begin{Lemma}
\label{rajoutfin}
Let $h$ and $F^*$ the two functions associated to $F$ a function $b$ {\bf regular}.  There exists constants $c>0 $ (depending  on $F$, $h$ or $F^*$) and $T>0$  such that for any $t\geq T$, $\alpha\in [1,\log t],\, \sigma \in [0,t_b] $ and $\mathtt{z}\leq (\log t)^{30}$
\begin{eqnarray}
\label{sansrhofirst} \E_{-\alpha}\Big[  \1_{\{ \overline{B}_{t_b }\leq 0,\, \overline{B}_{[\frac{t}{2},t_b] } \leq -\mathtt{z} \}}  h(B_{t_b}+\mathtt{z})F^*(\mathfrak{G}_{t,b,\sigma})  \Big]\leq  {c}\alpha t^{-\frac{3}{2}}.
\end{eqnarray} 

\end{Lemma}

%
%
\noindent{\it Proof of Lemma \ref{rajoutfin}.} By (\ref{trajB.21}), we can affirm that   for any $t\geq b>0$ large enough, $\alpha \in [1,\log t]$, $k,\, j,\, \geq 1$, $\mathtt{z}\in [1,(\log t)^{30}]$ and $\sigma\in[0,t_b]$, 
 \begin{equation} 
\label{Gcond1}  \E_{-\alpha}\Big[  \1_{\{ \overline{B}_{t_b}\leq 0,\, \overline{B}_{[\frac{t}{2},t_b]} \leq -\mathtt{z},\, B_{t_b}+\mathtt{z}\in [-(k+1),k] \}}\1_{\{   w^{(0,1)}_{\mathfrak{G}_{t,b,\sigma}(\cdot \ee^b)}(\frac{1}{j})\geq \frac{1}{4}\}}   \Big]    \leq  c_{23}(1+k)\frac{ \alpha}{t^{\frac{3}{2}}} e^{-c_{24}(b)j}.
\end{equation} 
According to ({\ref{et2rajout}), there exists $c_1(h)>0$ such that 
\begin{eqnarray}
\label{appli5.27} h(B_{t_b}+\mathtt{z})\leq \Big\{ \begin{array}{ll} e^{ B_{t_b}+\mathtt{z}} & \text{if  } B_{t_b}+\mathtt{z}\leq -c_1(h), \\
c_1(h) &\text{if   }  B_{t_b}+\mathtt{z} \geq -c_1(h).
\end{array}\Big.
\end{eqnarray}
By continuity of $y\mapsto \mathfrak{G}_{t,b,\sigma}(y)$,  
$$1\leq \1_{\{ w^{(0,1)}_{\mathfrak{G}_{t,b,\sigma}(\cdot \ee^b)}(1)\leq \frac{1}{4}   \}} + \sum_{j\geq 1} \1_{\{ w^{(0,1)}_{\mathfrak{G}_{t,b,\sigma}(\cdot \ee^b)}(j^{-1})\geq \frac{1}{4} \geq  w^{(0,1)}_{\mathfrak{G}_{t,b,\sigma}(\cdot \ee^b)}((j+1)^{-1})  \}}. $$
Thanks to (\ref{et1}), for any $j\geq 1$, on $\{\frac{1}{4} \geq  w^{(0,1)}_{\mathfrak{G}_{t,b,\sigma}(\cdot \ee^b)}(j^{-1})  \}$, $$F^*(\mathfrak{G}_{t,b,\sigma})\leq c j^{10}.$$
Combining these two inequalities with (\ref{appli5.27}), we get that
\begin{align*}
\E_{-\alpha}&\Big[  \1_{\{ \overline{B}_{t_b }\leq 0,\, \overline{B}_{[\frac{t}{2},t_b] } \leq -\mathtt{z} \}}   h(B_{t_b}+\mathtt{z})F^*(\mathfrak{G}_{t,b,\sigma})  \Big]
\\
\leq    &  c_1(h)\sum_{j=1}^\infty  (j+1)^{10}  \Big(\sum_{k=0}^{  c_1(h)} \E_{-\alpha} \Big[ \1_{\{  \overline{B}_{t_b }\leq 0,\, \overline{B}_{[\frac{t}{2},t_b] } \leq -\mathtt{z},\, B_{t_b}+\mathtt{z}\in [-(k+1),-k] \}} \1_{\{   w^{(0,1)}_{\mathfrak{G}_{t,b,\sigma}(\cdot \ee^b)}(j^{-1})\geq \frac{1}{4}  \}}   \Big]  
\\
& +\sum_{k= c_1(h)}^{\infty} e^{- k} \E_{-\alpha} \Big[ \1_{\{ \overline{B}_{t_b }\leq 0,\, \overline{B}_{[\frac{t}{2},t_b] } \leq -\mathtt{z}   ,\, B_{t_b}+\mathtt{z}\in [-(k+1),-k] \}} \1_{\{  w^{(0,1)}_{\mathfrak{G}_{t,b,\sigma}(\cdot \ee^b)}(j^{-1})\geq \frac{1}{4}  \}}   \Big]\Big).
\end{align*}
Finally according to (\ref{Gcond1}) we have for any $t>0$ large enough, $\alpha \in [1,\log t]$, $k,\, j,\, \geq 1$, $\mathtt{z}\in [1,(\log t)^{30}]$ and $\sigma\in[0,t_b]$, 
\begin{eqnarray}
\nonumber \E_{-\alpha}\Big[   \1_{\{ \overline{B}_{t_b }\leq 0,\, \overline{B}_{[\frac{t}{2},t_b] } \leq -\mathtt{z} \}} h(B_{t_b}+\mathtt{z})F^*(\mathfrak{G}_{t,b,\sigma})  \Big]&\leq &  \frac{c \alpha}{t^{\frac{3}{2}}}\sum_{j=1}^\infty \frac{(j+1)^{10} }{e^{c_{24}(b)j}}\Big[  c_1(h)^2  + \sum_{k= c_1(h)}^{\infty} \frac{(1+k)}{e^{ k}} \Big]
\\
\label{ineqPrat}&\leq& \frac{ c'\alpha}{t^\frac{3}{2}},
\end{eqnarray}
which ends the proof of Lemma \ref{rajoutfin}.

\begin{Remark}As a by product we have also shown the following affirmation.
 Fix $F$ a function {\bf $b$ regular}. For any $\epsilon>0$ there exists $M,\, T>0$  such that for any $t\geq T$, $\alpha\in [1, \log t]$, $\mathtt{z}\leq (\log t)^{30}$ and $\sigma\in [0,t_b]$  we have
\begin{equation}
\label{0.3}
\E_{-\alpha}\Big(  \1_{\{   \overline{B}_{t_b }\leq 0,\, \overline{B}_{[\frac{t}{2},t_b] } \leq -\mathtt{z} \}} h(B_{t_b}+\mathtt{z})  F^*\Big(\mathfrak{G}_{t,b,\sigma} \Big)  \1_{\{ w^{(0,1)}_{\mathfrak{G}_{t,b,\sigma}(\cdot \ee^b)}(\frac{1}{M})\geq \frac{1}{4}\}} \Big)\leq \frac{\epsilon \alpha}{t^{\frac{3}{2}}},
\end{equation}
and 
\begin{equation}
\label{0.4}
\E_{-\alpha}\Big(  \1_{\{  \overline{B}_{t_b }\leq 0,\, \overline{B}_{[\frac{t}{2},t_b] } \leq -\mathtt{z} \}} h(B_{t_b}+\mathtt{z})  F^*\Big(\mathfrak{G}_{t,b,\sigma} \Big)  \1_{\{  B_{t_b} +\mathtt{z}\leq -M\}} \Big)\leq \frac{\epsilon \alpha}{t^{\frac{3}{2}}}.
\end{equation} 
Indeed for (\ref{0.3}) as well as for (\ref{0.4}), it suffices to choose $M\geq c_1(h)$ large enough such that (see \ref{ineqPrat}))
$$ 
  c\sum_{j=M}^\infty \frac{(j+1)^{10} }{e^{c_{24}(b)j}}\Big[  c_1(h)^2  + \sum_{k= c_1(h)}^{\infty} \frac{(1+k)}{e^{ k}} \Big]+   c\sum_{j=1}^\infty \frac{(j+1)^{10} }{e^{c_{24}(b)j}}  \sum_{k= M}^{\infty} \frac{(1+k)}{e^{ k}}\leq  \epsilon.
$$ 
\end{Remark}

\noindent{\it  Proof of Lemma \ref{sansrho}.} For $t\geq \log l+b$,  as $\rho(\cdot) \in \mathcal{C}_R(l,\kappa_\d \log l,\log t)$,
\begin{eqnarray*}
||\mathfrak{G}_{t,b}^{\rho_x}-\mathfrak{G}_{t,b}||_\infty\leq  \underset{x\in \R_t,\,y \in B{(0,e^b)}}{\sup}|\rho(x+y\ee^{-t})-\rho(x)|   \leq \ee^{-\frac{t-b}{3}}.
\end{eqnarray*}
Recalling (\ref{et3bibi}), the quantity in \eqref{eqsansrho} is smaller than:
\begin{eqnarray*}
 \int_{\R_t} {e^{-\sqrt{2d}\rho(x)}}{t^\frac{3}{2}}    \E_{-\rho(x)}   \Big[ \1_{\{ \overline{B}_{t_b }\leq 0,\, \overline{B}_{[\frac{t}{2},t_b] } \leq -\mathtt{z} \}} c_9 e^{-\frac{t-b}{12}} h(B_{t_b}+\mathtt{z})F^*(\mathfrak{G}_{t,b})  \Big]dx.
\end{eqnarray*}
with $h$ and $F^*$ the two functions associated to the {\bf $b$ regular} function $F$. Now we conclude with Lemma \ref{rajoutfin} applied with $\sigma =t_b$.
\hfill$\Box$

\medskip

\medskip

 \noindent{\it  Proof of  (\ref{eqpassageaueps}).}  Let    $ b,\, \epsilon>0  $ and $F$ {\bf $b$ regular}. We have to study the expectation under  $\E_{-\alpha} $ of 
 \begin{eqnarray*}
  \1_{\{  \overline{B}_{t_b }\leq 0,\, \overline{B}_{[\frac{t}{2},t_b] } \leq -\mathtt{z} \}}\Big|F^{ (M)}(B_{t_b}+\mathtt{z},\mathfrak{G}_{t,b,\sigma})-F (B_{t_b} +\mathtt{z},\mathfrak{G}_{t,b})  \Big|.
 \end{eqnarray*}
Thanks to (\ref{0.3}) and (\ref{0.4}) we can choose  $M_1  $ large enough to restrain our study to the expectation of 
\begin{equation}
\label{restrain}
\1_{\{ w^{(0,1)}_{\mathfrak{G}_{t,b,\sigma}(\cdot \ee^b)}(\frac{1}{M_1})\leq \frac{1}{4},\,  B_{t_b}\geq -\mathtt{z}-M_1\}}  \1_{\{ \overline{B}_{t_b }\leq 0,\, \overline{B}_{[\frac{t}{2},t_b] } \leq -\mathtt{z} \}}  \Big| F^{ (M)}(B_{t_b} +\mathtt{z},\mathfrak{G}_{t,b,\sigma})-F (B_{t_b} +\mathtt{z},\mathfrak{G}_{t,b} ) \Big|,
\end{equation}
with $ t>b  $. Now we will choose $M>M_1$. On   $ \{ w^{(0,1)}_{\mathfrak{G}_{t,b,\sigma}(\cdot \ee^b)}(\frac{1}{M_1})\leq \frac{1}{4},\,  B_{t_b} +\mathtt{z}\geq  -M_1\} $, by the properties (\ref{et2rajout}), (\ref{et1}) and (iii) of $F$, we get  
\begin{equation}
\label{useful}
F ( B_{t_b}+\mathtt{z},\mathfrak{G}_{t,b,\sigma} )\leq h(B_{t_b} +\mathtt{z})F^*(\mathfrak{G}_{t,b,\sigma}) \leq c M_1^{10}:=M .
\end{equation}
 Then (\ref{restrain}) is equal to 
\begin{equation}
\label{restrain2}
\1_{\{ w^{(0,1)}_{\mathfrak{G}_{t,b,\sigma}(\cdot \ee^b)}(\frac{1}{M_1})\leq \frac{1}{4},\,  B_{t_b}\geq -\mathtt{z}-M_1\}}  \1_{\{ \overline{B}_{t_b }\leq 0,\, \overline{B}_{[\frac{t}{2},t_b] } \leq -\mathtt{z} \}}  \Big|F (B_{t_b} +\mathtt{z},\mathfrak{G}_{t,b,\sigma})-F (B_{t_b} +\mathtt{z},\mathfrak{G}_{t,b} ) \Big|\wedge 2M.
\end{equation}
We denote  $ || \Delta \mathfrak{G}_{\sigma}||_{\infty}:= \underset{y\in B(0,e^b)}{\sup} \Big|\mathfrak{G}_{t,b} (y)-\mathfrak{G}_{t,b,\sigma}(y) \Big| $, by the property (\ref{et3bibi}) of $F$, for any $\delta>0$,  we deduce that (\ref{restrain2}) is smaller than
\begin{eqnarray*}
&& \left|F_{L,b}(B_{t_b}+\mathtt{z},\mathfrak{G}_{t,b,\sigma})-F_{L,b}(B_{t_b}+\mathtt{z},\mathfrak{G}_{t,b} ) \right|\wedge 2M
 \\
&& \leq    2M \1_{\{ || \Delta \mathfrak{G}_{\sigma}||_{\infty}\geq \delta^4 \}}   +      \1_{\{  || \Delta \mathfrak{G}_{\sigma}||_{\infty}\leq \delta^4 \}} \delta  h_{L,b}(B_{t_b}+\mathtt{z})   F^*_b(\mathfrak{G}_{t,b,\sigma}(\cdot)).
\end{eqnarray*}
As $w^{(0,1)}_{\mathfrak{G}_{t,b,\sigma}(\cdot\ee^b)}(\frac{1}{M_1})\leq \frac{1}{4} $ and $B_{t_b}\geq -\mathtt{z}-M_1 $ we now use (\ref{useful}) to bound (\ref{restrain2}) by
\begin{equation}
 \label{takeexps}\1_{\{  B_{t_b}\geq -\mathtt{z}-M_1\}}\1_{\{   \overline{B}_{[\frac{t}{2},t_b] } \leq -\mathtt{z} \}}     \left(2M\1_{\{|| \Delta \mathfrak{G}_{\sigma}||_{\infty}\geq \delta^4 \}}   + M\delta \1_{\{  || \Delta \mathfrak{G}_{\sigma}||_{\infty}\leq \delta^4 \}}\right)
\end{equation}
Now we claim the following two assertions:

{\it -For any $L,\, b,\delta,\,M_1$ there exists $ T>0$ such that for any $t\geq T$, $\alpha \in [1,\log t]$, $\mathtt{z}\leq (\log t)^{30}$ we have
\begin{equation}
\label{eqcontroldiff01}
 \P_{-\alpha} \Big(   \overline{B}_{t_b }\leq 0,\, \overline{B}_{[\frac{t}{2},t_b] } \leq -\mathtt{z},\,  B_{t_b} +\mathtt{z}\geq  -M_1  \Big)\leq c_{12} \frac{ \alpha}{t^{\frac{3}{2}}}(1+M_1)^2,
\end{equation}

-For any $L,\, b,\delta,\,M_1$ there exists $\sigma,\, T>0$ such that for any $t\geq T$, $\alpha \in [1,\log t]$, $\mathtt{z}\leq (\log t)^{30}$ we have
\begin{equation}
\label{eqcontroldiff}
 \P_{-\alpha} \Big( || \Delta \mathfrak{G}_{\sigma}||_{\infty}\geq \delta,\,   \overline{B}_{t_b }\leq 0,\, \overline{B}_{[\frac{t}{2},t_b] } \leq -\mathtt{z},\,  B_{t_b} +\mathtt{z}\geq  -M_1  \Big)\leq c \frac{ \alpha}{t^{\frac{3}{2}}}(1+M_1)^2\exp(- \frac{c_{19}}{2}\delta^2 e^{2\sigma}).
\end{equation} }
So we take the expectation 

To conclude we notice that the assertion (\ref{eqcontroldiff01}) comes from (\ref{trajA.7}), whereas  (\ref{eqcontroldiff}) is a consequence of (\ref{trajB.20}). Indeed it suffices to notice that:
\begin{eqnarray*}
\{ || \Delta \mathfrak{G}_{\sigma}||_{\infty}\geq \delta \}&=&\{ \underset{y\in B(0,\ee^b)}{\sup}|\int_0^{t_b-\sigma} \mathtt{g}(\ee^{s-t} y)dB_s |\geq \delta\}
\\
&=&\{ \underset{|y|\leq \ee^{-t_b}}{\sup}\int_0^{t_b-\sigma} \mathtt{g}(\ee^s y)dB_s |\geq \delta\}\subset A_{t_b, t_b-\sigma, \delta}.
\end{eqnarray*}
\hfill$\Box$
\\

Now we tackle the proof of (\ref{eqpassageaueps2}). Let us introduce some notations:

- Let $(\mathtt{R}_{s})_{s \geq 0}$ be a three dimensional Bessel process starting from $0$.

- Let $(B_s)_{s\geq 0}$ be real Brownian motion and for any $\sigma>0 $ we denote $(B_s^{(\sigma)})_{s\geq 0}:= (B_{s+\sigma}-B_\sigma)_{s\geq 0} $.  \nomenclature[a4]{$(\mathtt{R}_s)_{s\geq 0}$}{: a three dimensional Bessel process starting from $0$}

- Let  $g$, $h$ be two processes, for any $ t_0\in \r^+$ the process $\mathtt{X}_\cdot(t_0,g,h)$ is defined by \nomenclature[a5]{$ \mathtt{X}_s(t_0,g,h)$}{: the concatenation at time $t_0$ of processes $ g $ and  $h  $}
\begin{eqnarray}
\label{defXX} \mathtt{X}_s(t_0,g,h)=\left\{\begin{array}{ll} g_s ,&\qquad \text{if  } s\leq t_0,
\\
g_{t_0}+ h_{t-t_0},&\qquad \text{if  } s\geq t_0.
\end{array}
 \right. 
\end{eqnarray}

- Let $\sigma>0 $ for any process $(g_s)_{s\leq \sigma}$ we set\nomenclature[a6]{$(\overset{\leftarrow \sigma}{g_s})_{s\leq \sigma} $}{$:=  (g_{\sigma-s}-g_\sigma)_{s\leq \sigma} $}
\begin{equation}
\label{defarrw} (\overset{\leftarrow \sigma}{g_s})_{s\leq \sigma}:= (g_{\sigma-s}-g_\sigma)_{s\leq \sigma}.
\end{equation}

\nomenclature[b97]{$\mathfrak{H}_{m,\sigma}$}{: the set of continuous function $F:\r\times \mathcal{C}([0,\sigma],\r)\to \r^+$ with $\underset{u\in \r,\, g\in \mathcal{C}([0,\sigma ],\r))}{\sup} F(u,g)\leq m$}
- We set  $\mathfrak{H}_{m,\sigma}$ the set of continuous functions $F:\r\times \mathcal{C}([0,\sigma],\r)\to \r^+$ with $\underset{u\in \r,\, g\in \mathcal{C}([0,\sigma ],\r))}{\sup} F(u,g)\leq m$. For $ g\in \mathcal{C}^1(\r^\d,\r) $ we denote by $ \nabla_{y}(g) $ the gradient of $g$ at $y\in \r^d $. At last we denote by $ \langle\cdot,\cdot\rangle $  the inner product in $\r^\d $.
 \\

%
%
%
%
%
%

Display (\ref{eqpassageaueps2}) is a consequence of the following Proposition which is proven in the Appendix.
\nomenclature[d3]{$C$}{$ := 2^\frac{3}{2}  \E(\ee^{- {\mathtt{R}^2_{1}} })$}
\begin{Proposition}
\label{renouv}
Let $B$ be a Brownian motion and let $ \mathtt{R}$ be a three dimensional Bessel process starting from $0$ independent of $B$. Let $m,\, \sigma \geq 0$ be two constants. For any $\epsilon>0$ there exists $T(m,\sigma ,\epsilon)\geq 0$ such that for any $t\geq T$, $1\leq \alpha,\,  \mathtt{z}\leq (\log t)^{30} $ and $F\in \mathfrak{H}_{m,\sigma}$
\begin{eqnarray}
\label{eqrenouv}
 \left| \frac{t^\frac{3}{2}}{\alpha}\E_\alpha\left(\1_{\{ \underline{B}_t\geq 0,\, \underline{B}_{[\frac{t}{2},t]}\geq \mathtt{z},\, B_t-\mathtt{z} \leq m\}}F\left( B_t-\mathtt{z}, (B^{(t-\sigma)}_{s})_{s\leq \sigma}\right)\right)-\right.\qquad\qquad\qquad 
 \\
\nonumber \left.   \sqrt{\frac{2}{\pi}}  \int_{0}^m  \int_{0}^u \E\left(     F(u, (\overset{\leftarrow \sigma}{\mathtt{X}_s}(T_{-\gamma}, B, \mathtt{R}))_{s\leq \sigma} )\right)d\gamma du    \right|\leq \epsilon    ,
\end{eqnarray}
where  $T_\gamma:=\inf\{ s\geq 0,\, B_s=\gamma\}$, $\gamma \in \r $.
\end{Proposition}

\medskip 

\noindent{\it  Proof of  \eqref{eqpassageaueps2}.} Fix $b,\, M,\, \sigma>0$ and $F$ a function {\bf $b$ regular}. Let us explicit the expectation in (\ref{eqpassageaueps2}). As $(B_s)_{s\geq 0}\overset{law}{=}(-B_s)_{s\geq 0} $ we have,
\begin{eqnarray*}
\E(\ref{eqpassageaueps2})&:=&\E_{-\alpha}\Big[\1_{\{\overline{B}_{t_b}\leq 0,\,\overline{B}_{[\frac{t}{2},t_b]} \leq -\mathtt{z}\}}F^{(M)}\Big(B_{t_b} +\mathtt{z},\mathfrak{G}_{t,b,\sigma}\Big)\Big]
\\
&=& \E_{\alpha}\Big[\1_{\{ \underline{B}_{t_b} \geq 0,\, \underline{B}_{[\frac{t}{2},t_b]}\geq \mathtt{z},\, B_{t_b}-\mathtt{z}\leq M \}} F\Big(-[B_{t_b}-\mathtt{z}], 
\\
&&\qquad\qquad\qquad \qquad\qquad\qquad  y\mapsto \int_{t_b-\sigma}^{t_b}(1-\mathtt{k}(e^{s-t}y))dB_s-\zeta_{t_b}(ye^{-t})+Z^0_{t_b}(ye^{-t})    \Big) \wedge M\Big].
\end{eqnarray*}
Moreover by integration by parts, the second argument of the function in $F$ can be rewritten as:
\begin{eqnarray*}
y\mapsto \big(1-\mathtt{k}(e^{-b}y)\big)\Big[B_{t_b}-B_{t_b-\sigma}  \Big] +   \int_{t_b-\sigma}^{t_b} [B_s-B_{t_b-\sigma}]\langle \nabla_{ye^{s-t}} k\cdot ye^{s-t}\rangle ds    -\zeta_{t_b}(ye^{-t})+Z^0_{t_b}(ye^{-t}),
\end{eqnarray*}
and we recall that the processes $B $ and $Z$ are independent. So $ \E(\ref{eqpassageaueps2}) $ is equal to 
\begin{eqnarray*}
\E_{\alpha}\Big[\1_{\{  \underline{B}_{t_b} \geq 0,\, \underline{B}_{[\frac{t}{2},t]} \geq \mathtt{z},\, B_{t_b}- \mathtt{z}\leq M\}}\Phi_{t_b}(B_{t_b}-\mathtt{z},\, (B_s^{(t-\sigma)})_{s\leq \sigma})\Big],
\end{eqnarray*}
with $ \Phi_{t_b}: \r\times \mathcal{C}([0,\sigma],\r) \to \r^+ $ , a continuous function, bounded by $M$ and defined by
\begin{eqnarray*}
(u,h) \mapsto \E\Big[ F(-u,y\mapsto   \big(1-\mathtt{k}(e^{t_b-t}y)\big)[h_{\sigma}-h_{0}]+\int_{t_b-\sigma}^{t_b} [h_{s-(t_b-\sigma)}-h_{0}]  \langle \nabla_{ye^{s-t}}k. ye^{s-t}\rangle  ds \Big.\qquad
\\
\Big. - \zeta_{t_b}(ye^{-t})+ Z^0_{t_b}(ye^{-t}) )\wedge M\Big].
\end{eqnarray*}
 Now we can apply Proposition \ref{renouv}, with  $t \leftrightarrow t_b>0$,  $\alpha    \leftrightarrow  \alpha ,\,  \mathtt{z} \leftrightarrow  \mathtt{z}\leq (\log t)^{30},\, \sigma \leftrightarrow  \sigma,\, m  \leftrightarrow M$ and  $F\leftrightarrow \Phi_{t_b}  $. Then for any $\epsilon>0$ there exists $T\geq 0$ such that for any $t\geq T$, $1\leq \alpha\leq (\log t)^{30}$
\begin{eqnarray}
\label{1(ii)}
\Big[|\frac{t^\frac{3}{2}}{\alpha}\E(\ref{eqpassageaueps2})- \sqrt{\frac{2}{\pi}}  \int_{0}^m  \int_{0}^u \E\Big[(     \Phi_{t_b}(u, (\overset{\leftarrow \sigma}{\mathtt{X}_s}(T_{-\gamma}, B, \mathtt{R}))_{s\leq \sigma} )\Big)d\gamma du    \Big|\leq \epsilon   .
\end{eqnarray}
Moreover, we observe that for any  $u>0,\gamma\leq u $,  
\begin{eqnarray*} 
&& \E\Big[(\Phi_{t_b}(u, \overset{\leftarrow \sigma}{\mathtt{X}}_s(T_{-\gamma}, B, R))_{s\leq \sigma} )\Big]= \E\Big[ F^{(M)}\Big(-u,  y\mapsto Z_{t_b}(ye^{-t})-\zeta_{t_b}(ye^{-t})
 \\
&& \qquad\qquad\qquad\qquad \qquad\qquad\qquad   -\int_{0}^{T_{-\gamma} \wedge \sigma }\big(1-\mathtt{k}(e^{-s}ye^{-b})\big) dB_s- \int_{T_{-\gamma}\wedge \sigma}^\sigma\big(1-\mathtt{k}(e^{-s}ye^{-b})\big) d \mathtt{R}_{s-T_{-\gamma}}\Big)\Big].
\end{eqnarray*}
Finally as $ (Z_{t_b}(ye^{-t})-\zeta_{t_b}(ye^{-t}))_{y\in B(0,e^{b})} $ is independent of $(B,R)$ and converges in law, when $t$ goes to infinity, to $ (Z(ye^{-b})-\zeta(ye^{-b}))_{y\in B(0,e^{b})} $ (see (\ref{Zlimit})), by combining with (\ref{1(ii)}) we deduce that: for any $\epsilon>0$ there exists $T\geq 0$ such that for any $t\geq T$, $1\leq \alpha\leq (\log t)^{30}$
\begin{eqnarray}
\Big|\frac{1}{\alpha t^{\frac{3}{2}}} \E_{-\alpha}\Big[ \1_{\{\overline{B}_{t_b}\leq 0,\, \overline{B}_{[\frac{t}{2},t_b]} \leq -\mathtt{z}\}} F^{(M)}\big(B_{t_b}+\mathtt{z},\mathfrak{G}_{t,b,\sigma}\big)\Big]-C_{M,\sigma}(F)\Big|\leq \epsilon,
\end{eqnarray}
with 
\begin{eqnarray}
\nonumber C_{M,\sigma}(F):=  \sqrt{\frac{2}{\pi}}\int_{0}^M\int_0^u\E\Big[  F^{(M)}\Big(-u, y\mapsto Z (ye^{-b})-\zeta(ye^{-b}) \qquad\qquad\qquad\qquad\qquad
\\
  -\int_{0}^{T_{-\gamma} \wedge \sigma }\big(1-\mathtt{k}(e^{-s}ye^{-b})\big) dB_s - \int_{T_{-\gamma}\wedge \sigma}^\sigma \big(1-\mathtt{k}(e^{-s}ye^{-b})\big) d \mathtt{R}_{s-T_{-\gamma}}\Big)\Big]d\gamma dU.
\end{eqnarray}
This completes the proof of (\ref{eqpassageaueps2})  \hfill $\Box $

\section*{Appendix}

\appendix

\section{Proof of the Proposition \ref{renouv}}
In the following we denote $\mathtt{R}_\cdot^{(x)}$ a Bessel three process starting from $x\geq 0$ ($\mathtt{R}_\cdot^{(0)}=\mathtt{R}_\cdot$). Our aim here, is to prove the Proposition \ref{renouv}. First let us state two results:

\begin{Proposition}[pp 255 in \cite{RYor99}]
\label{Yor1}
Let $ \mathtt{R}^{(x)}$ be a three dimensional Bessel process starting from $x>0$ and $\tau:= \inf\{s\geq 0,\,  \mathtt{R}^{(x)}_s=\underset{u\geq 0}{\inf} \mathtt{R}^{(x)}_u\}  $; the process $( \mathtt{R}^{(x)}_s,\, s\leq \tau)$  has the same law as $(B_t,\, t\leq T_\kappa)$, where $B$ is a Brownian motion starting from $x>0$ and $T_\kappa$ is the hitting time by $B$ of an independent random point $\kappa$ uniformly distributed on $[0,x]$. Moreover conditioned on $\{\mathtt{R}^{(x)}_{\tau}=y\} $, $(\mathtt{R}^{(x)}_{\tau+s}-y)_{ s\geq 0} $ is a three dimensional Bessel process starting from $0$ independent of $(\mathtt{R}^{(x)}_s)_{s\leq \tau} $.
\end{Proposition}

\begin{Lemma}
\label{proxform}
 Let $m>0$ and $\lambda_0\geq 0$. For any $\epsilon>0$ there exists $T(m,\lambda_0,\epsilon)>0$ such that for any $t\geq T$, $  b\geq t^{\frac{1}{4}}$, $\gamma\in [0,m] $,  $\lambda\leq \lambda_0$ and  $F\in \mathfrak{H}_{m,\lambda}$
\begin{equation}
\label{eqproxform}  \left| t^\frac{3}{2}\E\left( \frac{F(b-\mathtt{R}_t, ( \mathtt{R}_l)_{l\leq \lambda})}{ \mathtt{R}_t+\gamma}\1_{\{ \gamma \geq \mathtt{R}_t-b+m\geq 0\}}\right)-   \E\left( \int_{m-\gamma}^m   F(u,( \mathtt{R}_l)_{l\leq \lambda})du\right)\sqrt{\frac{2}{\pi}}b\ee^{-\frac{b^2}{2t}}\right|\leq  \epsilon + \epsilon b\ee^{-\frac{b^2}{2t}}.
\end{equation}
\end{Lemma}

\noindent{\it Proof of Lemma \ref{proxform}.} This is a slight extension of the local limit theorem for the three dimensional Bessel processes. Indeed let us assume that $F(u,g)=F(u)$. Recall that for a three dimensional Bessel process starting from $0$, $\P\left( \mathtt{R}_t\in dx\right)=\sqrt{\frac{2}{\pi t^3}}x^2\exp(-\frac{x^2}{2t})dx$, moreover as $b \gg m, \gamma$, we have
\begin{eqnarray}
\nonumber t^{\frac{3}{2}}\E\Big(F(b-\mathtt{R}_t)\frac{1_{\{\gamma\geq \mathtt{R}_t-b+m\geq 0\}}}{\mathtt{R}_t+\gamma} \Big) &=&  \sqrt{\frac{2}{\pi}} \int_{-m+b}^{\gamma-m+b} \frac{x^2}{x+\gamma}\exp(-\frac{x^2}{2t})F(b-x) dx
\\
\nonumber &=& \sqrt{\frac{2}{\pi}} \int_{m-\gamma}^m F(u)\frac{(b-u)^2}{b-u+\gamma} \exp(-\frac{(b-u)^2}{2t})du
\\
\label{fromzero} &=& \sqrt{\frac{2}{\pi}} b\ee^{-\frac{b^{2}}{2t}} \int_{m-\gamma}^m F(u) du (1+o(1)),
\end{eqnarray}
which proves (\ref{eqproxform}).  Of course, for any $A>0$, display (\ref{fromzero}) remains true for $ (\mathtt{R}^{(x)}_s)_{s\geq 0}$ uniformly in $x\in [0,A]$. Now let us prove (\ref{eqproxform}) for any function $F\in \mathfrak{H}_{m,\lambda}$. According to the Markov property at time $\lambda$,
\begin{eqnarray}
\nonumber t^\frac{3}{2}\E\left( \frac{F(b-\mathtt{R}_t, ( \mathtt{R}_l)_{l\leq \lambda})}{ \mathtt{R}_t+\gamma}\1_{\{ \gamma \geq \mathtt{R}_t-b+m\geq 0\}}\right)= \E\left( t^{\frac{3}{2}}  \E_{\mathtt{R}_\lambda} \left(    \frac{F(b-\mathtt{R}_{t-\lambda}, (g_l)_{l\leq \lambda}) }    {\mathtt{R}_{t-\lambda}+\gamma}\1_{\{  \gamma \geq \mathtt{R}_{t-\lambda}-b+m\geq 0\}}  \right)_{(g_l)_{l\leq \lambda}=(\mathtt{R}_l)_{l\leq \lambda}}  \right).
\end{eqnarray}
By letting $t$ going to infinity an applying (\ref{fromzero}) we obtain easily Lemma (\ref{proxform}).
\hfill$\Box$
\\

We can start the proof of the \noindent{\it Proof of Proposition \ref{renouv}.} Let $\tau_0:= t -\inf\{ s\geq 0,\, B_{t-s}=\underline{B}_{[\frac{t}{2},t]}\}$ and  $E_t:=\{t\geq B_{\frac{t}{2}}\geq t^{\frac{1}{4}}\}\cap \{ \tau_0 \geq t-\log t\}$. First we show that  uniformly in $1\leq \alpha,\, \mathtt{z}\leq (\log t)^{30}$, 
\begin{equation}
\label{choubka2}   \frac{t^\frac{3}{2}}{\alpha}\E_\alpha\left(m\1_{\{ \underline{B}_t\geq 0,\, \underline{B}_{[\frac{t}{2},t]}\geq \mathtt{z},\, B_t-\mathtt{z}\leq m\}} \1_{E^c_t}\right)= o(1),\qquad t\to \infty.
\end{equation}}
Clearly $  m\P_\alpha\left( B_{t/2}\geq t\right)=o(t^\frac{-3}{2}) $, thus to prove (\ref{choubka2}), it remains to study\\ $ \E_\alpha\left(m\1_{\{ \underline{B}_t\geq 0,\, \underline{B}_{[\frac{t}{2},t]}\geq   \mathtt{z},\, B_t-\mathtt{z}\leq m\}} \1_{\{B_{\frac{t}{2}}\leq t^{\frac{1}{4}}\}}\right) $ then $ \E_\alpha\left(\1_{\{ \underline{B}_t\geq 0,\, \underline{B}_{[\frac{t}{2},t]}\geq   \mathtt{z},\, B_t-  \mathtt{z}\leq m\}} \1_{\{\tau_0 \leq t-\log t\}}\right) $. By the Markov property at time $\frac{t}{2}$ then (\ref{2.6}) we get that
\begin{eqnarray}
\nonumber \frac{t^\frac{3}{2}}{\alpha}\E_\alpha\left(m\1_{\{ \underline{B}_t\geq 0,\, \underline{B}_{[\frac{t}{2},t]}\geq  \mathtt{z},\, B_t-\mathtt{z}\leq m\}} \1_{\{B_{\frac{t}{2}}\leq t^{\frac{1}{4}}\}}\right)  &\leq& \frac{t^\frac{3}{2}}{\alpha}\E_\alpha\left(m\1_{\{ \underline{B}_t\geq 0,\, B_{\frac{t}{2}}\leq t^{\frac{1}{4}}\}}\P_{B_{\frac{t}{2}}-\mathtt{z}}\left( \underline{B}_{\frac{t}{2}}\geq 0,\, B_{\frac{t}{2}}\leq m  \right)\right)
\\
\nonumber &\leq& c\frac{t^\frac{3}{2}}{\alpha}\E_\alpha\left(m\1_{\{ \underline{B}_t\geq 0,\, B_{\frac{t}{2}}\leq t^{\frac{1}{4}}\}} t^{\frac{1}{4}-\frac{3}{2}} \right)
\\
\label{choubkkka1} &\leq& c\frac{(1+\alpha)}{\alpha}m t^{\frac{3}{4}-\frac{3}{2}}=o(1).
\end{eqnarray}
Set $\tau_1:= \inf\{ s\geq 0,\, B_{s}=\underline{B}_{[\frac{t}{2},t]}\} $. By the Markov property at time $\frac{t}{3}$, then the property of time reversal of the Brownian motion we get that
\begin{eqnarray*}
&& \E_\alpha\left(\1_{\{ \underline{B}_t\geq 0,\, \underline{B}_{[\frac{t}{2},t]}\geq \mathtt{z},\, B_t-\mathtt{z}\leq m\}} \1_{\{\tau_0 \leq t-\log t\}}\right)
\\
&&\leq \E_\alpha\left(\1_{\{ \underline{B}_{\frac{t}{3}}\geq 0\}}\E\left(\1_{\{ \underline{B}_{\frac{2t}{3}}\geq -m,\,  0\leq x+B_{\frac{2t}{3}}-\mathtt{z}\leq m,\, \tau_1 \geq   \log t\}}\right)_{\Big| x=B_{\frac{t}{3}}}\right).
\end{eqnarray*}
In the second expectation, we apply the Markov property at time $\frac{t}{3}$, then the inequalities (\ref{2.5}) and (\ref{2.6}) and we get that
\begin{eqnarray}
 \nonumber \frac{t^\frac{3}{2}}{\alpha} \E_\alpha\left(\1_{\{ \underline{B}_t\geq 0,\, \underline{B}_{[\frac{t}{2},t]}\geq \mathtt{z},\, B_t-\mathtt{z}\leq m\}} \1_{\{\tau \leq t-\log t\}}\right)\leq \frac{t^\frac{3}{2}}{\alpha} \P_\alpha\left( \underline{B}_{\frac{t}{3}}\geq 0 \right) \underset{x\in \r}{\sup}\P\left( 0\leq x+B_{\frac{t}{3}}\leq m\right) \times
\\
\nonumber\P\left(  \underline{B}_{\frac{t}{3}}\geq -m,\, \inf\{ s\geq 0,\, B_{s}=\underline{B}_{[0,\frac{t}{3}]}\}  \geq   \log t\right) 
 \\
\label{choubkkka2} \leq \frac{t^\frac{3}{2}}{\alpha} \frac{c(1+\alpha)}{t^{\frac{3}{2}}\sqrt{\log t}}=o(1). \qquad\qquad\qquad\qquad\qquad\qquad 
\end{eqnarray}
From (\ref{choubkkka1}) and (\ref{choubkkka2}) on has (\ref{choubka2}). So we can restrict our study to
\begin{equation}
\label{ExpD}\E_\alpha\left(\1_{\{ \underline{B}_t\geq 0,\, \underline{B}_{[\frac{t}{2},t]}\geq \mathtt{z},\, B_t-\mathtt{z}\leq m\}}F\left( B_t-\mathtt{z}, (B^{(t-\sigma)}_{s})_{s\leq \sigma}\right)\1_{E_t}\right).
\end{equation}
By the Markov property at time $\frac{t}{2}$, the property of time reversal of $(B_s)_{s\leq \frac{t}{2}}$ and the equality $(B_s)_{s\leq \frac{t}{2}}\overset{(d)}{=}(-B_s)_{s\leq \frac{t}{2}}$,  the expectation in (\ref{ExpD}) is equal to
\begin{equation}
\label{premierMarkov}
\E_\alpha\left(\1_{\{\underline{B}_{\frac{t}{2}}\geq 0,\,t\geq  B_{\frac{t}{2}}\geq t^{\frac{1}{4}}\}} \Phi_{\mathtt{z},m,\sigma,t}(B_{\frac{t}{2}})\right) ,
\end{equation}
with $\Phi_{\mathtt{z},m,\sigma,t}(x)$ defined by (recall that $ \tau_1:=  \inf\{ s\geq 0,\, B_{s}=\underline{B}_{\frac{t}{2}}\}$)
\begin{equation}
\Phi_{\mathtt{z},m,\sigma,t}(x)=\E\left( F\left(-B_{\frac{t}{2}}+x-\mathtt{z},(B_{\sigma-s}-B_\sigma)_{s\leq \sigma}\right)\1_{\{  \underline{B}_\frac{t}{2}\geq B_\frac{t}{2}+ \mathtt{z}-x\geq -m,\, \tau_1\leq \log t\}} \right).
\end{equation}
 Using the notation (\ref{defarrw}) we get that
\begin{eqnarray*}
\Phi_{\mathtt{z},m,\sigma,t}(x) &=&\E_m\left( \frac{B_{\frac{t}{2}}}{B_{\frac{t}{2}}} F\left(m-B_{\frac{t}{2}}+x-\mathtt{z},(\overset{\leftarrow \sigma}{B}_s)_{s\leq \sigma} \right)\1_{\{  \underline{B}_\frac{t}{2}\geq B_\frac{t}{2}+ \mathtt{z}-x\geq 0,\, \tau_1\leq \log t\}} \right).
\end{eqnarray*}
We recognize the $h-$transform of the Bessel process, therefore with $\tau_2:= \inf\{ s\geq 0,\, \mathtt{R}_{s}=\underline{\mathtt{R}}_{\frac{t}{2}}\}$, we have
\begin{eqnarray*}
\Phi_{\mathtt{z},m,\sigma,t}(x) &=&\E_m\left( \frac{m}{ \mathtt{R}_{\frac{t}{2}}} F\left(m- \mathtt{R}_{\frac{t}{2}}+x-\mathtt{z},(\overset{\leftarrow \sigma}{R}_s)_{s\leq \sigma} \right)\1_{\{  \underline{R}_\frac{t}{2}\geq  \mathtt{R}_\frac{t}{2}+ \mathtt{z}-x\geq 0,\, \tau_2\leq \log t\}} \right).
\end{eqnarray*}
We define  $\tau:= \inf\{ s\geq 0,\, \mathtt{R}_{s}=\underline{\mathtt{R}}_{\infty}\}$. Observe that 
\begin{eqnarray}
\nonumber \left|\Phi_{\mathtt{z},m,\sigma,t}(x)-\E_m\left( \frac{m}{ \mathtt{R}_{\frac{t}{2}}} F\left(m- \mathtt{R}_{\frac{t}{2}}+x-\mathtt{z},(\overset{\leftarrow \sigma}{R}_s)_{s\leq \sigma} \right)\1_{\{  \underline{R}_\frac{t}{2}\geq  \mathtt{R}_\frac{t}{2}+ \mathtt{z}-x\geq 0,\,\tau_2= \tau\leq \log t\}} \right)\right|
\\
\label{thlkmlk} \leq   \E_m\left( \frac{m^2}{ \mathtt{R}_{\frac{t}{2}}}  \1_{\{  \underline{R}_\frac{t}{2}\geq  \mathtt{R}_\frac{t}{2}+ \mathtt{z}-x\geq 0,\, \tau\geq \frac{t}{2} \}} \right).
\end{eqnarray}
According to the Markov property at time $\frac{t}{2}$ and Proposition \ref{Yor1}, the expectation in the right hand side of (\ref{thlkmlk}) is smaller than:
\begin{eqnarray}
\nonumber  \E_m\left( \frac{m^2}{ \mathtt{R}_{\frac{t}{2}}}  \1_{\{  \underline{R}_\frac{t}{2}\geq  \mathtt{R}_\frac{t}{2}+ \mathtt{z}-x\geq 0\}} \P_{\mathtt{R}_{\frac{t}{2}}}\left( \underset{s\geq 0}{\inf }R_s\leq m \right) \right)&= & \E_m\left( \frac{m^3}{ \mathtt{R}_{\frac{t}{2}}^2}  \1_{\{  \underline{R}_\frac{t}{2}\geq  \mathtt{R}_\frac{t}{2}+ \mathtt{z}-x\geq 0\}} \right)
\\
\label{sisimseslm} &\leq & c t^{-\frac{3}{2}} \ee^{-\frac{x^2}{t}}=o(t^{-\frac{3}{2}} x\ee^{-\frac{x^2}{t}}).
\end{eqnarray}
Furthermore in order to use Lemma \ref{proxform}, we disintegrate the expectation in the left and side of (\ref{thlkmlk}) with respect of $\tau$ and $\underline{R}_{\infty}$ and apply Proposition \ref{Yor1}. Finally with (\ref{sisimseslm}), it stems
\begin{equation}
\label{Phiii}
\Phi_{\mathtt{z},m,\sigma,t}(x)=\int_{0}^m \E\left(\frac{F\left((\ref{Phiii})_{(a)},(\ref{Phiii})_{(b)} \right)}{ \mathtt{R}_{\frac{t}{2}-T_{\gamma-m}}+\gamma} \1_{\{ \gamma \geq  \mathtt{R}_{\frac{t}{2}-T_{\gamma-m}}+ \gamma+\mathtt{z}-x\geq 0,\, T_{\gamma-m}\leq \log t\}} \right)d\gamma+ o(t^{-\frac{3}{2}} x\ee^{-\frac{x^2}{t}}),
\end{equation}
with $T_{\gamma -m}:= \inf\{ s\geq 0,\, B_{s}=\gamma-m \}$ for $B$ a Brownian motion independent of  $ \mathtt{R}$ and
\begin{eqnarray*}
(\ref{Phiii})_{(a)}&:=&m- \mathtt{R}_{\frac{t}{2}-T_{\gamma-m}}-\gamma +x-\mathtt{z},
\\
(\ref{Phiii})_{(b)}&:=&( \overset{\leftarrow \sigma}{\mathtt{X} _s}(T_{\gamma-m}, B, \mathtt{R}))_{s\leq \sigma}.
\end{eqnarray*}
So from (\ref{Phiii}) we can write
\begin{eqnarray*}
\Phi_{\mathtt{z},m,\sigma,t}(x)&=&\int_{0}^m \E\Big( \E(...)_{\Big|...}\1_{\{ T_{\gamma-m}\leq \log t\}}    \Big)d\gamma + o(t^{-\frac{3}{2}} x\ee^{-\frac{x^2}{t}}),\qquad \text{with}
\\
\E(...)_{\Big|...}& =& \E\Big( \frac{F\left(b- \mathtt{R}_{\mathtt{t}}, ( \overset{\leftarrow \sigma}{ \mathtt{X}_s}(\lambda, g, R))_{s\leq \sigma} \right)}{ \mathtt{R}_{\mathtt{t}}+\gamma} \1_{\{ \gamma \geq  \mathtt{R}_{\mathtt{t}}+ \gamma+\mathtt{z}-x\geq 0,\}}\Big)_{\Big| \begin{array}{ll} \mathtt{t}=\frac{t}{2}-T_{\gamma-m} 
\\
b=m-\gamma+x-\mathtt{z}
\\
g=B,\, \lambda= T_{\gamma-m}
\end{array}}.
\end{eqnarray*}
For any $g\in \mathcal{C}(\r^+,\r)$, $\lambda>0$, $F\left(b- \mathtt{R}_{\mathtt{t}}, ( \overset{\leftarrow \sigma}{\mathtt{X}_s}(\lambda, g, \mathtt{R}))_{s\leq \sigma} \right)$ can be rewritten  $  F_{\lambda,\sigma,g}\left(b- \mathtt{R}_{\mathtt{t}}, ( R)_{s\leq (\sigma-\lambda)_+} \right)$  with 
$F_{\lambda,\sigma,g}$ a function in $  \mathfrak{H}_{m,(\sigma-\lambda)_+} $. So we can apply Lemma \ref{proxform} to $ \E(...)_{\Big|...}$, it allows us to affirm that
\begin{eqnarray*}
2^\frac{3}{2} (\frac{t}{2})^\frac{3}{2}\Phi_{\mathtt{z},m,\sigma,t}(x)= \frac{4}{\sqrt{\pi}}\int_0^m \int_{m-\gamma}^m\E\Big(  F\left(u, ( \overset{\leftarrow \sigma}{\mathtt{X} _s}(T_{\gamma-m}, B, \mathtt{R}))_{s\leq \sigma} \right)  (x+m-\gamma -\mathtt{z})\qquad\qquad\qquad 
\\
\ee^{\frac{-(x+m-\gamma-\mathtt{z})^2}{2[\frac{t}{2}-T_\gamma]}}  \1_{\{ T_{\gamma-m}\leq \log t\}}    \Big) du   d\gamma +o(1+x\ee^{-\frac{x^2}{t}}),
\end{eqnarray*}
Recall that for any $\mathtt{z}\leq (\log t)^{30},\, x\geq t^{\frac{1}{4}},\, T_\gamma\leq \log t $,
\begin{eqnarray*}
(x+m-\gamma -\mathtt{z})\ee^{-\frac{(x+m-\gamma-\mathtt{z})^2}{2[\frac{t}{2}-T_\gamma]}} = x\ee^{-\frac{x^2}{t }} +o(x\ee^{-\frac{x^2}{t}}).
\end{eqnarray*} 
Finally we get that $t^{\frac{3}{2}}\Phi_{\mathtt{z},m,\sigma,t}(x) $ is equal to
\begin{eqnarray}
\nonumber   =   \frac{4}{\sqrt{\pi}} x\ee^{-\frac{x^2}{t }}  \int_0^m \int_{m-\gamma}^m\E\left(  F\left(u, ( \overset{\leftarrow \sigma}{\mathtt{X} _s}(T_{\gamma-m}, B, \mathtt{R}))_{s\leq \sigma} \right)     \1_{\{ T_{\gamma-m}\leq \log t\}}    \right) du  d\gamma   + o(1+x\ee^{-\frac{x^2}{t}} )
\\
\label{part1} =  \frac{4}{\sqrt{\pi}} x\ee^{-\frac{x^2}{t}}  \int_{0}^m  \int_{0}^u \E\left(     F(u, \overset{\leftarrow \sigma}{\mathtt{X} _s}(T_{-\gamma}, B, \mathtt{R}))_{s\leq \lambda} )\right)d\gamma du     +o(1+x\ee^{-\frac{x^2}{t}}).
\end{eqnarray}
With an easy computation, we can obtain that
\begin{eqnarray}
\nonumber  \frac{4}{\sqrt{\pi}} \E_\alpha\left(\1_{\{\underline{B}_{t/2}\geq 0,\,  t\geq B_{t/2}\geq t^\frac{1}{4} \}}   B_{t/2}\ee^{-\frac{B_{t/2}^2}{t}}\right)&=&\frac{4}{\sqrt{\pi}} \alpha[\E_{\alpha}\left( \ee^{-\frac{ \mathtt{R}_{t/2}^2}{t} } \right) +o(1)]
\\
\label{part3} &=& \alpha \sqrt{\frac{2}{\pi}}(1+o(1)),\text{   uniformly in } \alpha\leq (\log t)^{30}.
\end{eqnarray}
Going back to (\ref{premierMarkov}), Proposition \ref{renouv} follows from a combination of (\ref{part1}) and (\ref{part3}).

\hfill$\Box$

\section{On the one dimensional Brownian motion $B$}
We refer to \cite{AShi10} and \cite{AShi10*} for the proof of the following Lemmas.
\begin{Lemma}
\label{Lem2.5}
 There exists a constant $c_{11}>0$ such that for any $x\geq 1$ and $t\geq 1$,
\begin{eqnarray}
\label{2.5} \P_{-x}\left(\overline{B}_t\leq 0\right)\leq c_{11}\frac{(1+x)}{\sqrt{t}},\qquad 
\P\left( B_t\in [x,x+1]\right) \leq \frac{c_{11}}{\sqrt{t}}.
\end{eqnarray}
\end{Lemma}

Recall the definitions of  $\BBi$ in (\ref{6.22}), $\Tra$ in (\ref{6.27}), $\DDi$ in (\ref{defclude}) and  $\Loza$ in (\ref{Dloza2}). From the Lemma \ref{Lem2.5} we deduce
\begin{Lemma}
\label{trajbrowseul}
There exists a constant $c_{12}>0$ such that

(i) for any  $b\geq a,z\geq 1$ and $ t\geq 1$ , 
\begin{equation}
\P_{-z}\left(B_t\in[-b,-a],\overline{B}_t\leq 0\right)\leq c_{12} z(1+b-a)(1+b)t^{-\frac{3}{2}},
\label{2.6}
\end{equation}

(ii) for any $a,\,  l,\, z,\,  L>1,\, \frac{t}{3}\geq a+l+1$ and $ m\in [t-a,t]$,
\begin{eqnarray}
\label{trajA.3} t^\frac{3}{2} \P\left(B\in \BBi^{z,L}_{t,a} \right) &\leq&  c_{12}  \E_z\left(B_{\log l}\1_{\{B_{\log l}\geq 0 \}} \right)(1+L) a^{-\frac{1}{2}},
\\
\label{trajA.4} {t^\frac{3}{2}}\P\left(B\in \Tra^{z,L}_{a,t}(m)\right)&\leq&  c_{12}  \E_z\left(B_{\log l}\1_{\{B_{\log l}\geq 0\}} \right),
\end{eqnarray}

(iii) for any $ l,\, z,\, L>1,\, t\geq l+1$ and $m\leq L$,
\begin{eqnarray}
\label{trajA.5bidbis} t^\frac{3}{2} \P \left(B\in \DDi_t^{z,L} \right)&\leq&  c_{12} z (1+L)^2  ,
\\
\label{trajA.6} t^\frac{3}{2} \P \left(B\in \Loza_t^{z,L}(m) \right)&\leq&  c_{12}\E_z\left(B^+_{\log l}\1_{\{\underline{B}_{\log l}^+\leq 1 \}} \right)  (1+L-m)  ,
\end{eqnarray}

(iv) for any $  t_b,\, \alpha>0,\,   \mathtt{z}\in [0,(\log t_b)^{30}]$ and $k\geq 0$,
\begin{eqnarray}
\label{trajA.7} t^\frac{3}{2} \P \left(  \overline{B}_{t_b}\leq \alpha,\, \overline{B}_{[\frac{t}{2},t_b]}\leq \alpha-\mathtt{z}  ,\, B_{t_b}+\mathtt{z}-\alpha\in [-(k+1),-k] \right)&\leq&  c_{12} z (1 +k).
\end{eqnarray}
\end{Lemma}
The proper proofs are minor adaptations of Lemma 2.2 \cite{AShi10*} for (\ref{2.6}); Lemma 2.4 in \cite{AShi10*} for (\ref{trajA.5bidbis}), (\ref{trajA.6}) and (\ref{trajA.7}); pp 14-15 in \cite{Aid11} for (\ref{trajA.3}) and (\ref{trajA.4}).

\begin{Remark}
\label{thoj}:  Each of these assertions can be proved by using the Markov property, Brownian time reversal and a combination of the inequalities in (\ref{2.5}).
\end{Remark}

In this section our aim is to extend Lemma \ref{Lem2.5} and Lemma \ref{trajbrowseul}.

\begin{Lemma}
\label{Browtronque}
There exists $c_{13}>0$ such that for any $t\geq 1,\, z\geq 0$, $u\in [0,t-1]$ and for any event $A(u)\in  \sigma\left( (B_{s+u}-B_u)_{s\in [0,1]}\right)$
\begin{eqnarray}
\label{eqBrowtronque1} 
 \P\left( \underline{B}_t\geq -z ,\,  A(u)  \right) &\leq& c_{13}\frac{1+z}{\sqrt{t}} \sqrt{\P\left( A(u)   \right)},
\\
\label{eqBrowtronque2}
 \P\left(B_t\in [z,z+1] ,\,    A(u) \right)& \leq&   \frac{c_{13}}{\sqrt{t}} \P\left( A(u)\right).
\end{eqnarray}
\end{Lemma}
\noindent{\it Proof of Lemma \ref{Browtronque}.}
First we prove (\ref{eqBrowtronque1}). If $u\geq \frac{t}{2}$,  by the Markov property at time $u$ and (\ref{2.5}) we get that
\begin{eqnarray*}
\P\left( \underline{B}_t\geq -z ,\,  A(u) \right)\leq \P\left( \underline{B}_u\geq -z \right)\P\left(A(u) \right)\leq c_{11} \frac{1+z}{\sqrt{t}}\P\left(A(u) \right)  .
\end{eqnarray*}
If $u\leq \frac{t}{2}$, by the Markov property at time $u+1$ and (\ref{2.5}) we get that
\begin{eqnarray*}
\P\left( \underline{B}_t\geq -z ,\,   A(u) \right) &\leq & \E\left(\1_{\{\underline{B}_{u+1}\geq -z\}}\1_{  A(u)}\P_{B_{u+1}}\left( \underline{B}_{t-(u+1)}\geq -z\right)\right)
\\
&\leq& \frac{c_{11}}{\sqrt{t}}\E\left((z+B_{u+1})\1_{\{ \underline{B}_{u+1}\geq -z,\,    A(u)   \}}\right).
\end{eqnarray*}
Observing that $B_{u+1}\leq B_u + \underset{s\leq 1} {\max} |B^{(u)}_s| $, we deduce that
\begin{eqnarray*}
\P\left( \underline{B}_t\geq -z ,\,   A(u) \right)&\leq& \frac{1}{\sqrt{t}}\E\left((z+B_u)\1_{\{ \underline{B}_{u}\geq -z \}}\1_{  A(u)  }\right) + \frac{1}{\sqrt{t}}\E\left(  \underset{s\leq 1} {\max} |B^{(u)}_s|; A(u) \right)
\\
&\leq & c\frac{1+z}{\sqrt{t}}  \sqrt{\P\left(A(u)\right)},
\end{eqnarray*}
where we have used the Cauchy-Schwarz inequality for the second term. So (\ref{eqBrowtronque1}) is proved. Now we prove (\ref{eqBrowtronque2}). If $u\leq \frac{t}{2}$, by the Markov property at time $u+1$ and (\ref{2.5}),
\begin{eqnarray*}
 \P\left(B_t\in [z,z+1] ,\,   A(u) \right)&\leq&  \E\left(\1_{ A(u)} \P_{B_{u+1}}\left( B_{t-u-1}\in [z,z+1] \right) \right)
 \\
 &\leq & \frac{c_{11}}{\sqrt{t}} \P\left( A(u)  \right).
\end{eqnarray*}
If $u\geq \frac{t}{2}$, then $t-u\leq \frac{t}{2}$ then we use the time reversal and Lemma \ref{Browtronque} follows. 
\hfill$\Box$

From the two previous results we can deduce:
\begin{corollary}
\label{jjkklllmmqq}
There exists a constant $c_{14}>0$ such that for any event $A(u)\in  \sigma\left( (B_{m+u}-B_u)_{m\in [0,1]}\right)$ and

(i) for any  $b\geq a,z\geq 1$ and $ t\geq 1$,
\begin{equation}
\P_{-z}\left(B_t\in[-b,-a],\overline{B}_t\leq 0 ,\,  A(u)\right)\leq  c_{14}z(1+b-a)(1+b)t^{-\frac{3}{2}}\sqrt{\P(A(u))},
\label{2.6666}
\end{equation}

(ii) for any $a,\, l,\, z,\,  L>1,\, \frac{t}{3}\geq l+1+a$, $\text{and } m\in [t-a,t]$
\begin{eqnarray}
t^\frac{3}{2} \P\left(B\in \BBi^{z,L}_{t,a},\,  A(u) \right) &\leq&  c_{14}  \E_z\left(B_{\log l}\1_{\{B_{\log l}\geq 0 \}} \right)(1+L) a^{-\frac{1}{2}}\sqrt{\P(A(u))},
\\
{t^\frac{3}{2}}\P\left(B\in \Tra^{z,L}_{t,a}(m) ,\,  A(u)\right)&\leq&  c_{14}  \E_z\left(B_{\log l}\1_{\{B_{\log l}\geq 0\}} \right)\sqrt{\P(A(u))},
\end{eqnarray}

(iii) for any $ l,\, z,\, L>1,\, t\geq l+1$ and $m\leq L$,
\begin{eqnarray}
\label{jkkjjjktrajB.13bidbis}t^\frac{3}{2} \P \left(B\in \DDi_t^{z,L},\,  A(u) \right)&\leq&  c_{14} z (1+L)^2 \sqrt{\P(A(u))} ,
\\
t^\frac{3}{2} \P \left(B\in \Loza_t^{z,L}(m),\,  A(u) \right)&\leq&  c_{14}\E_z\left(B^+_{\log l}\1_{\{\underline{B}_{\log l}^+\leq 1 \}} \right)  (1+L-m) \sqrt{\P(A(u))} ,
\end{eqnarray}

(iv) for any $  t_b,\, \alpha>0,\,   \mathtt{z}\in [0,(\log t_b)^{30}]$ and $k\geq 0$,
\begin{equation}
t^\frac{3}{2} \P \left(    \overline{B}_{t_b}\leq \alpha,\, \overline{B}_{[\frac{t}{2},t_b]}\leq \alpha-\mathtt{z}  ,\, B_{t_b}+\mathtt{z}-\alpha\in [-(k+1),-k] ,\, A(u) \right)\leq  c_{14} z (1+k) \sqrt{\P(A(u))}.
\end{equation}
\end{corollary}
\noindent{\it Proof of Corollary \ref{jjkklllmmqq}.} The result is an immediate consequence of the Remark \ref{thoj} and Lemma \ref{Browtronque}. Indeed we just have to reproduce the proofs for Lemma \ref{trajbrowseul} by replacing the inequalities in  (\ref{2.5}) by (\ref{eqBrowtronque1}) and (\ref{eqBrowtronque2}).
\hfill$\Box$

\section{On the fluctuations of the Gaussian processes $Z$ and $P$}
Recall that: - the process $(Z_s^0(y))_{s\in \r,y\in \r^d}$, is a centred Gaussian process with covariance:
\begin{eqnarray*}
\E\left( Z_u^0(y)  Z_v^0(z) \right):= \int_0^{u\wedge v} \left[ \mathtt{k}(\ee^s(y-z))-\mathtt{k}(\ee^sy)\mathtt{k}(\ee^{s}z)\right]ds,
\end{eqnarray*}

- the function $\mathtt{k}$ is symmetric and $\mathcal{C}^1$, in particular $\mathtt{k}'(0)=0$; 

- the function $\mathtt{g}:= 1-\mathtt{k}$ is $\mathcal{C}^1$, with $\mathtt{g}(0)=\mathtt{g}'(0)=0$,

- the function $\mathtt{g}'$ has a compact support included in $B(0,1)$, so there exists $c>0$ such that $\underset{y\in \r^\d}{\sup}|\mathtt{g}'(y)|\leq c |y|$.

\begin{Lemma}
\label{partieind}
There exist $c_{15},\,c_{16}>0$ such that

i) for any $\delta >0,\, T,\,t\geq 0$,
\begin{equation}
\label{eqpartieind}
 \P\left( \underset{|y|<\ee^{-T},\, s\in [0,t]}{\sup}  |Z_s^{0}(y)|  \geq \delta \right)\leq  c_{15}\exp({-c_{16}  \delta^2\ee^{2(T-t)}}),
\end{equation}

ii) for any $ b>0$ there exists $c_{17}(b)>0$ such that for any $\delta>0,\,t,\, j\geq 0$,
\begin{equation}
\label{eqpartieind2}
 \P\left( \underset{y,z\in B(0,\ee^b),\, |z-y|<j^{-1}}{\sup}  |Z_t^{0}(y\ee^{-t})-Z_t^{0}(z\ee^{-t})|  \geq \delta\right)\leq  \frac{c_{15}}{\delta^{2\d}} \exp({-c_{17}(b)(\delta j)^2}).
\end{equation}

\end{Lemma} 
\noindent{\it Proof of Lemma \ref{partieind}.}
Observe that 
\begin{eqnarray}
\label{decoupexp}
 \P\left( \underset{|y|<\ee^{-T}}{\sup}  \big|Z_\cdot^{0}(y)\big|_{t} \geq \delta \right)\leq \underset{k=0}{\overset{\lfloor t\rfloor }{\sum}}\P\left( \underset{|y|<\ee^{-T}}{\sup}  \big|Z_\cdot^{0}(y) \big|_{[k,k+1]} \geq \delta 2^{k-t}\right).
\end{eqnarray}
For $k\in [0,\lfloor t\rfloor]$, set $T_k:=([k,k+1])\times B(0,1) $. Recall Theorem 4.4.1 in \cite{Fer75} we introduce:
\begin{eqnarray*}
 \Gamma_k((u,x),(v,y))&:=& \E\left[Z_{u}^0({x\ee^{-T})}Z_{v}^0({y\ee^{-T}})\right], \text{ and } 
 \\
  \varphi_k(h)&:=& \underset{|u-v|\leq h}{\sup}\underset{ |x-y|\leq h}{\sup}\sqrt{\E\left( (Z_u^0({x\ee^{-T}})- Z_v^0({y\ee^{-T}})  )^2\right)}.
\end{eqnarray*} 
By an easy computation we have
\begin{eqnarray*}
\E\left( (Z_u^0({x\ee^{-T}})- Z_v^0({y\ee^{-T}}) )^2\right)  =\int_v^u[1-\mathtt{k}^{2}({\ee^{s-T} x})]ds+2\int_0^v\mathtt{g}({\ee^{s-T}}(x-y))ds-
\\
\int_0^v\left[\mathtt{k}({\ee^{s-T} y})-\mathtt{k}({\ee^{s-T} x})\right]^2ds.
\end{eqnarray*}
Therefore the Taylor expansion of $\mathtt{k}$ leads to   
\begin{eqnarray}
\label{previous} \varphi_k(h)&\leq & c {\sqrt{h}\ee^{k}}{\ee^{-T}},\quad \text{and  }\quad \underset{T_k^2}{\sup}\sqrt{\Gamma_k} \leq  {c}{\ee^{k-T}}.
\end{eqnarray}
Finally via the Theorem 4.4.1 of \cite{Fer75} and (\ref{previous}), we get that
\begin{eqnarray}
\nonumber&&\P\left( \underset{|y|<\ee^{-T}}{\sup}  \big|Z_.^0({y})\big|_{[k,k+1]} \geq \delta 2^{k-t}\right)
\\
\nonumber &\leq &  \P\left( \underset{|y|<\ee^{-T}}{\sup}  \big|Z_.^0({y})\big|_{[k,k+1]} \geq c \delta \ee^{T-t} (\frac{\ee}{2})^{t-k}\left[\underset{ T_k^2 }{\sup}\sqrt{\Gamma_k} +  \int_1^{\infty}\varphi_k(2^{-x^2})dx  \right]\right)
\\
\label{pioupiou} &\leq &  \frac{5}{2} 2^{2(d+1)} \int_{ c\delta \ee^{T-t} (\frac{\ee}{2})^{t-k} }^\infty   \ee^{-x^2}dx .
\end{eqnarray}
Going back to (\ref{decoupexp}) we obtain that $
 \P\Big( \underset{|y|<\ee^{-T}}{\sup}  \big|Z_\cdot^0({y})\big|_{t} \geq \delta \Big)\leq \underset{k=0}{\overset{\lfloor t\rfloor }{\sum}} c\ee^{-  (c\delta \ee^{T-t} (\frac{\ee}{2})^{t-k})^2}\leq c_{15}\ee^{-c_{16} \delta^2 \ee^{2(T-t)}}$, which proves inequality (\ref{eqpartieind}). Proof of inequality (\ref{eqpartieind2}) is similar, the details are omitted. Let us just mention that 
\begin{eqnarray*}
\varphi_t(h)&:=& \underset{ |x-y|\leq h}{\sup}\sqrt{\E\left( (Z_t^0({x\ee^{-t}})- Z_t^0({y\ee^{-t}})  )^2\right)}
\\
&=&\underset{ |x-y|\leq h}{\sup}\sqrt{ 2\int_0^t\mathtt{g}({\ee^{s-t}}(x-y))ds-\int_0^t\left[\mathtt{k}({\ee^{s-t} y})-\mathtt{k}({\ee^{s-t} x})\right]^2ds}
\\
&\leq & ch.
\end{eqnarray*} 
So instead of (\ref{pioupiou}), we can use here (see pp 54 in \cite{Fer75}, with $h=j^{-1},\, m=2j$, and $p=2$ )
\begin{eqnarray}
\nonumber &&  \P\left( \underset{y,z\in B(0,\ee^b),\, |z-y|<j^{-1}}{\sup}  |Z_t^{0}(y\ee^{-t})-Z_t^{0}(z\ee^{-t})|  \geq \delta \right)
 \\
\nonumber &\leq &  \P\left( \underset{x,y\in B(0,\ee^b),\, |x-y|<j^{-1}}{\sup}  |Z_t^0(y\ee^{-t})-Z_t^0(x\ee^{-t})| \geq c \delta j \varphi_t(j^{-1})\right)
\\
\label{pioupiou2} &\leq & c(b) \frac{5}{2} (2j) ^{2(d+1)}  \int_{c\delta  j }^\infty   \ee^{-x^2}dx .
\end{eqnarray}
\hfill$\Box$
\\

Now we shall estimate the fluctuations of the process $P_\cdot^0(y) $, i.e for $T,\,  t\geq\sigma >0,\, \delta>0,\, j\geq 1  $ and $b>0$,  we will control the events: \nomenclature[h3]{$A_{T,t,\delta} $}{$:=\{ \underset{\mid y \mid \leq \ee^{-T},\, s\in[0,t]}{\sup}\mid \int_0^s \mathtt{g}(\ee^v y)dB_v \mid \geq \delta \}  $}\nomenclature[h4]{$ B^{(\sigma)}_{j,t,\delta,b} $}{$:= \{ \underset{\mid y_1-y_2\mid \leq \frac{1}{j},\, (y_1,y_2)\in B(0,\ee^b)}{\sup}\mid \int_\sigma^t \mathtt{g}(\ee^{v-t} y_1)- \mathtt{g}(\ee^{v-t} y_2) dB_v \mid \geq \delta \} $}
\begin{eqnarray}
\label{Ev1} A_{T,t,\delta} &:=&\{ \underset{|y|\leq \ee^{-T},\, s\in[0,t]}{\sup}|\int_0^s \mathtt{g}(\ee^v y)dB_v |\geq \delta \},\qquad \text{and}
\\
\label{Ev2} B^{(\sigma)}_{j,t,\delta,b} &:=&\{ \underset{|y_1-y_2|\leq \frac{1}{j},\, (y_1,y_2)\in B(0,\ee^b)}{\sup}|\int_\sigma^t \mathtt{g}(\ee^{v-t} y_1)- \mathtt{g}(\ee^{v-t} y_2) dB_v |\geq \delta \}.
\end{eqnarray}
Event (\ref{Ev1}) appears in the proofs of Lemma \ref{LemL}, Proposition {\ref{tension*}}, inequality (\ref{eqcontroldiff}) 	and Lemma \ref{C.1}, whereas event (\ref{Ev2}) appears implicitly in (\ref{Gcond1}) and (\ref{trajB.21}).

We observe that for any $T,\, t>0,\, \delta>0$, $j\geq 1$, $b>0$ and $\sigma\in [0,t]$, 
\begin{eqnarray}
&& A_{T,t,\delta} \subset \underset{i=1}{\overset{\lfloor t+1 \rfloor}{\bigcup}} A_{T,t,\delta}(i),\qquad  B^{(\sigma)}_{j,t,\delta,b} \subset \underset{i=1}{\overset{\rfloor t+1\lfloor}{\bigcup}} B_{j,t,\delta,b}(i),
\end{eqnarray}
with for any $i\leq t+1 $, $A_{T,t,\delta}(i)$ and $ B_{j,t,\delta,b}(i)$ are measurable with respect to the sigma-field \\ $ \sigma\left( (B_{m+t-i}-B_{t-i})_{m\in [0,1]}\right)$ and defined by \nomenclature[h5]{$A_{T,t,\delta}(i) $}{$:=\Big\{ \underset{\mid y \mid\leq \ee^{-T}}{\sup} \underset{t-i\leq m\leq t-i+1}{\sup}\mid \int_{t-i}^m\mathtt{g}(\ee^{s}y)dB_s\mid \geq \delta 2^{-i} \Big\}  $} \nomenclature[h6]{$B_{j,t,\delta,b}(i) $}{$:= \Big\{ \underset{ \mid y-j\mid \leq j^{-1},\, \mid y_1\mid ,\mid y_2 \mid \leq \ee^b}{\sup} \underset{t-i\leq m\leq t-i+1}{\sup}\mid \int_{t-i}^m\mathtt{g}(\ee^{s}y_1)- \mathtt{g}(\ee^{s}y_2)dB_s \mid \geq \delta 2^{-i} \Big\} $}
\begin{eqnarray}
&&A_{T,t,\delta}(i):= \Big\{ \underset{|y|\leq \ee^{-T}}{\sup} \underset{t-i\leq m\leq t-i+1}{\sup}|\int_{t-i}^m\mathtt{g}(\ee^{s}y)dB_s|\geq \delta 2^{-i} \Big\},
\\
\label{mazedefinieition} &&B_{j,t,\delta,b}(i):= \Big\{ \underset{|y-j|\leq j^{-1},\, |y_1|,|y_2|\leq \ee^b}{\sup} \underset{t-i\leq m\leq t-i+1}{\sup}|\int_{t-i}^m\mathtt{g}(\ee^{s-t}y_1)- \mathtt{g}(\ee^{s-t}y_2)dB_s|\geq \delta 2^{-i} \Big\},
\end{eqnarray}

The following result is the core of this section
\begin{Lemma}
\label{majorAB}
There exist $c_{18},\, c_{19}>0$ such that for any $b>0$ there exists  $c_{19*}(b)>0$ such that for any 
$\delta>0,\, T,\, t\geq 0$ and $i\in \{1,...,\lfloor t+1\rfloor\}$
\begin{eqnarray}
\label{eqmajorAB}\P\left(A_{T,t,\delta}(i) \right)&\leq & c_{18} \exp\left(-c_{19} \delta^2(\frac{\ee}{2})^i\ee^{2(T-t)}\right),
\\
\label{eqmajorAB2}\P\left( B_{j,t,\delta,b}(i)\right)&\leq& \frac{c_{18}}{\delta^\d} \exp\left(-c_{19*}(b) \delta^2(\frac{\ee}{2})^i j\right).
\end{eqnarray}
Recall the definitions (\ref{Ev1}) and (\ref{Ev2}). By summing over $i\in \{1,...,\lfloor t+1\rfloor\} $ we deduce that there exists $c_{20}>0$ such that for any 
$\delta>0$, $ T\geq t\geq 0$ and $\sigma\in [0,t]$,
\begin{eqnarray}
\label{eqmajorAB3}\qquad  \P\left(A_{T,t,\delta}  \right) \leq   c_{20} \exp\left(-c_{19}  \delta^2 \ee^{2(T-t)}\right), \qquad \P\left( B^{(\sigma)}_{j,t,\delta} \right) \leq  \frac{c_{20}}{\delta^\d} \exp\left(-c_{19*}(b)  \delta^2  j\right).
\end{eqnarray}
\end{Lemma}
\noindent{\it Proof of Lemma \ref{majorAB}.} We start by (\ref{eqmajorAB}). By the Ito-formula for any $T,\, t$,  $s\geq s_1$, $|y|\leq \ee^{-T}$
\begin{eqnarray}
\nonumber|\int_{s_1}^s \mathtt{g}(\ee^{u-t}y)dB_u |= \left|  B^{(s_1)}_s\mathtt{g}(\ee^{s-t}y)- \int_{s_1}^s   <\nabla \mathtt{g}. y\ee^{u-t}>B^{(k)}_u du \right| \qquad\qquad
\\
 \leq c(s-s_1+1) \ee^{s-t}|y|\underset{u\leq s}{\sup}\, {|B^{(s_1)}_{u}|}.
\end{eqnarray}
Then 
\begin{eqnarray*}
\P\left(A_{T,t,\delta}(i) \right)&\leq& \P\left(3c \underset{|y|\leq \ee^{-T}}{\sup} |y| \underset{ m\in [0,1]}{\sup} \ee^{t-i-t} {|B^{(t-i)}_{m}|}\geq \delta 2^{-i}\right)
\\
&=&\P\left( \underset{ m\in [0,1]}{\sup}   {|B_{m}|}\geq \delta c'\ee^{T-t} \ee^i 2^{-i}\right) \leq c_{18} \exp\left(-c_{19} \delta^2(\frac{\ee}{i})^i\ee^{2(T-t)}\right).
\end{eqnarray*}

Concerning inequality (\ref{eqmajorAB2}) we use Lemme 4.1.3 pp 54 in \cite{Fer75} applied to the process:
\begin{eqnarray}
G_i(y,m):= \int_{t-i}^m\mathtt{g}(\ee^{s-t}y)dB_s,\qquad y\in B(0,\ee^b),\, m\in [0,1].
\end{eqnarray}
Indeed we first observe that for any $t,\, j,\, \delta>0$ and $i\in \{0,...,\lfloor t-1\rfloor\}$, 
\begin{eqnarray*}
B_{j,t,\delta,b}(i)\subset \left\{ \underset{|y_1-y_2|\leq j^{-1},\, |m_1-m_2|\leq j^{-1}}{\sup}|G_i(y_1,m_1)-G_i(y_2,m_2)|\geq \delta 2^{-i}     \right\}
\end{eqnarray*}
Then by an easy computation
\begin{eqnarray*}
\varphi(h)&:= &\underset{|y_1-y_2|\leq j^{-1},\, |m_1-m_2|\leq j^{-1}}{\sup}\sqrt{  \E\left( [G_i(y_1,m_1)-G_i(y_2,m_2)]^2\right)}
\\
&\leq &\sqrt{h}\exp({-i}),\qquad\qquad \forall h>0.
\end{eqnarray*}
So applying Lemme 4.1.3 pp 54 in \cite{Fer75} with $h=\frac{1}{j},\, m=2j\ee^b,\, p=2$, we get
\begin{eqnarray*}
\P\left( B_{j,t,\delta,b}(i) \right)&\leq& \P\Big( \underset{{\tiny \begin{array}{ll} |y_1-y_2|\leq j^{-1}
\\
 |m_1-m_2|\leq j^{-1}
 \end{array}}}{\sup}|G_i(y_1,m_1)-G_i(y_2,m_2)|\geq \frac{\delta 2^{-i}}{c\varphi(j^{-1})} [3\varphi(j^{-1})+ c' \int_1^\infty \varphi(j^{-1}2^{-u^2})]  \Big)
 \\
& \leq & c(b)j^{\d} \exp\left( -c'' \delta^2 j (\frac{\ee}{2})^i\right)\leq \frac{c_{18}}{\delta^{2\d}} \exp\left(-c_{19*}(b) \delta^2(\frac{\ee}{2})^i j\right).
\end{eqnarray*}
It ends the proof of Lemma \ref{majorAB}.\hfill$\Box$.
\\

Now we can state the following assertions ((\ref{c.21}), (\ref{trajB.16}),...,(\ref{trajB.21})) which are continuously used through the paper:
\\

Combining Corollary \ref{jjkklllmmqq} and (\ref{eqmajorAB}) we deduce that {\it There exists $c_{21}>0$ such that for any $T,\, t\geq 1,\, z\geq 0,\, \delta>0$
\begin{eqnarray}
t^\frac{3}{2} \P_z\left(B_t\in[a,b],\underset{s\leq t}{\inf}\,B_s\geq 0   ,\, A_{T,t,\delta}\right)
\nonumber&\leq&  t^\frac{3}{2} \sum_{i=1}^{\lfloor t+1\rfloor} \P_z\left(B_t\in[a,b],\underset{s\leq t}{\inf}\,B_s\geq 0   ,\,  A_{j,t,\delta}(i) \right)
\\
\nonumber&\leq& c z(1+b-a)(1+b)  \sum_{i=1}^{\lfloor t+1\rfloor}   \sqrt{c_{18} \exp\left(-c_{19} \delta^2(\frac{\ee}{2})^i\ee^{2(T-t)}\right)}
\\
\label{c.21} &\leq& c_{21} z(1+b-a)(1+b)  \exp({-\frac{c_{19}}{2}\delta^2\ee^{2(T-t)}}).
\end{eqnarray} }

Recall that $ A_{T,t,\delta} :=\{ \underset{|y|\leq \ee^{-T},\, s\in[0,t]}{\sup}|\int_0^s \mathtt{g}(\ee^v y)dB_v |\geq \delta \} $ (see (\ref{Ev1})). Similarly we can affirm that for some constant $c_{22}>0$ we have

(i) for any $a,\,  l,\, z,\,  L>1,\, t\geq l+1+a$ and $ m\in [t-a,t]$,
\begin{eqnarray}
\label{trajB.16} &&t^\frac{3}{2} \P\left(B\in \BBi^{z,L}_{t,a},\,  A_{T,t,\delta}  \right) \leq   c_{22}  \E_z\left(B_{\log l}\1_{\{B_{\log l}\geq 0 \}} \right)(1+L) a^{-\frac{1}{2}} \ee^{-\frac{c_{19}}{2}\delta^2\ee^{2(T-t)}}, 
\\
\label{trajB.17}  &&{t^\frac{3}{2}}\P\left(B\in \Tra^{z,L}_{t,a}(m) ,\,  A_{T-a+m,t-a+m,\delta}  \right)\leq  c_{22}  \E_z\left(B_{\log l}\1_{\{B_{\log l}\geq 0\}} \right)\ee^{-\frac{c_{19}}{2}\delta^2\ee^{2(T-t)}} ,
\end{eqnarray}

(ii) for any $ l,\, z,\, L>1,\, t\geq l+1$ and $m\leq L$
\begin{eqnarray}
\label{trajB.18bidbis} &&t^\frac{3}{2} \P \left(B\in \DDi_t^{z,L} ,\, A_{T,t,\delta} \right) \leq  c_{22} z (1+L)^2  \ee^{-\frac{c_{19}}{2}\delta^2\ee^{2(T-t)}} ,
\\
\label{trajB.19} & & t^\frac{3}{2} \P \left(B\in \Loza_t^{z,L}(m) ,\,  A_{T,t,\delta} \right)  \leq   c_{22}\E_z\left(B^+_{\log l}\1_{\{\underline{B}_{\log l}^+\leq 1 \}} \right)  (1+L-m)  \ee^{-\frac{c_{19}}{2}\delta^2\ee^{2(T-t)}} ,
\end{eqnarray}

(iii)  for any $  t_b,\, \alpha>0,\,   \mathtt{z}\in [0,(\log t_b)^{30}]$ and $k\geq 0$,
\begin{equation}
\label{trajB.20}  t^\frac{3}{2} \P \left(     \overline{B}_{t_b}\leq \alpha,\, \overline{B}_{[\frac{t}{2},t_b]}\leq \alpha-\mathtt{z}  ,\, B_{t_b}+\mathtt{z}-\alpha\in [-(k+1),-k] ,\, A_{T,t_b,\delta} \right)  \leq   c_{22} z (1+k)  \ee^{-\frac{c_{19}}{2}\delta^2\ee^{2(T-t_b)}}.
\end{equation}
\\

Finally let us prove the inequality used in (\ref{Gcond1}). We want bound  \\ $\P\left(  \overline{B}_{t_b}\leq \alpha,\, \overline{B}_{[\frac{t}{2},t_b]}\leq \alpha-\mathtt{z}  ,\, B_{t_b}+\mathtt{z}-\alpha\in [-(k+1),-k] ,\,  w^{(0,1)}_{\mathfrak{G}_{t,b,\sigma}(\cdot \ee^b)}(j^{-1})\geq \frac{1}{4}   \right) $. With Lemma \ref{rRhodesdec} observe that
\begin{eqnarray*}
&&\{w^{(0,1)}_{\mathfrak{G}_{t,b,\sigma}(\cdot \ee^b)}(j^{-1})\geq \frac{1}{4} \}\subset \{\sup_{x,y\in B(0,\ee^b),\, |x-y|<j^{-1}}\left| \int_{t_b-\sigma}^{t_b} \mathtt{g}(\ee^{s-t}y)-\mathtt{g}(\ee^{s-t}x)dBs \right| \geq 2^{-4}\}\cup
\\
&& \{\sup_{x,y\in B(0,\ee^b),\, |x-y|<j^{-1}} \left| Z^0_{t_b}(y\ee^{-t})-Z^0_{t_b}(x\ee^{-t}) \right|\geq 2^{-4} \}\cup \{\sup_{x,y\in B(0,\ee^b),\, |x-y|<j^{-1}} \left| \zeta_{t}(y\ee^{-t})-\zeta_{t}(x\ee^{-t}) \right|\geq 2^{-4} \}
\end{eqnarray*} 
Once $j$ large enough the last event is never realized. Therefore, recalling the definition (\ref{mazedefinieition}), we have
 \begin{eqnarray}
\nonumber &&  \P\left(  \overline{B}_{t_b}\leq \alpha,\, \overline{B}_{[\frac{t}{2},t_b]}\leq \alpha-\mathtt{z}  ,\, B_{t_b}+\mathtt{z}-\alpha\in [-(k+1),-k] ,\,  w^{(0,1)}_{\mathfrak{G}_{t,b,\sigma}(\cdot \ee^b)}(j^{-1})\geq \frac{1}{4}   \right)
\\
\nonumber && \leq  \sum_{i=1}^{\lfloor t+1\rfloor} \P\left(\overline{B}_{t_b}\leq \alpha,\, \overline{B}_{[\frac{t}{2},t_b]}\leq \alpha-\mathtt{z}  ,\, B_{t_b}+\mathtt{z}-\alpha\in [-(k+1),-k]   ,\,   B_{j,t,2^{-4}}(i)  \right)
\\
&&\nonumber     +\P_{-\alpha} \big( \overline{B}_{t_b}\leq 0,\, \overline{B}_{[\frac{t}{2},t_b]}\leq -\mathtt{z}  ,\, B_{t_b}+\mathtt{z}\in [-(k+1),-k]\big) \P\Big(  \underset{y,x\in B(0,\ee^b),\, |z-y|<j^{-1}}{\sup}  |Z_t^{0}(y\ee^{-t})-Z_t^{0}(x\ee^{-t})|  \geq 2^{-4}    \Big).
\end{eqnarray}
Using (\ref{eqpartieind2}) in Lemma \ref{partieind}, Corollary \ref{jjkklllmmqq} and (\ref{eqmajorAB2}) we deduce that there exists $c_{23}$ such that for any $b>0$ there exists $c_{24}(b)>0$ (a constant which depends on $b>0$) such that for any $t>0$ large enough, $\alpha \in [1,\log t]$, $k,\, j,\, \geq 1$, $\mathtt{z}\in [0,(\log t)^{30}]$ and $\sigma\in[0,t_b]$, 
\begin{eqnarray}
\nonumber &&\P\left(  \overline{B}_{t_b}\leq \alpha,\, \overline{B}_{[\frac{t}{2},t_b]}\leq \alpha-\mathtt{z}  ,\, B_{t_b}+\mathtt{z}-\alpha\in [-(k+1),-k] ,\,  w^{(0,1)}_{\mathfrak{G}_{t,b,\sigma}(\cdot \ee^b)}( j^{-1})\geq \frac{1}{4}   \right)
\\
\nonumber && \leq  c  \frac{\alpha (1+k)}{ t^{\frac{3}{2}}}\left[\Big( \sum_{i=1}^{\lfloor t+1\rfloor} c_{22}\sqrt{\P\left(  B_{j,t,2^{-4},b}(i)  \right)}  \Big)  +    \ee^{-c_{17}(b)(2^{-4} j)^2}\right]
\\
\label{trajB.21} && \leq c_{23} \alpha\frac{(1+k)}{t^{\frac{3}{2}}} \ee^{-c_{24}(b)j}.
\end{eqnarray}

\section{The $ L-$good particle}
Here we recall the definition of the ``good particles". It is convenient to introduce 
\begin{eqnarray*}
d_{i,l}^{L}(\rho(x)):=\left\{ \begin{array}{ll}   \log l^{\frac{2}{3}}\quad&\text{if  } i\in \{ 1,..., 5 \lfloor \log l\rfloor-1 \},
 \\
  \rho(x)-4e_i +D\quad &\text{if  } i\in \{ 5\lfloor \log l\rfloor ,...,\lfloor\frac{1}{2}t   \rfloor-1 \},
\\
a_t +\rho(x)+L-4e_i+D, \quad& \text{if  } i\in  \{\lfloor\frac{1}{2} t\rfloor,...,\lfloor t\rfloor\}.
\end{array}\right.
\end{eqnarray*}
where we recall that $e_s= s^{\frac{1}{12}}$ if $s\leq \frac{t}{2}$ and $e_s= (t-s)^{\frac{1}{12}}$ when $s\in [\frac{t}{2},t]$. Then, according to (\ref{Aket}), a particle $u\in [0,R]^\d $ is said to be $  L-\text{good}_i$ if   
\begin{eqnarray}
\label{defgood} &&\underset{y\in A_i(x)}{\sup}|Y_i(x)-Y_i(y)|\leq e_i+\frac{D}{2}     \text{   and  } Y_i(x)\leq d_{i,l}^L(\rho(x)),\qquad  i\in [1,\lfloor t\rfloor],
\end{eqnarray}
(see  (\ref{Aet}) for the definition of $ A_i(x)$).

\begin{Lemma}
\label{C.1}
Fix $L,R \geq  1$. For any $\epsilon>0$, there exists $D(L,\epsilon),l_0(L)$ large enough such that for any $l\geq l_0$ there exists $T(l,D)$ such that for any $ t\geq T$, $\rho(\cdot)\in \mathcal{C}_R(l, \kappa_\d \log l, +\infty)$, 
 \begin{equation}
\label{eqC.1}
 \int_{[0,R]^\d}\E\left(\frac{\1_{\{ Y_\cdot(x) \in \DDi_t^{\rho(x),L}\}}}{{\bf r_t}(x)^\d}\1_{\{ x\text{ not $L$-good} \}}\right)dx \leq  \epsilon\mathtt{I}_\d(\rho).
\end{equation}
\end{Lemma}

%
\noindent{\it Proof of Lemma \ref{C.1}.} Recall the definition of ${\bf r_t}(x)$ in (\ref{defrx}), for any $p\geq t$,  $\{ {\bf r_t}(x)< \frac{\ee^{-p}}{2} \}$ implies $ \{ w_{Y_{\cdot}(\cdot)}(\ee^{-p},x,t)\geq \frac{1}{4} \}$. Using Lemma \ref{rRhodesdec}, there exists ${\bf c}>0$ (as in the proof of Lemma \ref{maxfx}, see {\bf b)} in pp 16), ${\bf c}$ is a constant which depends only of $\mathtt{k}$, chosen in order to get rid of the deterministic part $\zeta_t^x$) such that for any $x\in [0,R]^\d$, $t\geq 0$, $p\geq t+{\bf c}$ and $ r\in (0,\ee^{-p}]$ , we have
\begin{eqnarray*}
 \{ w_{Y_{\cdot}(\cdot)}( r,x,t)\geq \frac{1}{4} \}\subset    \{ w_{P_{\cdot}^x(\cdot)}(r,x,t)  \geq 2^{-3} \} \cup  \{ w_{Z_{\cdot}^x(\cdot)}( r,x,t)  \geq 2^{-3} \}.
\end{eqnarray*}
So decomposing the value of ${\bf r_t}(x)$ in the intervals $[\ee^{-(t+{\bf c})},+\infty]$ and $[\ee^{-(p+1)},\ee^{-p}]$ with $p\geq t+{\bf c}$, for any $x\in [0,R]^\d$ one has
\begin{eqnarray}
\nonumber &&\E\left(\frac{\1_{\{ Y_\cdot(x) \in \DDi_t^{\rho(x),L}\}}}{{\bf r_t}(x)^\d}\1_{\{ Y_\cdot(x) \text{ not $L$-good} \}}\right)\leq c  \ee^{\d (t+{\bf c})} \P\left( Y_\cdot(x) \in \DDi_t^{\rho(x),L}   ,\,  x\text{ not $L$-good} \right) +
\\
\label{neeedto} && \qquad \qquad \qquad c\sum_{p\geq t+{\bf c}}\ee^{\d p}\E\Big(\1_{\{Y_\cdot(x) \in \DDi_t^{\rho(x),L},\,   x\text{ not $L$-good}\}}\left( \1_{\{ w_{Z_{\cdot}^x(\cdot)}( \ee^{-p},x,t)  \geq 2^{-3} \}} +\1_{\{ w_{P_{\cdot}^x(\cdot)}(  \ee^{-p},x,t)  \geq 2^{-3}\}} \right) \Big).
\end{eqnarray}
Then we need to :

A) decompose the event $ \{  x \text{ not $L$-good}\}$. Once $D$ large enough, for any $i\in [1,t]$, as $\mathtt{k}$ is Lipschitz,  $ \{ \underset{u\in A_i(x)}{\sup}|\zeta^x_i(u)|\geq  \frac{e_i+ \frac{D}{2}}{2}\}=\emptyset$, thus $ \{  x \text{ not $L$-good}\}$ is included in  the union from $i=1$ to $\lfloor t\rfloor $ of 
$$\{ Y_i(x) \geq d_{i,l}^{L}(\rho(x)) \} \cup \{ \underset{u\in A_i(x)}{\sup}|P^x_i(u)-Y_i(x)|\geq \frac{e_i+ \frac{D}{2}}{4} \}   \cup \{ \underset{u\in A_i(x)}{\sup}|Z^x_i(u)|\geq  \frac{e_i+ \frac{D}{2}}{4}\}, $$

B) by using the decomposition given by A), the events in (\ref{neeedto}) are either measurable according to $(Y_s(x))_{s\geq 0}$ either to $Z^x_\cdot(\cdot)$. Therefore, similarly to (\ref{numer1}) or also (\ref{h2goodk1}), we apply the Girsanov's transformation, with density $\ee^{\sqrt{2\d}Y_t(y)+\d t} $, to the two right hand terms of (\ref{neeedto}), recalling that $Y_\cdot(x)\in \DDi_t^{\rho(x),L}$ implies $\ee^{-\sqrt{2\d}Y_t(y) }\leq \ee^{-1} t^\frac{3}{2} \ee^{-\sqrt{2\d}\rho(x)}$.
\\

C) by using the the decomposition given by A), then the Girsanov's transformation of B), in the second term of the right and side of (\ref{neeedto}) appears naturally the following term:
\begin{eqnarray*}
&&c\sum_{p\geq t+{\bf c}}\ee^{\d p}\P\Big(  Y_\cdot(x) \in \DDi_t^{\rho(x),L},\,\underset{u\in A_i(x)}{\sup}|Z^x_i(u)|\geq  \frac{e_i+ \frac{D}{2}}{4}   ,\,  w_{Z_{\cdot}^x(\cdot)}( \ee^{-p},x,t)  \geq 2^{-3}  \Big)
\\
&&\leq ct^{\frac{3}{2}}\ee^{-\sqrt{2\d} \rho(x)}\P\big( B\in \DDi_t^{\rho(x),L}\big) \sum_{p\geq t+{\bf c}}\ee^{\d (p-t)}\P\Big(  \underset{u\in A_i(x)}{\sup}|Z^x_i(u)|\geq  \frac{e_i+ \frac{D}{2}}{4}   ,\,  w_{Z_{\cdot}^x(\cdot)}( \ee^{-p},x,t)  \geq 2^{-3}  \Big)
\end{eqnarray*}
To control the sum, we use the Cauchy-Schwarz inequality then Lemma \ref{partieind} and affirm that for any $x\in [0,R]^\d$, $i\in [0,t]$,
\begin{eqnarray*}
 &&\sum_{p\geq t+{\bf c}} \ee^{\d (p-t)}\P \Big(\underset{u\in A_i(x)}{\sup}|Z^x_i(u)|\geq  \frac{e_i+ \frac{D}{2}}{4},\,  w_{Z_{\cdot}^x(\cdot)}(  \ee^{-p},x,t)  \geq 2^{-3}  \Big)
 \\
& &\leq  \sum_{p\geq t+{\bf c}} \ee^{\d (p-t)} c_{15} \exp(-\frac{c_{16}}{2} 2^{-6}\ee^{2(p-t)})\sqrt{\P \Big(\underset{u\in A_i(x)}{\sup}|Z^x_i(u)|\geq  \frac{e_i+ \frac{D}{2}}{4}  \Big) }
\\
&& \leq  c\sqrt{\P \Big(\underset{u\in A_i(x)}{\sup}|Z^x_i(u)|\geq  \frac{e_i+ \frac{D}{2}}{4}  \Big) }.
\end{eqnarray*}

Finally, gathering A), B) and C), it stems that for any $x\in [0,R]^\d$, the first term plus the first part of the $\sum_{p\geq t+{\bf c}} ...$ in the right hand side of (\ref{neeedto}) is smaller than
\begin{eqnarray*}
  & c\ee^{-\sqrt{2\d} \rho(x)} t^{\frac{3}{2}}\left( \sum_{i=1}^t [(1)_A^{(i,x)}+ (2)_A^{(i,x)}+(3)_A^{(i,x)} ] \right),   & 
\end{eqnarray*}
with
\begin{eqnarray*}
&&(1)_A^{(i,x)}:=   \P\left(  B  \in \DDi_t^{\rho(x),L},\,   B_i \geq d_{i,l}^{L}(\rho(x))  \right),\quad 
\\
&&(2)_A^{(i,x)}:=   \P \Big(  B  \in \DDi_t^{\rho(x),L},\,     \underset{u\in A_i(0)}{\sup}|\int_0^i \mathtt{g}(\ee^v u)dB_v| \geq \frac{e_i+ \frac{D}{2}}{4}   \Big),
\\
&&(3)_A^{(i,x)}:= \P\Big(  B  \in \DDi_t^{\rho(x),L} \Big) \sqrt{\P \Big(\underset{u\in A_i(0)}{\sup}|Z^0_i(u)|\geq  \frac{e_i+ \frac{D}{2}}{4}  \Big)}.
\end{eqnarray*}
Similarly, the second part of the sum $\sum_{p\geq t+{\bf c}}...$ in the right hand side of (\ref{neeedto}) is smaller than
\begin{eqnarray*}
  & c\ee^{-\sqrt{2\d} \rho(x)} t^{\frac{3}{2}}\left( \sum_{i=1}^t [ (1)_B^{(i,x)}+(2)_B^{(i,x)} +(3)_B^{(i,x)} ] \right),    & 
\end{eqnarray*}
with
\begin{eqnarray*}
&&(1)_B^{(i,x)}:= \sum_{p\geq t+{\bf c}} \ee^{\d (p-t)} \P \Big(  B  \in \DDi_t^{\rho(x),L},\,   B_i \geq d_{i,l}^{L}(\rho(x)) ,\,    A_{p,t,2^{-3}}\Big),
\\
&&(2)_B^{(i,x)}:=   \sum_{p\geq t+{\bf c}} \ee^{\d (p-t)}  \P\Big(  B  \in \DDi_t^{\rho(x),L},\,    A_{p,t,2^{-3}} \underset{u\in A_i(0)}{\sup}|\int_0^i \mathtt{g}(\ee^v u)dB_v| \geq \frac{e_i+ \frac{D}{2}}{4}\Big)  ,
\\
&&(3)_B^{(i,x)}:= \sum_{p\geq t+{\bf c}}\P \Big(  B  \in \DDi_t^{\rho(x),L},\,    A_{p,t,2^{-3}}      \Big)\P \Big(\underset{u\in A_i(0)}{\sup}|Z^0_i(u)|\geq  \frac{e_i+ \frac{D}{2}}{4}  \Big) ,
\end{eqnarray*}
where we recall that $\{\underset{|u|\leq \ee^{-p} }{\sup}|\int_0^\cdot \mathtt{g}(\ee^v u)dB_v|_t  \geq 2^{-3} \} =A_{p,t,2^{-3}}$.
%
%
We start by studying $(1)_A^{(i,x)}$. We distinguish five cases:
\\

{\bf (i) $i\leq 5\log l$.} By definition $d_{i,l}^{L}(\rho(x))= \log l^{\frac{2}{3}}$, thus by the Markov property at time $i$ then (\ref{trajA.5bidbis}), we have
\begin{eqnarray*}
 (1)_A^{(i,x)}&\leq&  \E\left(\1_{\{ B_i\geq (\log l)^\frac{2}{3}\}}\P_{B_i-\rho(x)}\left( \overline{B}_{t-i}\leq  1,\, \overline{B}_{[\frac{t}{2}-i,t-i]}\leq a_t +L+1,\, B_{t-i}\geq  a_t -2 \right)\right)
\\
&\leq& c_{12}\E\left( (\rho(x)+1-B_i)_+ \1_{\{B_i\geq (\log l)^\frac{2}{3}\}}\right) \frac{(1+L)^2}{t^{\frac{2}{3}}}
\\
&\leq& c\ee^{-(\log l)^\frac{1}{6}}(1+L)^2\frac{\rho(x)+1}{t^{\frac{3}{2}}}.
\end{eqnarray*}
Then for all $l$ large enough, $\sum_{i=1}^{5\log l} (1)_A^{(i,x)} \leq c\log l\ee^{-(\log l)^{\frac{1}{6}}}(1+L)^2\frac{\rho(x)+1}{t^{\frac{3}{2}}}    \leq \epsilon \frac{ \rho(x)}{t^\frac{3}{2}} $.

{\bf (ii)  $5\log l\leq i\leq \frac{t}{3}$.} By definition $d_{i,l}^L(\rho(u))= \rho(u)-4i^{\frac{1}{12}}+D$, thus by the Markov property at time $i\geq 1$, then (\ref{trajA.5bidbis}) and (\ref{2.6}),
\begin{eqnarray*}
 (1)_A^{(i,x)} &\leq& c_{12}\frac{(1+L)^2 }{t^\frac{3}{2}}\E\left((1+\rho(x)-B_i)_+\1_{\{B_i\geq \rho(x) -4i^{\frac{1}{12}},\, \overline{B}_i\leq \rho(x)\}}\right)
\\
&\leq&  c  \frac{(1+L)^2}{t^\frac{3}{2}} (1+\rho(x)) \frac{i^\frac{1}{4}}{i^{\frac{3}{2}}}.
\end{eqnarray*}
Then for all $l$ large enough (depending on $L$), $\sum_{i=5\log l}^{t/3}   (1)_A^{(i,x)} \leq c \frac{(1+\rho(x)) }{t^\frac{3}{2}}  (1+L)^2\sum_{i=5\log l}^{t/3} {i^\frac{-5}{4}} \leq \epsilon \frac{ \rho(x)}{t^\frac{3}{2}}  $.

{\bf (iii)  $\frac{t}{3} \leq i\leq \frac{t}{2} $.} By definition $d_{i,l}^L(\rho(x))= \rho(x)-4i^{\frac{1}{12}}+D$, thus by the Markov property at time $i\geq \frac{t}{3}$ then by applying twice (\ref{2.6}), we get that 
\begin{eqnarray*}
 (1)_A^{(i,x)} &\leq& c\frac{(1+a_t)^2 }{t^\frac{3}{2}}\E\left((1+\rho(x)-B_i)_+\1_{\{B_i\geq \rho(x) -i^{\frac{1}{12}},\, \overline{B}_i\leq \rho(x)\}}\right)
\\
&\leq&  c  \frac{(1+a_t)^2}{t^\frac{3}{2}} (1+\rho(x)) \frac{i^\frac{1}{4}}{i^{\frac{3}{2}}}.
\end{eqnarray*}
Then for all $t$ large enough, $\sum_{i=t/3}^{t/2} (1)_A^{(i,x)}\leq c \frac{(1+\rho(x)) }{t^\frac{3}{2}}  (1+a_t)^2\sum_{i= t/3}^{t/2} {i^\frac{-5}{4}} \leq \epsilon \frac{ \rho(x)}{t^\frac{3}{2}}  $.

{\bf (iv)  $\frac{t}{2}\leq  i \leq   \frac{2t}{3}  $.} By definition $d_{i,l}^L(\rho(x))= a_t+\rho(x)+L-4e_i+D$, thus by the Markov property at time $i\geq \frac{t}{2}$ then with two times (\ref{2.6}), one has
\begin{eqnarray*}
 (1)_A^{(i,x)} &\leq& c_{12}\frac{(1+L)^2}{(t-i+1)^\frac{3}{2}}\E\left((1+\rho(x)+a_t+L-B_i)_+\1_{\{ B_i\geq a_t+\rho(x)+L-4e_i+D,\, \overline{B}_{[\frac{t}{2},i]}\leq a_t+L+\rho(x)\}}\right)
\\
&\leq&   c\frac{(1+L)^2}{(t-i+1)^\frac{3}{2}} e_i\P\left(   B_i\geq a_t+\rho(x)+L-4e_i+D,\, \overline{B}_i\leq \rho(x),\, \overline{B}_{[\frac{t}{2},i]}\leq a_t+L+\rho(x)   \right) 
\\
\\
&\leq & c\frac{(1+L)^2}{(t-i+1)^\frac{3}{2}}(4e_i-a_t)   \frac{e_i^2}{t^\frac{3}{2}}(1+\rho(x)) .
\end{eqnarray*}
As $e_i=(t-i)^\frac{1}{12}$, we have $
\sum_{i=\frac{ t}{2}}^{\frac{2t}{3}} (1)_A^{(i,x)} \leq \frac{1+\rho(x)}{t^\frac{3}{2}}  \sum_{i=\frac{t}{2}}^{\frac{2t}{3}}  \frac{c(1+L)^2}{(t-i+1)^\frac{5}{4}}  \leq \epsilon \frac{ \rho(x)}{t^\frac{3}{2}},$

{\bf (v)  $\frac{2t}{3}\leq  i \leq   t  $.} By definition $d_{i,l}^L(\rho(x))= a_t+\rho(x)+L-4e_i+D$, thus by the Markov property at time $i\geq \frac{2t}{3}$ then (\ref{2.6}) and (\ref{trajA.5bidbis}),
\begin{eqnarray*}
 (1)_A^{(i,x)} &\leq& c_{12}\frac{(1+L)^2}{(t-i+1)^\frac{3}{2}}\E\left((1+\rho(x)+a_t+L-B_i)_+\1_{\{ B_i\geq a_t+\rho(x)+L-4e_i+D,\, \overline{B}_{[\frac{t}{2},i]}\leq a_t+L+\rho(x)\}}\right)
\\
&\leq&   c\frac{(1+L)^2}{(t-i+1)^\frac{3}{2}} e_i\P\left(   B_i\geq a_t+\rho(x)+L-4e_i+D,\, \overline{B}_i\leq \rho(x),\, \overline{B}_{[\frac{t}{2},i]}\leq a_t+L+\rho(x)   \right)\1_{\{4e_i\geq D \}}
\\
\\
&\leq & c\frac{(1+L)^2}{(t-i+1)^\frac{3}{2}}   \frac{e_i^3}{t^\frac{3}{2}}(1+\rho(x))\1_{\{ 4e_i\geq D\}}.
\end{eqnarray*}

%
%
%
As $e_i=(t-i)^\frac{1}{12}$, we have $
\sum_{i=\frac{ 2t}{3}}^{t} (1)_A^{(i,x)} \leq \frac{1+\rho(x)}{t^\frac{3}{2}}  \sum_{i=\frac{2t}{3}}^{t-D^{12}}  \frac{c(1+L)^2}{(t-i+1)^\frac{5}{4}}  \leq \epsilon \frac{ \rho(x)}{t^\frac{3}{2}},$ once $D$ large enough ($D$ depends on  $L$).

Finally we conclude that for any $\epsilon>0$, there exist $l_0$ and $D$ large enough such that for any $l\geq l_0,\, t\geq \ee^l$, $x\in[0,R]^\d$, $\rho(\cdot)\in \mathcal{C}_R(l, \kappa_\d \log l, +\infty)$
\begin{equation}
\sum_{i=1}^{t} (1)_A^{(i,x)}  \leq \epsilon \frac{ \rho(x)}{t^\frac{3}{2}}.
\end{equation}
which ends the study of $ (1)_A^{(i,x)}$. The study of $(1)_B^{(i,x)}$ is quite similar. Indeed it consists to reproduce the case {\bf(i)} to {\bf(v)} by using (\ref{c.21}), (\ref{trajB.18bidbis}) instead of respectively (\ref{2.6}), (\ref{trajA.5bidbis}) it provides the following assertion: for any $\epsilon>0$, there exist $l_0$ and $D$ large enough such that for any $l\geq l_0,\, t\geq \ee^l$, $x\in[0,R]^\d$, $\rho(\cdot)\in \mathcal{C}_R(l, \kappa_\d \log l, +\infty)$
\begin{eqnarray}
\label{bn23}  \sum_{i=1}^{t}  (1)_B^{(i,x)}   &\leq& \epsilon \frac{ \rho(x)}{t^\frac{3}{2}} \sum_{p\geq t +{\bf c}}\ee^{\d t (p-t)} \exp({-\frac{c_{19}}{2} 2^{-6}\ee^{2(p-t)}}) \leq  \epsilon c\frac{ \rho(x)}{t^\frac{3}{2}} .
\end{eqnarray}
The details of the proof of (\ref{bn23}) are omitted.

Now we study $(2)_A^{(i,x)}$ and $ (2)_B^{(i,x)}$. First observe that
\begin{eqnarray*}
\{\underset{u\in A_i(0)}{\sup}  |\int_0^i \mathtt{g}(\ee^s u)dB_s| \geq \frac{e_i+ \frac{D}{2}}{4} \}&\subset & \overset{i}{\underset{j=1}{\bigcup}} \{ \underset{u\in A_i(0)}{\sup} |\int_{i-j}^{i-j+1}  \mathtt{g}(\ee^s u)dB_s|\geq \frac{ e_i+ \frac{D}{2}}{2^{j+2}}\}
\end{eqnarray*}
Recall that $\mathtt{g}\in \mathcal{C}^1$, thus by the Ito formula we can rewrite for any $j\in [1,i]$, $u\in A_i(0)$, $\int_{i-j}^{i-j+1}\mathtt{g}(\ee^s u)dB_s= \mathtt{g}(\ee^{i-j+1}u)(B_{i-j+1}-B_{i-j})-\int_{i-j}^{i-j+1}(B_s-B_{i-j})<\nabla_{\ee^{s}u}(g),\ee^s u>ds$. Recall also that $\mathtt{g}$ is Lipschitz with $\mathtt{g}(0)=0$ and $|\nabla_{\ee^{s}\cdot }(g)|$ is bounded. As $u\in A_i(0)$ implies $|u\ee^{i-j+1}|\leq c\ee^{3-j}$, for all $i\in [1,\lfloor t\rfloor]$ (for $i=1$ recall that $R$ is fixed), we deduce
\begin{eqnarray*}
\{\underset{u\in A_i(0)}{\sup}  |\int_0^i \mathtt{g}(\ee^s u)dB_s| \geq \frac{e_i+\frac{D}{2}}{4} \}&\subset & \overset{i}{\underset{j=1}{\bigcup}} \{  c \ee^{3-j} \underset{s\in [0,1]}{\sup} |B_{s+i-j}-B_{i-j}|   \geq \frac{ e_i+ \frac{D}{2}}{2^{j+2}}\}
\\
&\subset & \overset{i}{\underset{j=1}{\bigcup}} \{  \underset{s\in [0,1]}{\sup} |B_{s+i-j}-B_{i-j}| \geq \frac{1}{c'}\frac{e_i+ \frac{D}{2}}{4}(\frac{ \ee^1}{2})^{j}\},
\end{eqnarray*}
and for any $j\in \{1,...,i\}$,  $ \{   \underset{s\in [0,1]}{\sup} |B_{s+i-j}-B_{i-j}| \geq \frac{1}{c'}\frac{e_i+ \frac{D}{2}}{4}(\frac{ \ee^1}{2})^{j} \}$ is measurable with respect to $\sigma((B_{s+j}-B_j)_{s\in [0,1]} )$. Then according to (\ref{jkkjjjktrajB.13bidbis})
\begin{eqnarray}
 \nonumber (2)_A^{(i,x)}&\leq&  c_{14}  \frac{\rho(x)}{t^\frac{3}{2}} (1+L)^2 \overset{i}{\underset{j=1}{\sum}} \P\big(  \underset{s\in [0,1]}{\sup} |B_{s+i-j}-B_{i-j}| \geq \frac{1}{c'}\frac{e_i+ \frac{D}{2}}{4}(\frac{ \ee^1}{2})^{j}  \big)^\frac{1}{2}   
 \\
\label{bouit2} &\leq& c  \frac{\rho(x)}{t^\frac{3}{2}}(1+L)^2 \ee^{-c''({ \frac{D}{2}+e_i})}.
\end{eqnarray}
Similarly by the Cauchy-Schwarz inequality, then (\ref{jkkjjjktrajB.13bidbis}) and (\ref{trajB.18bidbis}) we get that
\begin{eqnarray}
\nonumber&&(2)_B^{(i,x)}\leq  \P\Big(   B  \in \DDi_t^{\rho(x),L},\,     \underset{u\in A_i(0)}{\sup}|\int_0^i \mathtt{g}(\ee^v u)dB_v| \geq \frac{e_i+ \frac{D}{2}}{4}\Big)^{\frac{1}{2}} \times\sum_{p\geq t +{\bf c}} \ee^{\d (p-t)}  \P\Big(  B  \in \DDi_t^{\rho(x),L},\,   A_{p,t,2^{-3}}\Big)^\frac{1}{2}
\\
\nonumber&&\qquad\qquad  \leq   \sqrt{(2)_A^{(i,x)} } \sqrt{(1+L)^2\frac{\rho(x)}{t^\frac{3}{2}} } \sum_{p\geq t +{\bf c}} \ee^{\d (p-t)}\ee^{-\frac{c_{19}}{4}2^{-6}\ee^{2(p-t)}}
\\
\label{bouit1}  &&\qquad \qquad\leq c  (1+L)^2\frac{\rho(x)}{t^\frac{3}{2}} \ee^{-c''({ \frac{D}{2}+e_i})}.
\end{eqnarray}
Combining (\ref{bouit1}) and (\ref{bouit2}), we get that for any $\epsilon>0$, there exist $l_0$ and $D$ large enough such for any $l\geq l_0,\, t\geq \ee^l$, $x\in[0,R]^\d$, $\rho(\cdot)\in \mathcal{C}_R(l, \kappa_\d \log l, +\infty)$
\begin{eqnarray}
\label{bn37} \sum_{i=1}^{t} [(2)_A^{(i,x)} + (2)_B^{(i,x)}  ] \leq   \frac{\rho(x)}{t^\frac{3}{2}}  c (1+L)^2  \ee^{-c''  \frac{D}{2}} \sum_{i=1}^t \ee^{-c'' e_i}\leq  \epsilon  \frac{ \rho(x)}{t^\frac{3}{2}}.
\end{eqnarray}
It remains to treat $(3)_{B}^{(i,x)}$. By (\ref{trajA.5bidbis}) and (\ref{trajB.18bidbis}) one has 

\begin{eqnarray*}
&& t^\frac{3}{2} \P \left(B\in \DDi_t^{\rho(x),L} \right) \leq  c_{12} \rho(x) (1+L)^2,
\\
&&t^\frac{3}{2} \P \Big(  B  \in \DDi_t^{\rho(x),L},\,    A_{p,t,2^{-3}}    \Big)  \leq  c_{22} \rho(x) (1+L)^2  \ee^{-\frac{c_{19}}{2}2^{-6}\ee^{2(p-t)}}.
\end{eqnarray*}
Moreover from Lemma \ref{partieind}, we see that
\begin{eqnarray*}
\P\Big(\underset{u\in A_i(0)}{\sup}|Z^0_i(u)|\geq  \frac{e_i+ D}{2}  \Big) &\leq& \P \Big(\underset{|u|\leq  c\ee^{2-i}}{\sup}|Z^0_i(u)|\geq  \frac{e_i+ D}{2}  \Big)
\\
&\leq&  c_{15}\exp({-c_{16} c ({e_i+ D})^{2}  \ee^{-4}}) .
\end{eqnarray*}
Combining these three inequalities we get that for any $\epsilon>0$, there exist $l_0$ and $D$ large enough such for any $l\geq l_0,\, t\geq \ee^l$, $x\in[0,R]^\d$, $\rho(\cdot)\in \mathcal{C}_R(l, \kappa_\d \log l, +\infty)$
\begin{eqnarray}
\label{bn98}   \sum_{i=1}^{\lfloor t \rfloor } [(3)_A^{(i,x)} + (3)_B^{(i,x)}  ] \leq    \frac{\rho(x)}{t^\frac{3}{2}}  c  (1+L)^2 \ee^{-c'  D} \sum_{i=1}^{\lfloor t \rfloor}\ee^{-c'' e_i} \leq  \epsilon  \frac{ \rho(x)}{t^\frac{3}{2}}.
\end{eqnarray}
Finally we deduce Lemma \ref{C.1} by gathering (\ref{bn23}) (\ref{bn37}) and (\ref{bn98}). \hfill$\Box$
\\

Observe that the event $\{ x\, \text{good}_i\, \forall i\in [2,\lfloor t\rfloor]\} $ does not depend of $R\geq 1$. Then as a by product of the previous proof we have the following corollary:
\begin{corollary}
\label{C.3}
For some constants $c,\, c'>0$, there exist $D(L,\epsilon),l_0$ large such that for any $l\geq l_0$, $\exists T(l,D)$ so that the following inequalities hold
\begin{equation}
\label{bibtC.3}
c'\ee^{-\d t} \rho(x)\ee^{-\sqrt{2\d}\rho(x)} \leq \P\left( Y_\cdot(x)\in  \DD_{t}^{\rho(x),0}, \,x\,   \text{good}_i\, \forall i\in [2,\lfloor t\rfloor]\right)\leq \P\left( Y_\cdot(x)\in  \DD_{t}^{\rho(x),0}\right) \leq c \ee^{-\d t} \rho(x)\ee^{-\sqrt{2\d}\rho(x)},
\end{equation}
provided that  $t  \geq T$, $R\geq 1$ and $\rho(x) \in [\kappa_\d\log l,\log t]$. 
\end{corollary}
\noindent{\it Proof of Corollary \ref{C.3}.} By applying the Girsanov's transformation with density $\ee^{\sqrt{2\d}Y_t(y)+\d t}$ 
\begin{eqnarray*}
t^{\frac{3}{2}}\ee^{-\sqrt{2\d}\rho(x)} \P\left( B\in  \DD_{t}^{\rho(x),0} \right) \leq \P\left( Y_\cdot(x)\in  \DD_{t}^{\rho(x),0} \right) \leq t^{\frac{3}{2}}\ee^{-\sqrt{2\d}\rho(x)} \P\left( B\in  \DD_{t}^{\rho(x),0} \right).
\end{eqnarray*}
From (\ref{trajA.5bidbis}) and (2.10) pp 6 in \cite{Aid11} we have also for any $t  \geq 1$ and $\rho(x) \in  [\kappa_\d\log l,\log t]$, 
\begin{eqnarray*}
\alpha_4 \frac{\rho(x)}{t^\frac{3}{2}}  \leq \P\left( B\in  \DD_{t}^{\rho(x),0} \right)\leq c_{12} \frac{\rho(x)}{t^\frac{3}{2}} .
\end{eqnarray*}
Finally it stems that
\begin{eqnarray}
\label{coro3}
 c'\ee^{-\d t} \rho(x)\ee^{-\sqrt{2\d}\rho(x)} \leq \P\left( Y_\cdot(x)\in  \DD_{t}^{\rho(x),0} \right) \leq c \ee^{-\d t} \rho(x)\ee^{-\sqrt{2\d}\rho(x)}.
\end{eqnarray}
It proves the upper bound. For the lower bound we just remark that
\begin{eqnarray*}
&&\P\left( Y_\cdot(x)\in  \DD_{t}^{\rho(x),0}, \,x\,   \text{good}_i\,  \text{for some}\,  i\in [2,\lfloor t\rfloor]\right) 
\\
&&\geq \P\left( Y_\cdot(x)\in  \DD_{t}^{\rho(x),0}\right) -\P\left( Y_\cdot(x)\in  \DD_{t}^{\rho(x),0} \,x\,   \text{not good}_i\,  \text{for some}\,  i\in [2,\lfloor t\rfloor]\right) .
\end{eqnarray*}
We choose $D$ large enough such that 
\begin{equation}
\label{coro2} \P\left( Y_\cdot(x)\in  \DD_{t}^{\rho(x),0},\,x\,   \text{not good}_i\,  \text{for some}\,  i\in [2,\lfloor t\rfloor]  \right) \leq  \epsilon \ee^{-\d t} \rho(x)\ee^{-\sqrt{2\d}\rho(x)}.
\end{equation}
We combine (\ref{coro2}) with (\ref{coro3}) to conclude.
\hfill$\Box$
\begin{figure}[h!]
\centering
\caption{}
\label{Tux32}
\includegraphics[interpolate=true,width=18cm,height=18cm]{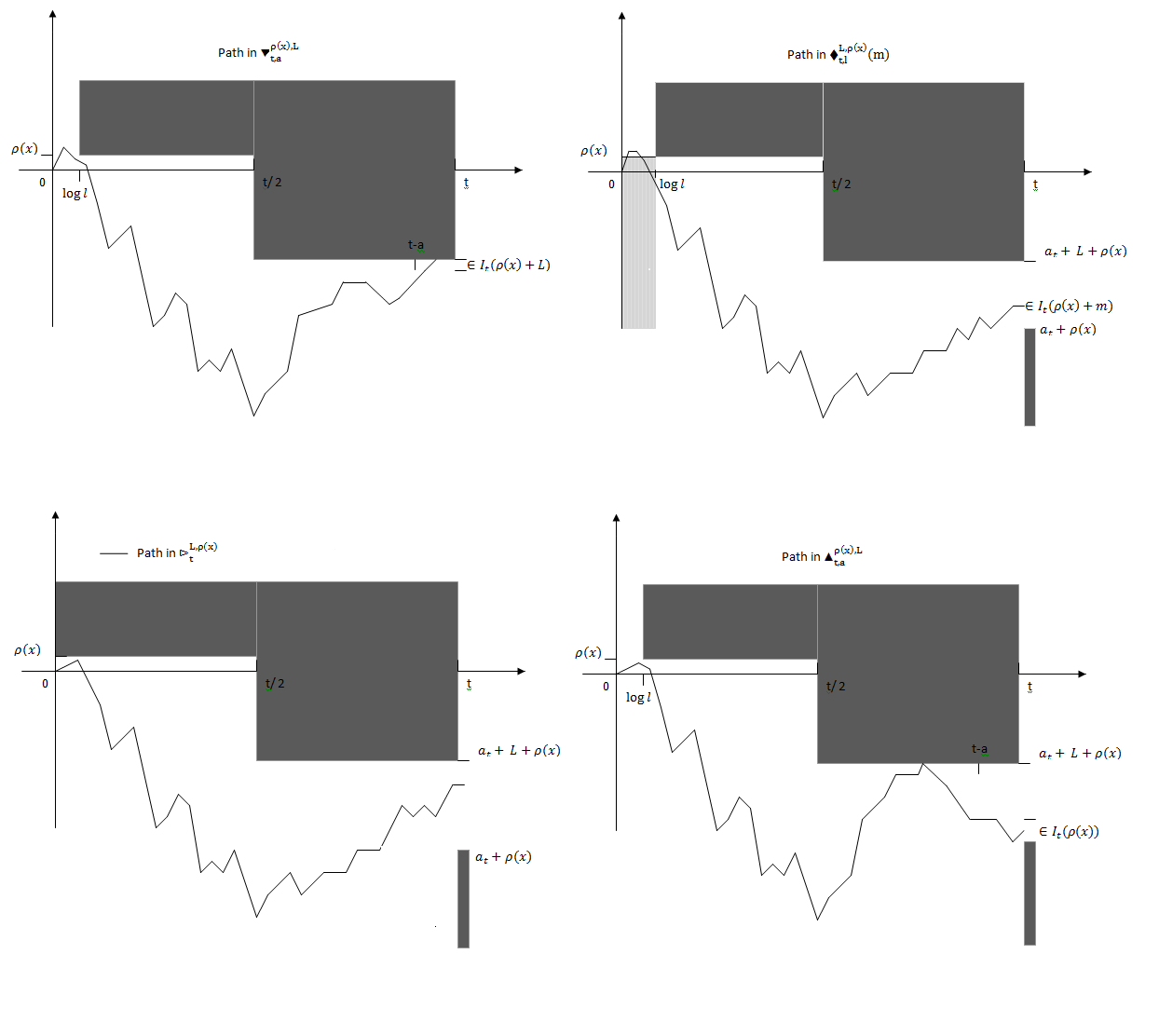} 
\end{figure}
\newpage

\paragraph*{Acknowledgement:} I would like to thank my advisor Yueyun Hu for his very useful advice. I would also like to thank the referee for his helpful comments.

\printnomenclature

\bibliographystyle{plain}
\bibliography{bibli}

\end{document}